\documentclass[mnsc,nonblindrev]{informs4a}
\SingleSpacedXI

\usepackage{endnotes}
\let\footnote=\endnote

\usepackage{tikz}
\usetikzlibrary{patterns}

\usepackage{float, subfigure}
\usepackage{wrapfig}

\usepackage[colorlinks=true, allcolors=blue]{hyperref}
\newif\ifproofread

\proofreadfalse

\newcommand{\changemarker}[1]{%
\ifproofread
\textcolor{black}{#1}%
\else
#1%
\fi
}

\usepackage{multirow}
\usepackage{booktabs}
\usepackage{array}

\usepackage[sort]{natbib}
 \bibpunct[, ]{(}{)}{,}{a}{}{,}%
 \def\bibfont{\small}%

\definecolor{myGreen}{RGB}{0,200,0}
\definecolor{darkGreen)}{rgb}{0.0, 0.5, 0.0}
\definecolor{myPink}{RGB}{255,153,153}
\definecolor{myCyan}{RGB}{0,206,209}

\def\hm#1{{\bf \color{magenta}#1}}

\def\kg#1{{ \bf \color{red}#1}}

\TheoremsNumberedThrough     % Preferred (Theorem 1, Lemma 1, Theorem 2)
\ECRepeatTheorems

\EquationsNumberedThrough    % Default: (1), (2), ...
\newcommand{\bA}{ \mathbf{A} }
\newcommand{\ba}{ \mathbf{a} }
\newcommand{\bb}{ \mathbf{b} }

\newcommand{\bc}{ \mathbf{c} }

\newcommand{\bG}{ \mathbf{G} }

\newcommand{\bh}{ \mathbf{h} }
\newcommand{\bI}{ \mathbf{I} }

\newcommand{\bP}{ \mathscrsfs{D} }

\newcommand{\bp}{ \mathbf{p} }
\newcommand{\bq}{ \mathbf{q} }

\newcommand{\bx}{ \mathbf{x} }

\newcommand{\by}{ \mathbf{y} }

\newcommand{\bzero}{ \mathbf{0} }

\newcommand{\hbx}{\hat{\mathbf{x}}}

\newcommand{\cA}{\mathcal{A}}
\newcommand{\cB}{\mathcal{B}}
\newcommand{\cC}{\mathcal{C}}
\newcommand{\cV}{\mathcal{V}}
\newcommand{\cN}{\mathcal{N}}
\newcommand{\cL}{\mathcal{L}}
\newcommand{\cNt}{\widetilde{\mathcal{N}}}
\newcommand{\cLt}{\widetilde{\mathcal{L}}}
\newcommand{\cI}{\mathcal{I}}

\newcommand{\cF}{\mathcal{F}}

\newcommand{\cS}{\mathcal{S}}

\newcommand{\cX}{\mathcal{X}}
\newcommand{\cXt}{\widetilde{\mathcal{X}}}
\newcommand{\cQ}{\mathcal{Q}}

\newcommand{\cK}{\mathcal{K}}
\newcommand{\cH}{\mathcal{H}}

\DeclareSymbolFont{rsfs}{U}{rsfs}{m}{n}
\DeclareSymbolFontAlphabet{\mathscrsfs}{rsfs}
\newcommand{\cKbad}{\mathcal{K}^{-}}
\newcommand{\cKgood}{\mathcal{K}^{+}}

\newcommand{\FO}{\mathbf{FO}}

\begin{document}
%%%%%%%%%%%%%%%%
\proofreadtrue  

%\setcounter{page}{1}
%\DOI{10.1287/xxxx.0000.0000}%

\RUNAUTHOR{Mahmoudzadeh and Ghobadi}

%\RUNTITLE{Learning from Good and Bad Decisions: A Data-driven IO Approach}
%\RUNTITLE{Learning Acceptability Criteria from Good and Bad Decisions: An IO Approach}
\RUNTITLE{\changemarker{Expert-guided IO for Constraint Inference}}

% \TITLE{Bundling Information Goods of Decreasing Value}
%\TITLE{Learning from Good and Bad Decisions: \\ A Data-driven Inverse Optimization Approach}
\TITLE{\changemarker{Expert-Guided Inverse Optimization for \\Convex Constraint Inference}}
%Expert-Guided Inverse Constraint Inference in \\Convex Optimization}
%\TITLE{\changemarker{Acceptability-aware Inverse Optimization for Convex Constraint Inference \\
%Feasibility-aware Inverse Optimization for Convex Constraint Inference\\
%Decision-aware Inverse Optimization for Convex Constraint Inference\\
%Expert Decision-aware Inverse Optimization for Convex Constraint Inference\\
%Expert-Guided Inverse Optimization for Convex Constraint Inference}}
%Old title: \TITLE{{Learning Acceptability Criteria from Good and Bad Decisions: An Inverse Optimization Approach}}
%\note{mytitlenote}
%Data-driven Inference of Feasible Regions}

\ARTICLEAUTHORS{%
\AUTHOR{Houra Mahmoudzadeh\footnote{Both authors contributed equally to this manuscript.}}
\AFF{Department of Management Science and Engineering, University of Waterloo, ON, Canada. \EMAIL{houra.mahmoudzadeh@uwaterloo.ca}}
%\AFF{Department of Bread Spread Engineering, Dairy University, Cowtown, IL 60208, \EMAIL{slippery@dairy.edu}} %, \URL{}}
\AUTHOR{Kimia Ghobadi${^{\text{\normalsize{*}}}}$}
\AFF{Malone Center for Engineering in Healthcare, Center for Systems Science and Engineering, Department of Civil\\ and Systems Engineering, Johns Hopkins University, Baltimore, MD, USA. \EMAIL{kimia@jhu.edu}}

} % end of the block

\ABSTRACT{%
Conventional inverse optimization inputs a solution and finds the parameters of an optimization model that render a given solution optimal. The literature mostly focuses on inferring the objective function in linear problems when accepted solutions are provided as input. In this paper, we propose an inverse optimization model that inputs several accepted and rejected solutions and recovers the underlying convex optimization model that can be used to generate such solutions. The novelty of our model is two-fold: 
First, we focus on inferring the parameters of \changemarker{the underlying convex feasible region}. Second, the proposed model learns the convex constraint set from a set of past observations that are either accepted or rejected by an expert. The resulting inverse model is a mixed-integer nonlinear problem that is complex to solve. To mitigate the inverse problem complexity, we employ variational inequalities and the theoretical properties of the solutions to derive a reduced formulation that retains the complexity of its forward counterpart. 
Using realistic breast cancer patient data, we demonstrate that our inverse model can utilize a subset of past accepted and rejected treatment plans to infer clinical criteria that can lead to nearly guaranteed acceptable treatment plans for future patients. 
}%

%\KEYWORDS{deterministic inventory theory; infinite linear programming duality;
%  existence of optimal policies; semi-Markov decision process; cyclic schedule}

%\KEYWORDS{inverse optimization, constraint inference, model learning, convex optimization, radiation therapy.}
%%%Keywords have to be from a list.
%\KEYWORDS{Optimization, Inference, Linear Programming, Nonlinear Programming, Healthcare Treatment.} 
% \HISTORY{Started to write June 30, 2020. } <----- literally 2 years!

\maketitle

\section{Introduction}
In the era of big data, learning from past \changemarker{expert} decisions and their corresponding outcomes, whether good or bad, %provides an invaluable opportunity for decision-makers facing similar situations. 
provides an invaluable opportunity for improving future decision-making processes. 
While there is considerable momentum %in the computer science community 
to learn from data through artificial intelligence, machine learning, and statistics, the field of operations research has not been using this valuable resource to its full potential in learning from past decisions to inform future decision-making processes. One of the emerging methodologies %in OR 
that can benefit from this abundance of data is inverse optimization~\citep{Ahuja01}. 

%Regular (forward) optimization models a system and determines an optimal solution that represents a 
A regular (forward) optimization problem models a system and determines an optimal solution that represents a decision for the system. On the contrary, inverse optimization aims to recover the optimization model that made a set of given observed solutions (or decisions) optimal. %These observations may consist of decisions made by experts when the underlying logic behind the decisions is not explicitly known. 
%given a set of solutions (decisions), conventional inverse optimization aims to recover the optimization model that made the observed solution optimal. %These observations may consist of decisions made by experts when the underlying logic behind the decisions is not explicitly known. 
For instance, in radiation therapy treatment planning for cancer patients, radiation oncologists make decisions on whether the quality of personalized plans generated through a treatment planning system is acceptable for each patient. In this case, an inverse model would be able to learn the implicit logic behind the oncologist's decision-making process. 
%In this case, an inverse model would learn the implicit logic that the oncologist employs in making such decisions. 
%In this case, the underlying optimization problem is the decision-making process of the expert (oncologist) to balance the complex tradeoffs between cancer treatment (objective) while adhering to a set of clinical guidelines (constraints).
%
Traditionally, the input to inverse models almost exclusively constitutes `good' solutions that are optimal or near-optimal, regardless of feasibility. Little attention has been paid to learning from `bad' solutions that must be avoided. In inverse optimization, learning from both {`good' and `bad'} observed solutions can provide invaluable information about the patterns, preferences, and restrictions of the underlying forward optimization model.  %\kg{we can connect to the next paragraph by talking about the meaning of good/bad in the obj first.}
%With the increased exploration of data-driven inverse models in recent years, larger sets of (sometimes real-world) observations can be considered in inverse settings. These extended sets of observations are increasingly more likely to include infeasible or unacceptable data points, whether due to errors (e.g., measurement errors) or based on expert domain knowledge (e.g., in dieting or cancer treatment).

%infeas

% nonlinear 
%In the literature, i

The inverse optimization literature largely focuses on inferring the objective coefficients of a forward model. The parameters often denote the utility or preferences of a decision-maker, when feasibility conditions are known. Very little attention has been paid to recovering the constraints. % of an optimization model. %Inferring constraints would determine the feasibility conditions in an optimization model.
There are three fundamental differences between recovering an objective function of an optimization problem when the feasible set is known and recovering the feasible region itself. 
%\begin{itemize} \setlength \itemsep{0pt}
    %\item 
    Firstly, when recovering the objective, an inverse model aims to satisfy optimality conditions for observed solutions (projections), regardless of their feasibility. On the contrary, when constraints are to be inferred, the focus shifts to satisfying feasibility conditions for the observed solutions while maintaining optimality conditions on a subset of observations that are optimal. This difference results in a mathematically-complex model that is harder to formulate and solve and has therefore largely been ignored in the literature, particularly in the presence of unfavorable solutions. %Specifically, we identified two key conceptual challenges when inferring the constraint set of a forward problem.  Firstly, 
    %\item 
    Secondly, when inferring constraints, any solution can be made optimal by inferring tailored constraints that render the observed solution optimal. Hence, traditional measures of (near) optimality that are used for objective coefficients, e.g., duality gap measures, are not applicable.
    %\item %contrary to objective inference, when inferring constraints, good and bad decisions can be viewed as solutions that are deemed acceptable and unacceptable with respect to a constraint set, respectively. %, as opposed to a set of %all-feasible  observations that are near or far from optimality, respectively, regardless of their feasibility status. 
     Thirdly, while good and bad decisions have various interpretations and implications when recovering objective functions, in the context of recovering constraints, they can be viewed as \changemarker{solutions that are deemed} accepted (feasible) and rejected (infeasible) \changemarker{by an expert}, respectively, \changemarker{which} guide the true shape of the underlying feasible region. %Hence, in inferring the objective, all observations can be made optimal or near-optimal regardless of their feasibility, whereas accepted and rejected decisions need to be treated differently when inferring constraints. 
Inverse optimization is well-studied for inferring linear optimization models \citep{chan2021IOreview}. This focus is mostly due to the tractability and existence of optimality guarantees in linear programming. % given that strong duality conditions of linear programming can be used.  % that are used in writing the inverse optimization model. 
In practice, however, %the underlying forward optimization is not always linear. N
nonlinear models are sometimes better suited for characterizing complex systems and capturing past solutions' attributes. 
%nonlinear models are sometimes required to better characterize  complex data structures and past observations. 
%constraint in general not inferred, only objective mostly. 
The literature also largely focuses on inferring the utility function of decision-makers, which can be interpreted as the objective function of an optimization problem when the feasible region is known. Inferring the feasible region itself, on the contrary, has not received much attention, which may be attributed to the fact that inverse models for constraint inference are nonlinear, even when the forward problem is linear. For linear problems, there have been recent attempts for recovering the forward feasible region through inverse optimization~\citep{ghobadi2014optimization, chan2018inverse}, however, these studies do not generalize to nonlinear problems \changemarker{and assume all observations are acceptable, hence do not incorporate expert-guided rejected observations} in constraint inference.

%+ then RT example
% RT
In radiation therapy treatment planning for cancer patients, a large pool of historical treatment plans exists that can be used in an inverse learning process. A plan is often designed to meet a set of pre-determined and often conflicting criteria, which are referred to as clinical guidelines. These guidelines %are blanket statements and not personalized, which means that they 
may be too strict or too relaxed (sometimes simultaneously), % for individual patients. %The oncology team typically resorts to designing a plan by identifying the guidelines that can be improved and the guidelines that have to be relaxed in an iterative process. 
and as a result, some plans that satisfy the original guidelines may be rejected by oncologists and some seemingly infeasible plans may be accepted. This application lends itself well to using inverse optimization for inferring the true underlying clinical guidelines for patient populations, which can lead to more efficient treatment planning and improved quality of treatment. 
In other words, the goal is to learn the parameters of an implicit optimization problem that can be solved to produce acceptable plans for future patients. % that replicate historically accepted plans. %are immediately accepted. % patients. 
% The radiation therapy treatment plan is often designed to meet the soft and hard constraints--the clinical objectives and guidelines for the patients. Clinical guidelines are often general and not personalized to patients, which means they could be too strict or too relaxed (sometimes simultaneously) for individual patients. The oncologist team typically resorts to designing the plan by identifying the guidelines that can be improved and the guidelines that have to be relaxed in an iterative process. As a result, some plans that are feasible may be rejected and some infeasible plans may be accepted by the team. Using inverse optimization to infer the true underlying clinical guidelines for patient populations will lead to more efficient treatment planning and possibly improved quality of the plans for patients. 
%constraint
While much attention has been paid to understanding the tradeoff balance in the objective of cancer treatment using inverse optimization, the problem of understanding the feasible region and constructing proper clinical guidelines remains under-explored. %~\citep{ajayi2022RT} 
An incorrect guideline or constraint in the optimization model can lead to a significantly different feasible region and affect the possible optimal solutions that the objective function can achieve.

%\subsection{Challenges with inferring the constraint set of a forward optimization problem}

%contribution:
In this paper, we focus on recovering \changemarker{convex} constraints of an optimization model through a novel inverse optimization framework. %for general convex problems. %We write our models when for a given objective function and a set of observations, whether optimal, feasible, or infeasible to the underlying optimization problem. 
%** AVOID SAYING "LABELLED DATA" to stay away from classification. **
Our model inputs a set of past observed decisions, that are either accepted or rejected by a rational expert decision-maker, and uses it to infer the underlying optimization problem that makes %these past decisions 
them feasible or infeasible, respectively.  
We further propose a reformulation of our inverse optimization model using variational inequalities to mitigate its complexity and improve solvability. We demonstrate the merit of our framework using the problem of radiation therapy treatment planning for breast cancer patients where we impute the underlying \changemarker{convex} constraints that represent the implicit guidelines that the expert decision-maker had in mind. % correctly characterize acceptable and non-acceptable treatment plans. 
The results can aid in standardizing clinical guidelines that can be used to produce acceptable plans, and hence, improving 
%the the quality of initial drafts of the treatment plan and 
the efficiency of the planning process.  %It what follows, we first review the 
Given a (forward) optimization problem with a set of partially known parameters, inverse optimization inputs a set of given solution(s) and recovers the \changemarker{problem} parameters \citep{Ahuja01}. The input solution is often a single observation that is optimal~\citep{ Iyengar05, ghate2020inverse} or \changemarker{near-optimal~\citep{Chan14, naghavi2019inverse} in which case the inverse model minimizes some measure of the optimality gap of the single input. Recently, with more focus on data-driven models, multiple noisy observations have also been considered as the input to inverse models~\citep{Keshavarz11,Troutt06,Troutt08,Chow12, Bertsimas15,  esfahani2018IncompleteInfo, Chan15, babier2021ensemble, zhang1999further, Aswani16, gupta2023efficient}, in which case, not all of the input observations can be (near-)optimal, and some measure of the collective data is often %assumed to be near-optimal
optimized instead. Some studies also consider uncertainty in data that affect the inverse models \citep{ghobadi2018robust}, or infer the structure of solutions to the inverse model instead of reporting a single inverse solution~\citep{tavasliouglu2018structure}.} For a comprehensive review of inverse optimization, we refer the readers to the review paper by~\cite{chan2021IOreview}.

Inverse optimization has been extensively considered for inferring linear optimization models. 
% This focus is mostly due to the tractability and existence of optimality guarantees in inverse linear models given that strong duality conditions of linear programming can be used.  % that are used in writing the inverse optimization model. 
% In practice, however, the underlying forward optimization is not always linear. Nonlinear models can often better characterize  complex data structures and past observations. 
%However, 
When the underlying forward optimization is assumed to be nonlinear, sufficient conditions for optimality of observations cannot be guaranteed unless the model falls under specific classes such as convex optimization for which Karush-Kuhn-Tucker (KKT) conditions are sufficient for optimality~\citep{boyd2004convex}. %For these models, the Karush-Kuhn-Tucker (KKT) conditions can be utilized to derive the required optimality conditions in inverse settings. 
\citet{zhang2010inverseSeparable} recover the objective function for linearly-constrained convex separable models. \citet{zhang2010InverseQuadratic} propose an inverse conic model that infers quadratic constraints and shows that it be efficiently solved using the dual of the obtained semi-definite programs. 
\citet{Keshavarz11} consider general convex models and use past observations to recover the objective function parameters by minimizing the optimality errors in KKT conditions. 
\citet{aswani2019data} also propose an NP-hard inverse optimization model to recover parameters for convex optimization models with noisy data and provide a duality-based reformulation. % with a regularization scheme that smooths discontinuities in the formulation. % and results in statistically consistent parameters.  
% \citet{Bertsimas15} propose an inverse variational inequality problem to model inverse equilibrium. %This last one seems irrelevant. 
% \kg{need to read some more literature to add here. \citep{Bertsimas15} solving conic forward??}
%\kg{keep?} While these studies have advanced the theory of inverse optimization for inferring nonlinear forward models, the focus has been the inference of the utility function of decision-makers, which translates to inferring the objective function of the forward models, and constraint inference has not received much attention in the literature. 

% \hm{Switch with the next paragraph. End on const inference being at a donkey stage.}
%%%%%%%%%%%%% constraint inference
%\hm{make sure to cite Archis's paper that infers some constraint parameters - I think transition probabilities for MDP. \cite{ghate2022MDPinverse}}
\changemarker{The current inverse optimization literature mainly focuses on inferring the objective function of the underlying forward problem \citep{chan2021IOreview}. Constraint inference, on the contrary, has received little attention. Recovering the right-hand side of the constraint parameters alongside the objective parameters has been explored by \citet{dempe2006inverse, Chow12} and \citet{cerny2016inverse}. Similarly, \citet{birge2017inverse} recover the right-hand side parameters so that a given observation becomes optimal utilizing properties of the specific application and \citet{dempe2006inverse, guler2010capacity, saez2018short} make a given observation optimal or near-optimal. \citet{chan2018inverse} perturb the nearest facet to make a given observation optimal and hence, find the left-hand side parameter of a linear constraint when the right-hand side parameters are known. \citet{ghate2022MDPinverse} infer the unknown transitional probabilities in Markov Decision Processes which are part of the left-hand-side of the constraint set. \cite{Ghate2024SDP_IO_Const} infer left-hand side constraint parameters of a semi-definite program. \cite{aswani2019data} and \cite{gupta2023efficient} write KKT conditions to infer full convex problems, but only test their models to infer the objective parameters, due to the complexity of constraint inference. Closest to our work is the study by \cite{ghobadi2021inferring} in which the full set of the constraint parameters is inferred in a linear model where the objective function and a set of feasible observations are given. Their method utilizes properties of linear optimization and does not generalize to convex constraints.} 
\changemarker{When imputing the objective function parameters in forward models with known feasible regions, infeasible observations (as near-optimal decision) have been considered by \citet{babier2021ensemble, chan2022pathways, ahmadi2020inverseLearning}, and \citet{shahmoradi2019quantile}, among others. }
%However, since these studies focus on recovering the objective function, the infeasible solution provide similar information as the feasible solution and are considered 'good' observations that contain information about the direction of the objective function.  
%
%%However, because the feasible region is pre-determined in these studies, any infeasible observation is still treated as a `good' decision that needs to be made near-optimal for the inferred objective. %
%\kg{needs more explanation of their work.} 
%In these works, the authors utilize the infeasible observations in a similar manner as the feasible observations to extract information about the utility function and the objective tradeoffs of the forward optimization. 
In these studies, infeasible and feasible observations are used in a similar manner to extract information about the utility function and provide objective parameter trade-offs in the forward problem~\citep{shahmoradi2019quantile, ahmadi2020inverseLearning}. Outliers and irrelevant observations may also be removed from the data prior to the optimization or as part of it. % \kg{others? chan? }. 
%
% When imputing constraint parameters, on the contrary, infeasible observations can be treated as `bad' decisions. Hence, they can provide an additional layer of information by %identifying the locations that are not part of the feasible region, hence, 
% enabling a better determination of where the feasible region should reside, which areas should be excluded from it, and which constraint parameters would provide a better fit. %Chan \& Ghobadi’s methods did not allow for the inclusion of infeasible 
%
When imputing constraint parameters, on the contrary, infeasible observations provide additional information to guide the shape of the underlying feasible region by signaling areas that must be made infeasible. To our knowledge, such undesired observations that must be avoided have not been explored in an inverse setting for constraint inference. 

% Paper structure: 
% \newpage

% \hm{maybe move this paragraph to the end of 1.1 and lead to the motivation in the next section.}
Inverse optimization has found a wide range of applications including energy \citep{ brucker2009inverse, birge2017inverse, fernandez2021forecasting}, dietary recommendations \citep{ghobadi2018robust, shahmoradi2019quantile, ahmadi2020inverseLearning, ahmadi2022MLIOClusteringDiet, SHAHMORADI20222ClustringIO}, finance \citep{li2021FinanceIO,roland2016finding, yu2023learning}, and healthcare systems \citep{chan2022pathways}, to name a few. %and cancer treatment planning \citep{babier2018inverse, ajayi2022RT}. %\kg{need more? there are more applications if we need them.}
%%%%%%%%%%%%% cancer care application
In particular, radiation therapy treatment planning for cancer has been studied in the context of inverse optimization~\citep{Chan14,chan2018trade,ajayi2022RT,chan2022pathways, boutilier2015IMRT, babier2020RT, babier2018inverse, Goli15, lee2013predicting,ghate2020RTFractionation}. %need to add a lot more here 
%Explain work on inverse for cancer planning. Taewoo, Rafid, etc. %
% The current literature mostly focuses on recovering the objective function and better understanding the underlying complex tradeoffs between different objective terms in radiation therapy treatment plans. 
%Most are linear, with the exception of ... . They all input good solutions only. 
For instance, both \citet{chan2018trade} and \citet{sayre2014automaticRT} input accepted treatment plans to recover the appropriate weights for a given set of convex objectives using inverse optimization. \citet{gebken2019inverse} finds the objective weights for unconstrained problems using singular value decomposition. Personalization for different patient groups has been explored by \citet{boutilier2015IMRT} by recovering the utility functions appropriate to each group. 
\cite{ajayi2022RT} employs inverse optimization for feature selection to identify a sparse set of clinical objectives for prostate cancer patients.  
%explain that they focus on objective inference mostly. Taewoo's objective selection for cancer treatment with Andrew.  
% need to add even more .. better explanation of the literature. 
%The current literature mostly 
These studies all focus on understanding the underlying tradeoffs between different objective terms in radiation therapy treatment plans and only
%These studies exclusively
use accepted treatment plans as an input to their inverse models, regardless of their feasibility with respect to the clinical guidelines. 
% When the accepted plans are infeasible with respect to guidelines, %they are treated as 'good' plans and
% the inverse models aim to infer the forward objective that makes them near optimal. 
In Section~\ref{sec:cancer}, we present some of the challenges of finding acceptable treatment plans based on current clinical guidelines, which motivate the methodologies proposed in this paper. 

% \kg{Stopped here -- Sept 14}

\subsection{Cancer Treatment Motivation} \label{sec:cancer}

In 2023, there were an estimated 1.96 million new cancer cases diagnosed and 609,820 cancer deaths in the United States, and approximately 60\% of them received radiation therapy as part of their treatment \citep{cancerstats}. The radiation therapy treatment planning process is a time-consuming process that %has to be repeated for each patient and 
often involves manual planning by a treatment planner and/or oncologist.
%RT define, clinical metrics and physician's perspective in setting these goals
%maybe figure here. 
%motivate the problem
%Consider the {\bf radiation therapy treatment planning} problem for {\bf cancer patients}. 
The input of the planning process is a medical image (e.g., CT, MRI) which includes contours that delineate the cancerous region (i.e., tumor) and the surrounding healthy organs at risk (OAR). The goal is to find the direction, shape, and intensity of radiation beams such that a set of clinical metrics on the tumor and the surrounding healthy organs is satisfied. In current practice, there are clinical guidelines on these radiation metrics. % for %radiation delivered %dose %thresholds 
%\kg{dose gradients, conformity, homogeneity, dose-volume criteria for specific organs or regions of the body, etc.} 
% the tumor and the surrounding healthy organs. 
% The guidelines are often blanket statements, too loose for some patients who can have better guidelines or too tight for some that their approved plans (good plans) are infeasible for the guidelines. There are patients who meet the current clinical guidelines but the plans are rejected (bad plans) because the oncologists believe that better plans can be achieved. For the other patients, new guidelines should be used to align with their needs. 
However, these guidelines are not universally agreed upon and often differ per institution. Additionally, adherence to these guidelines is at the discretion of oncologists. % based on their expert opinion. %the specific needs of each patient.  % as they see fit for each patient.  
% Often times the guidelines are blanket statements that are either too restrictive and cannot be met for all patients, and plans that do not meet g

% \kg{Just listing the points we want to make: clinical guidelines are too broad (leading to unnecessary back and forth + resulting in suboptimality for some patients + manual relaxation for patients who cannot meet the criteria- hard patients), not standard and vary by institution and sometimes individual oncologists, not adhered too because too broad/tight, criteria not always "optimization friendly" to be included sometimes resulting in linear or convex interpretations of it manually (instead of inferring the best linear/convex approximation of it). Anything else?}

Planners often try to find a treatment plan that meets these clinical guidelines and forward it to an oncologist who will, in turn, either accept or reject the plan. If the plan is rejected, the planner receives a set of instructions on which metrics to adjust. 
This iterative process can lead to unnecessary back and forth between the planner and the oncologist and may involve manual relaxation of the required clinical criteria. %, and can result in suboptimal plans for patients. 
This process continues until the plan is accepted by the oncologist. 

% The final approved plan may or may not meet all the clinical guidelines simultaneously as there are trade-offs between different metrics and the clinical guidelines are not personalized for each patient.  %Either remove this last sentence or revise the next paragraph because it is repetitive... % I'd say rewrite this last sentence to better connect with the next paragraph
 %remove the underlined part? and merge with the next paragraph? 
%Suppose we have a set of approved treatment plans from previous patients. Even though there are clinical guidelines on acceptable thresholds for different metrics, in reality, 

\begin{figure}[htbp]
    \centering
    \subfigure[Feasible region of the clinical guidelines \label{fig:motivation_a}]{
    % \includegraphics[width=0.47\textwidth, trim={1cm 1.2cm 1cm 0.5cm},clip]{images/motivation1-small.eps}}
    % \subfigure[Implicit acceptance criteria\label{fig:motivation_b}]{
    % \includegraphics[width=0.47\textwidth , trim={1cm 1.2cm 1cm 0.5cm},clip]{images/motivation2-small.eps}
    \includegraphics[width=0.47\textwidth, trim={1.5cm 1.2cm 1.5cm 0.5cm},clip]{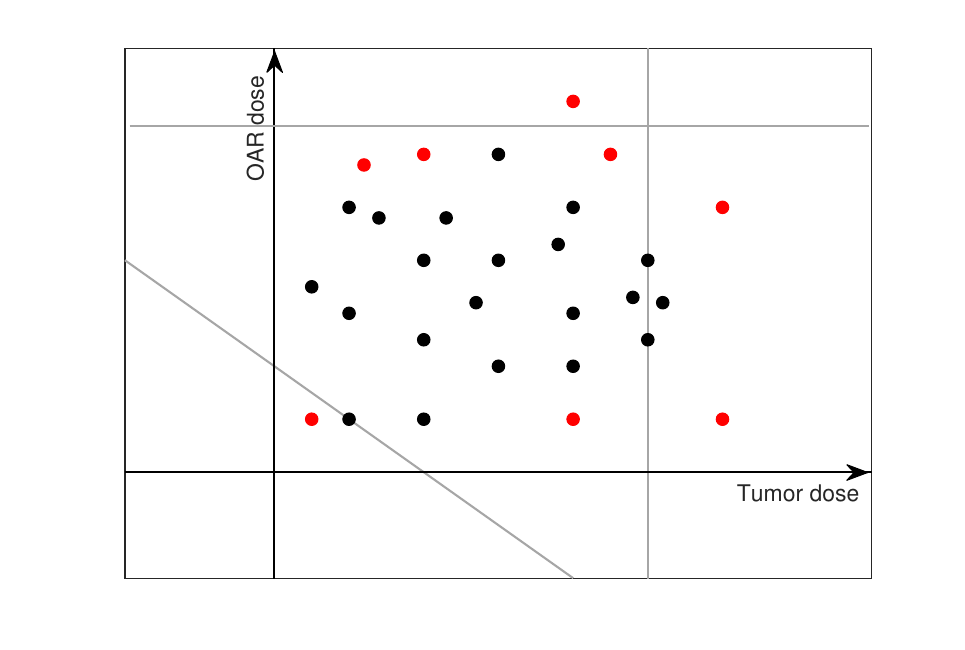}}
    \subfigure[Implicit feasible region of oncologists \label{fig:motivation_b}]{
    \includegraphics[width=0.47\textwidth , trim={1.5cm 1.2cm 1.5cm 0.5cm},clip]{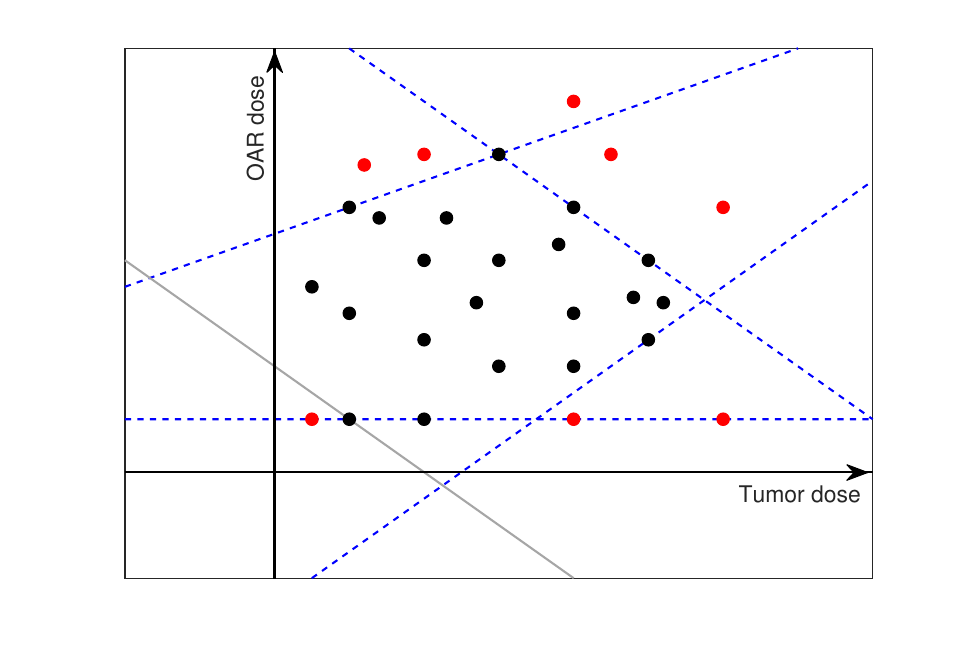}}
        \caption{Simplified schematic representation of convex feasible regions based on guidelines (gray solid) versus implicit (blue dashed) constraints of oncologists. The black and red dots denote accepted and rejected plans, respectively.
        % \kg{Now I'm thinking that maybe real data would be better .. lol! Maybe the desmos figure or if we can get a figure for the false positive and false negatives in the 105 patients (or the 5 patients)}. \hm{I know, I keep changin my mind too, but we can't have a 2D representation of the real data and the inferred/known constraints like this here. Also if we bring the data here, we'd have to explain what dvh criteria are etc.} \kg{It can be 2d if we reduce it to acceptale (red) , rejected (black), and then 2 of the criteria (and we mention that it's only 2 of the criteria. Or maybe showing the FN/FP rates??}
        }
    \label{fig:guidelines}
\end{figure}

As we will show later in Section~\ref{sec:case}, most accepted treatment plans do not meet all the clinical guidelines simultaneously, typically because there are trade-offs between different metrics, %the guidelines are not personalized for each patient, 
and some radiation dose limits are too restrictive for some patients. This will lead to an increased back-and-forth between the planner and the oncologist. %Conversely, there may also exist plans that meet the guidelines but are not approved because the oncologist may find the guidelines too relaxed for some patients and believe better plans are achievable, which may also lead to an increased back and forth between the planner and the oncologist. 
Conversely, sometimes the practically accepted plans follow much tighter constraints than the guidelines because the oncologists may believe that some guidelines are too relaxed and better plans are achievable.  
In mathematical programming terminology, the implicit (unknown) feasible region, based on which oncologists make an acceptance/rejection decision, often does not align with the feasibility/infeasibility of plans for the guidelines. 
%is unknown. 
A hypothetical schematic of accepted and rejected plans with respect to two metrics (Tumor dose and OAR dose) is shown in Figure~\ref{fig:motivation_a}. It can be seen that some points do not meet the guidelines but are accepted and others meet all guidelines but are rejected. There may be other complex criteria that capture the trade-off between the OAR dose and the tumor that oncologists consider when deciding on the acceptability of a treatment plan, as shown in~\ref{fig:motivation_b}. Understanding these trade-offs results in more practical and standardized guidelines that allow the planners to accurately represent the feasible region of the treatment optimization problem. Solving such an optimization problem that considers the learned guidelines can help planners produce acceptable treatment plans for future patients. 
%
% that learns from past data to generate acceptable plans for future patients.  
% \hm{We are assuming that the oncologist has an implicit set of criteria in mind and is solving an optimization problem to find which plans are optimal for each patient. The goal of IO is to find the parameters of this underlying optimization problem such that a planner can solve that optimization problem and find plans that are immediately approved by the oncologist. (Look at Mohajerani's 2023 routing paper) }

%\hm{give stats on what percentage meet or don't meet...} 
%The implicit constraints, however, always correctly classify the accepted/rejected plans. %\hm{explain OAR}

%Motivation Figure~\ref{fig:motivationFO}. 

% \hm{Cite optimization papers that focus on RT. \cite{bertsimas2020novel} }

In the radiation therapy treatment planning problem, we demonstrate that %by considering the historically-approved plans as ``feasible points'', we can employ our inverse optimization framework to find the constraints parameters. 
our inverse framework can learn from both accepted and rejected plans and infer %the constraint parameters that shape the underlying feasible region. \kg{am I right that we basically want to map the accepted/rejected to feas/infeasible? Since accep/rej is what eventually matters to oncologists. So reducing FT and FN between them hence, improving the process and possibly the quality of plans??} \hm{Yup!} 
the true underlying criteria based on which the accept/reject decisions are made. 
Finding such constraints enables us to better understand the implicit logic of oncologists in approving or rejecting treatment plans. In doing so, we help both oncologists and planners by ({\it i}\,)~standardizing the guidelines and care practices, ({\it ii}\,)~generating more realistic criteria on the trade-offs between clinical metrics based on past observations as opposed to simple upper/lower bounds on individual metrics, ({\it iii}\,) improving the quality of the initial plans according to the oncologist's opinion %given the clear guidelines and hence, 
and hence, reducing the number of iterations between planners and oncologists, and ({\it iv}\,) improving the quality of the final plans by preventing low-quality solutions that otherwise satisfy the acceptability thresholds, especially for automated treatment planning methods that heavily rely on provided radiation thresholds and may result in infeasibility if clinical thresholds are not personalized. 

%\hm{2023-10-04: We stopped here. NEXT TIME: start editing from here onward! }

\subsection{Contributions} %\kg{This part needs review at the end: }
%This paper aims to recover the underlying feasible region of a general convex forward problem based on both good (feasibility acceptable) and bad (feasibility unacceptable) observations. 
This paper aims to recover the \changemarker{underlying convex feasible region} of a forward problem based on both accepted and rejected observations of past decisions. 
The goal is to provide a streamlined process to replace the current ill-fitting optimization models that are used as guidelines in practice to find the optimal solutions to a forward problem, and hence, the resulting solutions undergo iterative revisions with \changemarker{practical expert-driven} guidelines.  
% We note that the goal of IO fundamentally differs from conventional classification models, due to the nature of the underlying optimization decision-making problem. 
\changemarker{Using historical expert opinions, the} goal of the proposed IO methodology is to learn a structured optimization model that can later be used in a forward setting to generate decisions that mimic successful decisions that have been accepted by experts in the past. The optimal solutions of such a learned model will facilitate and streamline the approval process by avoiding historically unacceptable solutions. %output of our inverse model is used in future decision-making through solving a convex optimization problem based on the inferred guidelines.
%
%To our knowledge, there has been no work that uses both feasible and infeasible observations in inferring the constraints of a linear or general convex forward problem.
% \hm{Contrast with classification - We are not ONLY trying to classify observations. We are planning to infer the underlying optimization model so that it can be used to generate new plans for new patients. we need to describe the overall big picture better. Our goal is to update current ill-fitting optimization models so that they would fit the observed data more accurately. These optimization models (the guideline standards) will be used prospectively.}
% \hm{again, contrast with classification and maybe cite a review paper. say why those cannot be used for our purpose. We need to infer a model for future decision making, not just classification. Explain role of FO: We are inferring an FO (with an obj function, constrary to classification methods) that will render the observations feas or infeas.  }
The specific contributions of this paper are as follows. 
%over
%\hm{Bullet point for contributions here \\}
%\hm{
%This paper makes the following specific contributions
\begin{enumerate} \setlength \itemsep{0pt}
    \item We propose an inverse model to recover a fully- or partially-unknown \changemarker{convex} feasible region of a forward optimization problem. %first nonlinear convex opt constraint inference
    \item The proposed inverse model inputs both accepted and rejected decisions that guide the shape of the imputed feasible region for the forward problem.  %first feasible and infeasible learning in const inference
    \item We propose a reformulation of the proposed model using variational inequalities to reduce its computational complexity. %, the LP-relaxation of which retains the complexity of the forward model.  %reduced reformulation based on variational inequalities
    \item We demonstrate an application of the proposed methodology in standardizing the radiotherapy clinical guidelines for cancer treatment. %application in inferring clinical guidelines.
\end{enumerate}
%}

In the remainder of this paper, we first define our forward optimization problem mathematically and present the proposed inverse optimization model in Section~\ref{sec:method}. Next, in Section~\ref{sec:solution}, we present a reformulation of the proposed inverse model to mitigate its computational complexity. Finally, we apply our methods to an example of deriving clinical guidelines for radiation therapy treatment planning for breast cancer patients in Section~\ref{sec:case}, and conclude the paper in Section~\ref{sec:conclusions}.
    
% Make a new figure that shows the guidelines and two patients. One that meets all the guidelines but is labeled as infeasible. Another that does not meet some of the guidelines but is labeled as feasible. The guidelines are a subset of items shown in figure 1.Then maybe move Figure 1 back to its original place.  

% Add a table on current guidelines for breast cancer RT. Point out a specific case on the barplot that meets the guidelines but was rejected. 

% % Motivation Figure~\ref{fig:motivationFO}. 
% \begin{figure}
%     \centering
%     \subfigure[\label{fig1}]{
%     \includegraphics[width=0.47\textwidth]{images/motivation1-small.eps}}
%     \subfigure[\label{fig2}]{
%     \includegraphics[width=0.47\textwidth]{images/motivation2-small.eps}}
%     \caption{Think: can we remake this figure to include nonlinear? }
%     \label{fig:motivationFO}
% \end{figure}

%-----------------------
\section{Methodology}  \label{sec:method}
%-----------------------
In this section, we first formulate \changemarker{a forward optimization problem with convex constraints}, where all or some of the constraints are unknown. We then define the inverse problem mathematically where a set of accepted and rejected observations are given, and the goal is to find constraint parameters that correctly classify these observations while enforcing optimality conditions on a preferred solution. We then present our \changemarker{expert-guided inverse optimization model and characterize the properties of its solutions.}

\subsection{Problem Definition} \label{sec:problemdef}
Let $\cI$ be the set of all constraints in a forward optimization problem. We denote the set of known nonlinear and linear constraints by $\cN$ and $\cL$, respectively, and the set of unknown nonlinear and linear constraints to be inferred by and $\cNt$ and $\cLt$, respectively. Note that $\cN \cup \cL \cup \cNt \cup \cLt = \cI$ and $\cN, \cL,\cNt, \cLt$ are mutually exclusive sets. We note that the known constraints are a trusted subset of the guidelines that need to be satisfied by all future solutions, which can potentially be an empty set. Assume that the decision variable is $\bx \in \mathbb{R}^m$. \changemarker{Consider a rational decision maker with a differentiable monotone convex objective function $f(\bx; \bc)$ and let $g_n(\bx; \,\bq_n), \, \forall \cN$ be differentiable concave functions on $\bx$}. The convex forward optimization (FO) model can be formulated as: 
\begin{subequations}
\begin{align}  
 \FO: \qquad \underset{\bx}{\text{minimize}} \quad &f(\bx; \bc) \label{eq:FO-obj}\\ 
 \text{subject to} \quad & g_n(\bx; \,\bq_n) \ge \bzero\,, && \forall n \in \cN \cup \cNt \label{eq:FO-nonlin}\\ 
 & \ba_{\ell}'\, \bx \geq b_{\ell}\,.   && \forall \ell \in \cL \cup \cLt \label{eq:FO-lin}
\end{align} \label{eq:FO}
\end{subequations}
%\kg{Kimia to think about flipping the $g$ sign. if g is convex then $g>0$ would be a non-convex set. we need to either flip the direction of inequality or declare that g is concave, meaning that -g is convex. -- I have to come back to this later}
%
% \hm{probably move definitions to the previous section.} 
\noindent Note that because $g_n(\bx, \bq_n)$ is concave, the \changemarker{nonlinear constraints~\eqref{eq:FO-nonlin}
 and linear constraints~\eqref{eq:FO-lin} form a convex feasible region. }%$g_n(\bx; \bq_n) \geq \bzero$ corresponds to a convex set. 
For brevity of notations, let the set of all known constraints be defined as $\cX=\{\bx \in \mathbb{R}^m \mid g_n(\bx; \,\bq_n) \ge \bzero, \, \forall n \in \cN, \,\, \ba_{\ell}' \, \bx \geq b_{\ell}, \, \forall \ell \in \cL\}$, the region identified by the known constraints of FO. \changemarker{Given that the objective function ~\eqref{eq:FO-obj} minimizes a convex monotone function $f(\bx;\bc)$, the optimal solution of FO is always on the boundary of its feasible region.}

% \hm{add a definition for the set of obs, x: }
%We define the set defined by such constraints as follows.  
%\kg{$a_\ell$ is the vector of the $\ell$th row of $A$. \\}
Assume that the structures %is there a better word? type? properties? class?
of the functions $g_n(\bx; \bq_n), \, \forall n \in \cN$ are known, and the goal is to find unknown parameters  $\bq_n, \, \forall n \in \cNt$. Note $\ba_{\ell} \in \mathbb{R}^m$ and $b_{\ell} \in \mathbb{R}$ are parameters of fixed size while each $\bq_n$ might be a vector of a different length $\bq_n \in \mathbb{R}^{\phi_n}$, where $\phi_n$ depends on the type of nonlinear function that is to be inferred. % (the number of terms we are inferring the coefficients of).  
For example, $g_1(\bx; \bq_1) = q_{11} x_1^2 + q_{12} x_2^2 + q_{13} x_1 x_2 + q_{14}$ is a two-dimensional quadratic function with four unknown parameters $q_{11},\dots,q_{14}$ to be inferred.  In what follows, we describe the proposed inverse methodology for imputing the  constraint parameters of the FO model. 
%---H: the example is still fine when g is concave. 
% \hm{Somewhere early on we need to say something like this: The purpose of IO is to find the parameters q,a,b such that when solving the FO problem, we are able to produce solutions that the expert (oncologist) would approve as an acceptable plan. }

%\kg{add motivation figure here later. A simpler figure than Figure 1. Impose the known constraints (as opposed to guidelines) and inferred (imputed) on the same figure. }

\subsection{Inverse Problem Formulation}\label{sec:IOFormulation}
Let $\bx^k,\,  k \in \cK$ denote a set of given solutions corresponding to past decisions, where $\cK=\cKgood \cup \cKbad$ and $\cKgood$ and $\cKbad$ denote accepted and rejected observed decisions, respectively.  The goal of the inverse problem is to find a set of parameters $\bq_n$, $\ba_{\ell}$, and $b_{\ell}$ such that solving the corresponding FO model will result in an optimal solution that mimics the accepted observation and avoids the rejected ones. 
%In this section, we propose an inverse formulation that inputs these past decisions and infers the set of linear and nonlinear constraints of the forward problem such that all previously accepted constraints are inside the inferred forward feasible region and the rejected observations become infeasible. The inverse model also identifies a preferred solution based on the objective function of the forward problem and infers the forward feasible region such that this preferred point would be optimal for FO. 
To this end, the inverse optimization problem infers a convex feasible region that would render all the past observations $\bx^k,\, k\in \cK^+$ as feasible for FO, and contrarily, all $\bx^k,\, k\in\cK^-$ as infeasible for FO. 

% \hm{Maybe an additional intro sentence. } *********************
%Before introducing the inverse model, we first outline a few standard assumptions on the structure of the problem and the observed data to ensure that the problem is well-defined and feasible. ********************
Let $\cH$ be the convex hull of all accepted observations $\bx^k, \, k \in \cKgood$. %denotes decisions that have been accepted (good) and $\cKbad$ are those that were rejected (bad). 
We assume that the decision maker is rational and is solving an implicit convex optimization problem, which results in the observed data being well-posed and the inverse problem being feasible, as shown in Assumption~\ref{assumption:wellposed}.  
%The assumption is that the decision-maker was solving an implicit convex optimization problem to make the past decisions. 
%
\begin{assumption} \label{assumption:wellposed}
    The sets of given observations are well-posed and the forward problem is convex, %and the known constraints of the forward problem are  well-posed, % to each other and the underlying convex forward optimization,
    i.e.,  
    \begin{itemize}
        \item[\normalfont{(a)} ] $\exists \, \hat{\bq}_n \in \mathbb{R}^{\phi_n} $, such that $g_n(\bx_k;\hat{\bq}_n) \geq \bzero, \quad \forall n \in \cNt, \,\,  k \in \cKgood$, %The sets of known constraints and good observations are well-defined with respect to each other. 
        \item[\normalfont{(b)}  ] $\bx^k \in \cX, \quad \forall k \in \cKgood$, %All accepted observations satisfy the known constraints. 
        \item[\normalfont{(c)}  ] $\not \exists k \in \cKbad$ such that $\bx^k \in \cH$. %, \, k\in \cKbad$. %, where $\cH$ is the convex hall of $\bx^k$. %There exists an underlying convex feasible region that includes all accepted observations and none of the rejected ones. 
    \end{itemize} 
\end{assumption}

\begin{figure}[htbp]
    \centering 
\begin{tikzpicture}[thick,scale=1.2, every node/.style={scale=1.1}]
% % %----axes
\draw[->, >=stealth] (0,0) -- (5,0);
\draw[->, >=stealth] (0.5,-.5) -- (0.5,3);
% \draw[lightgray] (0, .5) -- (5, .5); %known constraint
% % %----c vector:
% \draw[->]  (4.5,2.5) coordinate -- (4.2,2.2) node[anchor= north west] {$\bc$};
% \draw[->, blue]  (4.25,2.35) coordinate -- (3.95,2.05) ;
\draw[->, blue]  (2.2,-.5) coordinate -- (1.8,-.65) ;
% \draw (4.6, 2.4) -- (4.4, 2.6);
% \draw[ dashed , blue] plot [smooth] coordinates{(3,2.75) (4.25,2.35) (4.5, 1.25)} node[anchor= south west] {$f(\bx;\bc)$} ; 
\draw[ dashed , blue] plot [smooth] coordinates{(0.7,1.05) (1.95,0.65) (2.2, -.95)} node[anchor= south west] {$f(\bx;\bc)$} ; 
\draw[ dashed , blue] plot [smooth] coordinates{(0,0.75) (1.25,.35) (1.5, -1.25)}; 
%node[anchor= south west] {$f(\bx;\bc)$} ; 
% \draw[ dashed , blue] plot [smooth] coordinates{(-.7,0.35) (.5,-.05) (.8, -1.15)}; %node[anchor= south west] {$f(\bx;\bc)$} ; 
%
%----------observations
    \draw [fill] (3,1.8) circle [radius=0.05];
    % \node[anchor=east] at (3,2) (coord) {$\bx^1$};
    \draw [fill, blue] (1.8,.8) circle [radius=0.05];
    % \node[anchor=west] at (1.8,.8) (coord) {$\bx^2$};
     \node[anchor=south west, blue] at (1.8,.8) (coord) {$\bx^0$};
    \draw [fill] (2.7,2.3) circle [radius=0.05];
    \draw [fill] (4,2) circle [radius=0.05];
    % \node[anchor=south west] at (4,1.8) (coord) {$\bx^3$};
    \draw[fill] (2.8,1.2) circle [radius=0.05];
    % \node[anchor=south west] at (2.8,1.2) (coord) {$\bx^4$};
    \draw[fill] (3.7,0.8) circle [radius=0.05];
    % \node[anchor=east] at (3.7,0.8) (coord) {$\bx^5$};
    \draw[fill] (2.2,1.6) circle [radius=0.05];
    % \node[anchor=south east] at (2.2,1.5) (coord) {$\bx^6$};
    \draw[fill] (1.4,2.2) circle [radius=0.05];
    % \node[anchor=south east] at (1.4,1.8) (coord) {$\bx^7$};
    %------convex hull: -------
    \draw[dashed] (1.4,2.2)--(2.7,2.3);
    \draw[dashed](2.7,2.3)--(4,2);
    \draw[dashed] (4,2)--(3.7,0.8);
    \draw[dashed] (3.7,0.8)--(1.8,0.8);
    \draw[dashed] (1.8,0.8)--(1.4,2.2);
    \node[anchor=south east, gray] at (3.7, 1.1){\Large $\cH$};
    \fill[pattern=north east lines, pattern color=lightgray, opacity=0.5] (1.4,2.2)--(2.7,2.3)--(4,2)--(3.7,0.8)--(1.8,0.8)--cycle;  
    %------rejected ones:------
    \draw [fill, red] (2.6,.4) circle [radius=0.05];
    \draw [fill, red] (1.2,.8) circle [radius=0.05];
    \draw [fill, red] (1.3,2.8) circle [radius=0.05];
    \draw [fill, red] (2.3,2.5) circle [radius=0.05];
    \draw [fill, red] (3.3,2.7) circle [radius=0.05];
    \draw [fill, red] (4.7,1.2) circle [radius=0.05];
    \draw [fill, red] (1.1,1.5) circle [radius=0.05];
\end{tikzpicture}
\caption{The preferred solution $\bx^0$ has the best objective value among the convex hull of accepted observations.} %\hm{needs update}
\label{fig:preferred}
\end{figure}
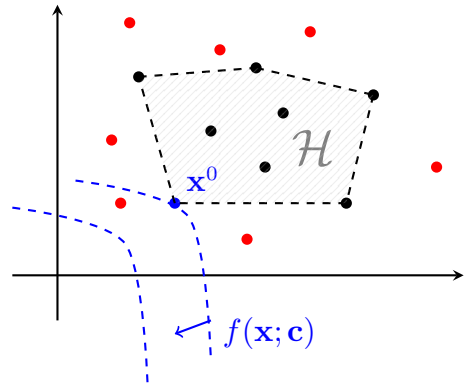

% \hm{Houra stopped here.}
Assumption~\ref{assumption:wellposed} ensures that the data and model are well-defined and that it is possible to construct a convex feasible region for FO. It states that (a) there exists a convex set that encapsulates the accepted observations, (b) the known constraints can indeed make the accepted observations feasible, and (c) the given input data do not contradict each other. These assumptions are not limiting given that the decision-maker is rational and is solving an implicit convex optimization problem. \changemarker{Note that if the data is not well-posed, a data pre-processing step and/or a notion of noise or error can be introduced in the inverse model. %, as often done in the literature {\kg{(add citations)}}. 
We further elaborate on such additions in Section~\ref{sec:future}.} %, meaning that the known constraints can indeed make the accepted observations feasible, and there exists a convex set that encapsulates the accepted observations. 
Because the objective function $f(\bx; \bc)$ is known in FO, we can identify the points in the convex hull of all accepted observations that provide the best objective value in FO. Definition~\ref{def:x0} characterizes such a point as the ``preferred solution''. %in the convex hull of all accepted observations. 
An example of the preferred solution for a convex objective function is visualized in Figure~\ref{fig:preferred}, where the blue dashed lines indicate iso-cost lines of the objective function and $\bx^0$ is the preferred solution. This concept is formally introduced in Definition~\ref{def:x0}. 

%and Definition~\ref{def:nominal} characterizes a potential solution that renders all accepted observations feasible and the rejected ones infeasible. 
\begin{definition}\label{def:x0} The~\emph{preferred} solution $\bx^0 \in \mathbb{R}^m$ is defined as %  {\hm{change to search over $\cH$}}
\begin{equation*}
% \bx^0 \in \argmin_{\bx^k, \, k\in\cKgood}\{ f(\bx; \bc) \}\,.    
\bx^0 \in \argmin_{\bx \in \cH}\{ f(\bx; \bc) \}\,.    
\end{equation*}
\end{definition}
% \begin{definition}\label{def:x0} The~\emph{preferred} solution $\bx^0 \in \cKgood$ is defined as   {\hm{change to search over $\cH$}}
% %
% \begin{equation*}
% \bx^0 \in \argmin_{\bx^k, \, k\in\cKgood}\{ f(\bx; \bc) \}\,.    
% \end{equation*}
% \end{definition}
%For FO to be convex, this preferred observation needs to be on the boundary a convex feasible set for FO that contains all accepted observations. Assumption~\ref{assumption:x0} is therfore logically necessary, as follows. 
%Visualization of preferred solution in Figure~\ref{fig:preferred}. 
%
% \begin{assumption}\label{assumption:x0}
% $\bx^0 \in \ext(\cH)$.  \hm{remove}
% \end{assumption}
% Assumption~\ref{assumption:x0} states that the preferred observation is well-behaved, meaning that it is not an interior point of the convex hull of observations. Without this assumption, the FO will be infeasible as there exists no convex feasible set that makes $\bx^0$ optimal. 

\noindent Note that depending on the type of objective function and the shape of the convex hull $\cH$, there may be multiple observations that satisfy the definition of a preferred solution, in which case, we arbitrarily label one of them as $\bx^0$. The preferred solution is not necessarily one of the observations, but it is always on the boundary of the convex hull of all observations.

The goal of this paper is to compute a set of linear and nonlinear constraints for the FO problem such that the accepted/rejected observations are inside/outside the inferred feasible region, respectively, and a preferred solution is an optimal solution for the FO model with the inferred feasible set.
% \begin{enumerate}
%     \item all the accepted observations are inside the inferred feasible region,
%     \item all rejected observations are outside of the inferred feasible region,
%     \item a preferred solution is optimal when the inferred feasible region is used in the FO model.
% \end{enumerate}
Hence, the intersection of the known constraints and the inferred constraints must include all the accepted observations and none of the rejected ones, providing a separation between the accepted and rejected points. 
\noindent We will refer to such a set of inferred constraints as a ``nominal set'', as formally defined in Definition~\ref{def:nominal}. 
\begin{definition}\label{def:nominal}
A convex set $\cXt$ is a \textit{nominal} set if 
\begin{align*}
\bx^k \in \cX \cap \cXt  &\qquad \forall k\in \cKgood,  %\label{eq:nominalGood} 
\\  
\bx^k \not \in \cX \cap \cXt &\qquad  \forall k\in \cKbad. %\label{eq:nominalBad}
\end{align*}
\end{definition}

\begin{figure}[htbp]
    \centering
\begin{tikzpicture}[thick,scale=1.2, every node/.style={scale=1.2}]
%[thick, scale=.8]
% % %----axes
\draw[->, >=stealth] (0,0) -- (4,0);
\draw[->, >=stealth] (0.5,-.5) -- (0.5,3);
% \node[anchor=south east] at (1.2,1.5) (coord) {$\bx^6$};
\node[blue] at (2.5,1.55) (coord) { $\cXt$};
%----------observations
    {\draw [fill] (2,2) circle [radius=0.07];
    \draw [fill] (3,2) circle [radius=0.07];
    \draw[fill] (1.8,1.2) circle [radius=0.07];
    \draw[fill] (2.7,0.8) circle [radius=0.07];
    \draw[fill] (1.2,1.6) circle [radius=0.07];
     \draw [fill] (1.4,.8) circle [radius=0.07];
    } %\pause
    {\draw [fill, red] (.8,.8) circle [radius=0.07];
    \draw [fill, red] (1,2) circle [radius=0.07];
    \draw [fill, red] (3,3) circle [radius=0.07];
    \draw [fill, red] (3.5,1) circle [radius=0.07];
    \draw [fill, red] (2,-.2) circle [radius=0.07];% we CAN have infeas obs in nominal.
    \draw [fill, red] (2.5,.2) circle [radius=0.07];
    }
%---inferred constraint
\draw[dashed] plot [smooth] coordinates{(1.1,-.2) (1,1.5) (2,2.5) (3.2,2.1) (2.7,-.5)} ; 
\draw[pattern=north east lines, pattern color=blue, opacity=0.5] plot [smooth] coordinates{(1.1,-.2) (1,1.5) (2,2.5) (3.2,2.1) (2.7,-.5)} ;
\node[anchor= west,blue] at (0.9,.5) (coord) {\small Nominal set};
%----known constraint
%\draw[lightgray] (0, .32) -- (4.2, .32); %known constraint
%\node[anchor=south east,gray] at (3.8,.3) (coord) {\small Known constraint};
% \fill[pattern=north west lines, pattern color=gray, opacity=0.5] (0, .32) -- (0,2.5) -- (3.7, 2.5) -- (3.7, .32) -- cycle;
\end{tikzpicture}  
%----------------------------------------------------------------------
%----------------------------------------------------------------------
%-------------------------------------------------------------------
\begin{tikzpicture}[thick,scale=1.2, every node/.style={scale=1.2}]
\node[anchor=east,gray] at (0,1.5) (coord) { \Large $\cap$};
%[thick, scale=.8]
% % %----axes
\draw[->, >=stealth] (0,0) -- (4,0);
\draw[->, >=stealth] (0.5,-.5) -- (0.5,3);
% \node[anchor=south east] at (1.2,1.5) (coord) {$\bx^6$};
\node[gray] at (2.5,1.555) (coord) { $\cX$};
%----------observations
    {\draw [fill] (2,2) circle [radius=0.07];
    \draw [fill] (3,2) circle [radius=0.07];
    \draw[fill] (1.8,1.2) circle [radius=0.07];
    \draw[fill] (2.7,0.8) circle [radius=0.07];
    \draw[fill] (1.2,1.6) circle [radius=0.07];
     \draw [fill] (1.4,.8) circle [radius=0.07];
    } %\pause
    {\draw [fill, red] (.8,.8) circle [radius=0.07];
    \draw [fill, red] (1,2) circle [radius=0.07];
    \draw [fill, red] (3,3) circle [radius=0.07];
    \draw [fill, red] (3.5,1) circle [radius=0.07];
    \draw [fill, red] (2,-.2) circle [radius=0.07];
    \draw [fill, red] (2.5,.2) circle [radius=0.07];
    }
%---inferred constraint
%\draw[red, dashed] plot [smooth] coordinates{(1.1,-.2) (1,1.5) (2,2.5) (3.2,2.1) (2.7,-.5)} ; 
%\draw[pattern=north east lines, pattern color=blue, opacity=0.5] plot [smooth] coordinates{(1.1,-.2) (1,1.5) (2,2.5) (3.2,2.1) (2.7,-.5)} ;
%----known constraint
\draw[gray] (0, .32) -- (4.2, .32); %known constraint
\node[anchor=south east,gray] at (3.8,.3) (coord) {\small Known constraint};
 \fill[pattern=north west lines, pattern color=gray, opacity=0.5] (0, .32) -- (0,2.5) -- (3.9, 2.5) -- (3.9, .32) -- cycle;
\end{tikzpicture}  
%-----------------------------------------------------------------------------
%-----------------------------------------------------------------------------
%-----------------------------------------------------------------------------
\begin{tikzpicture}[thick,scale=1.2, every node/.style={scale=1.3}]
\node[anchor=east,gray] at (0.2,1.5) (coord) {\large $=$};
%[thick, scale=.8]
\draw[->, >=stealth] (0,0) -- (4,0);
\draw[->, >=stealth] (0.5,-.5) -- (0.5,3);
% \node[anchor=south east] at (1.2,1.5) (coord) {$\bx^6$};
\node[black] at (2.5,1.55) (coord) { $\cX\cap\cXt$};
%----------observations
    {\draw [fill] (2,2) circle [radius=0.07];
    \draw [fill] (3,2) circle [radius=0.07];
    \draw[fill] (1.8,1.2) circle [radius=0.07];
    \draw[fill] (2.7,0.8) circle [radius=0.07];
    \draw[fill] (1.2,1.6) circle [radius=0.07];
     \draw [fill] (1.4,.8) circle [radius=0.07];
    } %\pause
    {\draw [fill, red] (.8,.8) circle [radius=0.07];
    \draw [fill, red] (1,2) circle [radius=0.07];
    \draw [fill, red] (3,3) circle [radius=0.07];
    \draw [fill, red] (3.5,1) circle [radius=0.07];
    \draw [fill, red] (2,-.2) circle [radius=0.07]; % we CAN have infeas obs in nominal.
    \draw [fill, red] (2.5,.2) circle [radius=0.07];
    }
%---inferred constraint
%\draw[red, dashed] plot [smooth] coordinates{(1.1,-.2) (1,1.5) (2,2.5) (3.2,2.1) (2.7,-.5)} ; 
%\draw[pattern=north east lines, pattern color=blue, opacity=0.5] plot [smooth] coordinates{(1.1,-.2) (1,1.5) (2,2.5) (3.2,2.1) (2.7,-.5)} ;
%----known constraint
 \draw[gray] (0, .32) -- (4, .32); %known constraint
% % \node[anchor=south east, lightgray] at (4,.32) (coord) {\footnotesize Known constraint};
%  \fill[pattern=north west lines, pattern color=gray, opacity=0.5] (0, .32) -- (0,2.5) -- (3.7, 2.5) -- (3.7, .32) -- cycle;
\draw[pattern=north east lines, pattern color=blue, opacity=0.5] plot [smooth] coordinates{(1.03,.32) (1,1.5) (2,2.5) (3.2,2.1) (2.92,.32)} ;
 \draw[pattern=north west lines, pattern color=gray, opacity=0.5] plot [smooth] coordinates{(1.03,.32) (1,1.5) (2,2.5) (3.2,2.1) (2.92,.32)} ;
% %  }
\end{tikzpicture}    
\caption{An illustration of the intersection of a nominal set and known constraints.}% The blue dashed area is the nominal set. The gray area is the known constraint.}
    \label{fig:nominalset}
\end{figure}
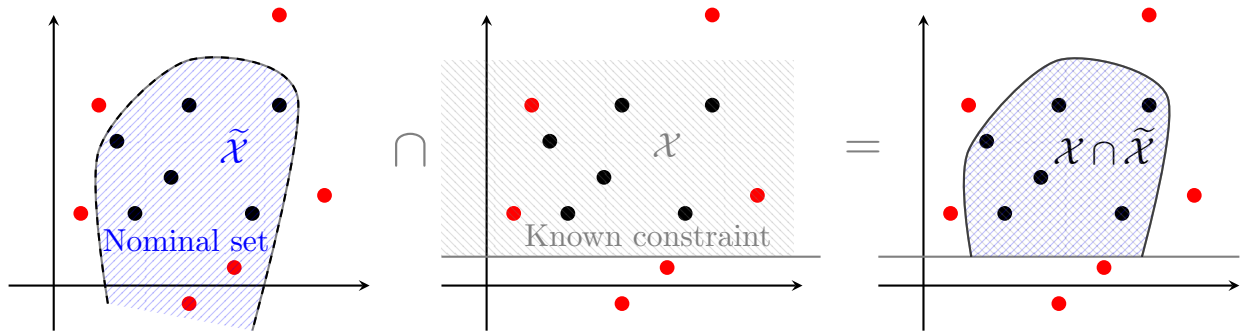

A simplified two-dimensional schematic of a nominal set is depicted in Figure~\ref{fig:nominalset}. \changemarker{We note that different sets of constraints can result in the same nominal set geometrically.} 
% Based on Definition~\ref{def:nominal}, a \emph{nominal} set is any convex set that when added to the known constraints renders all the accepted observations feasible and the rejected ones infeasible. A simple example of a nominal set and its interesection with the known constraints is provided in Figure~\ref{fig:nominalset}. %Based on Assumption~\ref{assumption:ObsAreFeas} and Definition~\ref{def:nominal}, 
% Based on these definitions, %to find a feasible set that satisfies the properties of it is sufficient 
Hence, the goal of the inverse problem is to find constraint parameters $\ba_{\ell}, b_\ell,\, \forall \ell \in \cLt$, and $\bq_n, \, \forall n \in \cNt$ such that the resulting inferred feasible set $\cXt$ %\[\cX = \left\{\,\bx\in \mathbb{R}^n \, \middle\vert \, g_n(\bx; \,\bq_n) \ge \bzero, \, \forall n \in \cNt, \,\,\, \ba_{\ell}' \bx \ge b_{\ell}, \, \forall \ell \in \cLt \, \right\}\] 
is a nominal set, and the preferred solution $\bx^0$ is an optimal solution for 
\begin{align*}  
\underset{\bx}{\text{minimize}} \quad &f(\bx; \bc) \\ 
 \text{subject to} \quad & \bx \in \cX \cap \cXt.
\end{align*}
%
% To find a set of constraint parameters with the properties above, 
To impute such constraints, we propose an \changemarker{expert-guided} inverse optimization (GIO) formulation that imposes feasibility constraints on the accepted observations, ensures the infeasibility of the rejected points, and enforces optimality conditions on the preferred solution $\bx^0$.  The \changemarker{GIO} model can be written as follows. 
\begin{subequations}\label{eq:DIO}
\begin{align} 
\mathbf{\changemarker{GIO}}: \quad \underset{\ba, b, \bq, \lambda, \mu,y}{\text{Maximize}} \quad  &  
\bP \left(\bq_1, \dots, \bq_{|\cNt|}, \bA, \bb; (\bx^{1}, \dots, \bx^{|\cK|}) \right) \label{eq:IOobj} 
\\ 
%  \text{Maximize} \quad & \max(distance([\bA, \bb]_\text{inferred},\bx^{k'})) & \\
 \text{subject to} \quad & g_n(\bx^k; \,\bq_n) \ge \bzero, \quad \forall k\in \cKgood ,\,  n\in \cNt  \label{eq:IOprimalfeasNL} \\ 
&   \ba_\ell'\, \bx^k \ge b_\ell,   \qquad \forall k\in \cKgood, \ell \in \cLt \label{eq:IOprimalfeasL}  \\
& \nabla f(\bx^0; \bc) + \sum_{n\in \cN\cup \cNt} \lambda_n \nabla g_n\, (\bx^0,\, \bq_n) + \sum_{\ell \in \cL \cup \cLt}\mu_\ell \, \ba_{\ell} \, = \bzero, \label{eq:IOstationarity}  \\
% \quad & \lambda_n \,\, g_n(\bx^0, \bq_n) = 0 \qquad \forall n\in \cN\cup \cNt, \label{eq:IOcs}  %CS for equal linear consrts.\\
& \lambda_n \,\, g_n(\bx^0, \bq_n) = 0, \qquad \forall n\in \cN\cup \cNt \label{eq:IOcsNL}  \\
& \mu_\ell \,\, (b_\ell - \ba_\ell' \, x^0)= 0, \qquad \forall \ell\in \cL\cup \cLt \label{eq:IOcsL}  \\
% & f(\bx^k; \bc) \geq f(\bx^0; \bc) \quad \forall k \in \cKgood \qquad \text{(half-space)} \label{eq:IOhalfspace} \\
& \ba_{\ell}'\,\bx^k \le b_{\ell} - \epsilon + M y_{\ell k}, \qquad \forall \ell \in \cL \cup \cLt, \,  k \in \cKbad \label{eq:IOinfeasL}\\
& g_n(\bx^k; \,\bq_n) \le \bzero - \epsilon + M y_{nk}, \qquad \forall n \in \cN \cup \cNt, \, k \in \cKbad \label{eq:IOinfeasNL}\\
& \sum_{i \in \cI } y_{ik} \leq \, \mid \cI \mid - 1, \qquad \forall k \in \cKbad \label{eq:IOinfeasSum}\\
%& \mathbf{1}' \by^{\cL}_k + \mathbf{1}' \by^{\cN}_k \leq |\cI| - 1, \quad \forall k \in \cK' \\
& \ba_{\ell} \in \cA_{\ell} \,, b_{\ell} \in \cB_{\ell}, \qquad \forall \ell \in \cLt  \label{eq:IOnormL} \\
& \bq_n \in \cQ_n \,, \qquad \forall n \in \cNt  \label{eq:IOnormNL} \\ %||\bq_n||= 1
& \lambda_n,\mu_\ell \leq 0, \qquad \forall n\in \cN\cup \cNt, \, \ell\in \cL\cup \cLt  \label{eq:IOstationaritySign}\\
%& y_\ell, y_n \in \{0,1\}, \qquad \forall n\in \cN\cup \cNt, \, \ell\in \cL\cup \cLt
& y_{ik} \in \{0,1\}, \qquad \forall i \in \cI, \, k \in \cKbad. \label{eq:IOBinaryVars}
\end{align}
\end{subequations}

%\hm{How many constraints is (2d)? One? m? isn't gradient a vector? KG: So, we decided it'd be m)??}
%\kg{key results from this section: 1) infeasible obs, 2) nonlinear half-space}
\noindent The objective function~\eqref{eq:IOobj} maximizes a measure of distance between the constraint parameters and the observations. An example of the objective function can be maximizing the total distance between the inferred constraints and all the infeasible observations using a desirable distance matrix $\bP$. We provide more  details on this objective function example in Section \ref{sec:distance}. Constraints \eqref{eq:IOprimalfeasNL} and \eqref{eq:IOprimalfeasL} enforce primal feasibility conditions. Constraints~\eqref{eq:IOstationarity} capture the stationarity conditions. Complementary slackness for the linear and nonlinear constraints of FO are captured in \eqref{eq:IOcsNL} and \eqref{eq:IOcsL}, respectively. Constraints~\eqref{eq:IOinfeasL}--\eqref{eq:IOinfeasSum} ensure that at least one constraint is violated by each of the rejected observations. Constraints~\eqref{eq:IOnormL}--\eqref{eq:IOnormNL} provide a set of desirable conditions on the coefficients of the imputed constraints such as normalization or convexity conditions. As an optional step, any other desirable condition on the parameters can also be included in $\cQ$, and similar conditions on the linear constraint parameters can also be considered as $\ba_{\ell} \in \cA, \, b_{\ell} \in \cB$. Lastly, constraints~{\eqref{eq:IOstationaritySign}--\eqref{eq:IOBinaryVars}} indicate sign and binary declarations. 
%
% {\red{somewhere around here we need to formally define $\cXt$ in terms of the inferred constraint parameters. We keep using it without context.}}
%\subsection{DIO Properties: }
%In this section, we discuss some of the properties of 

We next show that any optimal solution produced by the \changemarker{GIO} model exhibits the desired properties of an inferred feasible region for FO. %Specifically, Proposition~\ref{prop:DIOoptimal} shows that\changemarker{GIO}produces a nominal set that makes $\bx^0$ optimal for FO. 

\begin{proposition} \label{prop:DIOoptimal}
% Any feasible solution of \changemarker{GIO} is a nominal set such that %makes $\bx^0$ optimal for $\FO$, i.e., $\bx^0 \in \underset{x\in\cX \cap \cXt}{\argmin}\,\{f(\bx; \bc)\}$. %satisfies all four requirements Req. 1-4.
Any feasible solution of \changemarker{GIO} corresponds to a nominal set $\cXt$  
%$\cXt=\{\bx \in \mathbb{R}^m \mid g_n(\bx;,\bq_n) \geq \bzero, \forall n \in \cNt; \, \ba_{\ell}'\bx \geq b_{\ell}, \forall \ell \in \cLt \}$ 
such that %makes $\bx^0$ optimal for $\FO$, i.e., 
${\bx^0 \in \underset{\bx \in\cX \cap \cXt}{\argmin}\,\{f(\bx; \bc) \}}$.
\end{proposition}
As Proposition~\ref{prop:DIOoptimal} states, any solution of \changemarker{GIO} has the properties of a nominal feasible set for FO and makes $\bx^0$ optimal for the forward problem. To fulfill this requirement, \changemarker{GIO} inserts at least one conflicting constraint per rejected observation such that the rejected observation becomes infeasible for FO while ensuring none of the accepted observations are cut off. Assuming that the forward model allows us to infer as many constraints as needed to do so, then it is always possible to find a solution for \changemarker{GIO}. %Proposition~\ref{prop:DIOfeas} formalizes this concept. 

\begin{proposition} \label{prop:DIOfeas}
For sufficiently large $|\cLt| + |\cNt|$, \changemarker{GIO} is guranteed to be feasible. 
%the feasible region of \changemarker{GIO} is guranteed to be non-empty.
%For $g_n$, use the well-defineness assumption to show that there exists a $q_n$
\end{proposition}
\noindent Proposition~\ref{prop:DIOfeas} states that \changemarker{GIO} is feasible when the number of inferred constraints is sufficiently large. Next, in Remark~\ref{remark:minNumConstr},  we construct an upper bound on the minimum number of constraints needed to make \changemarker{GIO} feasible. 
%Remark~\ref{remark:minNumConstr} constructs an upper bound on the minimum number of constraints stated in Proposition~\ref{prop:DIOfeas} to make \changemarker{GIO} feasible. 

% Depending on the number of rejected observations and their spatial distribution, a small number of constraints may be sufficient to cut off a large number of rejected observations. However, in the worst case, we would need one constraint per rejected observation to guarantee that %the feasible region of the \changemarker{GIO} model is nonempty, %. We note that we need at least one inferred constraint that makes 
% each rejected point is infeasible for FO while ensuring the feasibility of all accepted points. We also need at least one inferred constraint ensuring that the preferred solution $\bx^0$ is optimal for FO. 
% Remark~\ref{remark:minNumConstr} provides a bound for the number of inferred constraints required for \changemarker{GIO}.

%\vspace{0.5em}
\begin{remark} \label{remark:minNumConstr}
An upper bound to the minimum number of required constraints in \changemarker{GIO} is {$|\cKbad|+1$}. 
%And the minimum number of constraints in adversarial example is $\min \{ |\cKbad|+1,|\CH|\},$ where $\cH$ is the convex hull of $\bx^k, \forall k \in \cKgood \}$ and $|\CH|$ denotes the number of facets in $\CH$.  
%$\min \{ |\cKbad|$, number of facets of convex hull of $\bx^k, \forall k \in \cKgood \}$
\end{remark}  %\vspace{0.5em}

\noindent Depending on the number of rejected observations and their spatial distribution, a small number of constraints may be sufficient to cut off a large number of rejected observations. However, in the worst case, we would need one constraint per rejected observation to guarantee that %the feasible region of the \changemarker{GIO} model is nonempty, %. We note that we need at least one inferred constraint that makes 
each rejected point is infeasible for FO while ensuring the feasibility of all accepted points. We also need at least one inferred constraint ensuring that the preferred solution $\bx^0$ is optimal for FO. We note again that one constraint may serve multiple purposes, which would result in a lower number of constraints needed in practical settings. For instance, a single constraint may cut a large number of rejected points out of the inferred feasible region. 
The proposed \changemarker{GIO} model~\eqref{eq:DIO} can be a very complex mixed-integer nonlinear problem which can pose a challenge for state-of-the-art commercial solvers.
%and we found that commercial solvers often fail to find any good-quality solutions even for small problem instances. %of solutions found by solvers. 
In Section~\ref{sec:solution}, we propose a method for mitigating the complexity of the proposed inverse problem.

% \hm{lead to problem complexity}

\section{Mitigating Inverse Problem Complexity}\label{sec:solution}
The complexity of the \changemarker{GIO} model~\eqref{eq:DIO} depends on the complexity of the FO model~\eqref{eq:FO}. %, %e.g., it can be a mixed-integer nonlinear nonconvex problem, which can pose a challenge for state-of-the-art commercial solvers. % to find good-quality solutions.
% In our empirical analysis, commercial solvers often fail to find good-quality solutions even in small problem instances. 
% The complexity of \changemarker{GIO} model~\eqref{eq:DIO} is largely due to 
Additionally, it includes KKT conditions of complementary slackness and stationarity, which can be nonlinear and nonconvex. 
% \hm{emphasize that there is not much work on constraint inference because it is computationally intractable due to highly nonlinear IO formulation, regardless of the FO complexity. However, there is hope! in LP, there is a  method that can be generalized to convex, and can partially mitigate the complexity of the resulting IO problem. }
% The \changemarker{GIO} formulation introduced in Section \ref{sec:IOFormulation} infers constraints based on observations for convex nonlinear problems. The formulation itself however is also nonlinear and more complex than the original FO problem. For instance, constraints \eqref{eq:eIOcsNL} may or may not be convex depending on the choices of the $g_n$. 
% When inferring linear constraints in forward LP models, \citep{ghobadi2021inferring} introduced a method to leverage the structure of a \kg{class of solutions} in inverse models that infer constraints. In their framework, they introduced an additional linear constraint that ensured strong duality and hence simplified their formulation to a problem that could be linear.  
% eliminate the strong duality constraint of the inverse optimization models when constraints are inferred. 
% \hm{From lit rev: The resulting inverse model is nonlinear, but their method leverages the structure of a class of solutions to provide an equivalent tractable reformulation that can be reduced to a linear program. }
%
In the context of recovering linear constraint sets, prior work used linear programming duality to replace the optimality conditions % in their linear forward problem 
and provide a tractable reformulation~\citep{ghobadi2021inferring}. %To do so, they add a linear half-space as a known constraint that would eliminate the need for writing strong duality and dual feasibility conditions in the inverse model. 
\changemarker{To mitigate the complexity of inferring nonlinear constraints, we propose a tractable reformulation for solving GIO leveraging the concept of variational inequalities \citep{harker1990finite, kinderlehrer2000introduction,bertsimas2015data} to identify a linear constraint that can enforce optimality conditions instead of the nonlinear KKT conditions, generalizing the approach of \cite{ghobadi2021inferring} to nonlinear cases.}

In what follows, we first introduce a few definitions and discuss preliminaries for constructing the reformulation. We then discuss the properties of an optimal inverse solution and provide problem-specific sufficient optimality conditions to replace the KKT criteria. Lastly, we present a reduced reformulation of the inverse problem.  %\changemarker{Throughout this section, we assume that the given known objective function $f(\bx;\bc)$ is convex and monotone.}%, % which is easier to solve than the original \changemarker{GIO} model. %retains the complexity of its FO counterpart. 

% \hm{*** Houra to edit next: }

% \hm{Do we need these subsections?}
%\subsection{Sublevel Set and Tangent Hyperplane}
\subsection{Preliminaries and Definitions}

% \hm{
% Variational inequality; for our objective $f(\bx;\bc)$\\
% Let $\cX$ be a convex set in $\mathbb{R}^n$ and a function $f(\bx;\bc): \mathbb{R}^n\rightarrow \mathbb{R}$ be a differentiable convex function. Then, the variational inequality problem, denoted as $VI(f, \cX)$  is to find $\bx^0$ such that }
% \\
% variational inequality: $\nabla f(\bx^0; \bc)'(\bx - \bx^0) \geq 0$
% \\

% tangent half-space: $\cC = \{\bx \in \mathbb{R}^n  \mid (\nabla f(\bx^0;\bc))' \,\bx \geq f(\bx^0; \bc)  \}.$ 
% \\

% sublevel set: $\cV = \{\bx \mid f(\bx;\bc) \geq f(\bx^0; \bc)\}.$
% \\

% If $f=\bc'\bx$ then $\nabla f(\bx^*)'(\bx - \bx^*)\geq 0 \Rightarrow  c'(\bx-\bx^0)\geq 0$ all of the above are just our old half-space. 

Recall that the preferred solution $\bx^0$ is optimal for FO \changemarker{and has a} %, meaning that it has a 
better objective value for  $f(\bx;\bc)$ than any other accepted observation. %, i.e., %. Consider the convex objective function of FO, \[ \text{minimize}_{\bx} \,\, f(\bx; \bc). \] 
%This means that 
%$f(\bx^0;\bc) \leq f(\bx^k; \bc), \, \forall k \in \cKgood$. 
%This relationship can be captured using the concept of variational inequalities and formalized in the following definition. 
\changemarker{This characteristic relates to the well-studied concept of variational inequalities \citep{harker1990finite}, which we use to define a sublevel set in definition~\ref{def:halfspace}. }% we define a sublevel set of the objective as the set of all solution that are dominated by $\bx^0$.}
%In a simplified two-dimensional setting, if we draw the iso-cost objective function curve at $f(\bx;\bc) = f(\bx^0;\bc)$, all accepted observations fall on one side of this curve. Definition~\ref{def:halfspace} formally defines this space on one side of the curve as the sublevel set of the objective function at the preferred solution\changemarker{, using the concept of variational inequalities \citep{harker1990finite}}. 
%
%
\begin{definition} \label{def:halfspace}
The set $\cV$ is a \textit{sublevel set} of $f(\bx;\bc)$ at $\bx^0$ defined as 
\[\cV = \{\bx \in \mathbb{R}^n \mid f(\bx;\bc) \geq f(\bx^0; \bc)\}. \] 
\end{definition}
%\hm{can remove this sentence. We just aid the same thing up there:} %Note that the sublevel set $\cV$  always contains all accepted observations, because any point $\hat{\bx} \not \in \cV$ would have a better objective value than $\bx^0$, i.e.,  $f(\hat{\bx};\bc) < f(\bx; \bc)$, which contradicts with $\bx^0$ being the preferred solution.\\
\changemarker{Note that the sublevel set of the objective function at the preferred solution $\bx^0$ contains all accepted observations because $f(\bx^0;\bc) \leq f(\bx^k; \bc), \, \forall k \in \cKgood$. %, because any point $\hat{\bx} \not \in \cV$ would have a better objective value than $\bx^0$, i.e.,  $f(\hat{\bx};\bc) < f(\bx; \bc)$, which contradicts with $\bx^0$ being the preferred solution. 
Solving $\underset{\bx \in \cX}{\text{minimize }} f(\bx;\bc)$ to optimality is equivalent to finding $\bx^0$ that satisfies these variational inequality conditions. Since $f(\bx;\bc)$ is convex and differentiable, we use its first-order condition to define a  
%Next, we define a 
linear tangent to the sublevel set at $\bx^0$, as stated in Definition~\ref{def:tangent}.}

\begin{figure}[htbp]
    \centering
\begin{tikzpicture}[thick,scale=1.2, every node/.style={scale=1.3}]
% \node[anchor=east,gray] at (0.2,1.5) (coord) {\large $=$};
%[thick, scale=.8]
\node[anchor=east,blue] at (1.2,2.65) (coord) {\large $\cV$};
\draw[->, >=stealth] (0,0) -- (4,0);
\draw[->, >=stealth] (0.5,-.5) -- (0.5,3);
%----------observations--------
    {\draw [fill] (2,2) circle [radius=0.07];
    \draw [fill] (3,2) circle [radius=0.07];
    \draw[fill] (1.8,1.2) circle [radius=0.07];
    \draw[fill] (2.7,0.8) circle [radius=0.07];
    \draw[fill] (1.2,1.6) circle [radius=0.07];
     \draw [fill, blue] (1.4,.76) circle [radius=0.07];
    } %\pause
    {\draw [fill, red] (.8,.8) circle [radius=0.07];
    \draw [fill, red] (1,2) circle [radius=0.07];
    \draw [fill, red] (3,3) circle [radius=0.07];
    \draw [fill, red] (3.5,1) circle [radius=0.07];
    \draw [fill, red] (2,-.2) circle [radius=0.07]; % we CAN have infeas obs in nominal.
    \draw [fill, red] (2.5,.2) circle [radius=0.07];
    }
%---inferred constraint--------
% \draw[red, dashed] plot [smooth] coordinates{(1.1,-.2) (1,1.5) (2,2.5) (3.2,2.1) (2.7,-.5)} ; 
%\draw[pattern=north east lines, pattern color=blue, opacity=0.5] plot [smooth] coordinates{(1.1,-.2) (1,1.5) (2,2.5) (3.2,2.1) (2.7,-.5)} ;
%----known constraint--------
\draw[gray] (0, .32) -- (4, .32); %known constraint
\draw[pattern=north east lines, pattern color=gray, opacity=0.4] plot [smooth] coordinates{(1.03,.32) (1,1.5) (2,2.5) (3.2,2.1) (2.92,.32)} ;
%--------sublevel set of f(x;c) --------
% \draw[->, blue]  (0.1,.5) coordinate -- (-.3,.15) ;
\draw[ dashed , blue] plot [smooth] coordinates{(0.1,1.15) (1.55,0.65) (2.1, -.95)} node[anchor= south east] {\footnotesize $f(\bx;\bc)$}; 
% \draw[pattern=north west lines, pattern color=blue, opacity=0.4] plot [smooth] coordinates{(0.1,1.15) (1.55,0.65) (2.1, -.95) (4, 0.5) (4,4) (0,4) } ;
\fill[pattern=north west lines, pattern color=blue,  opacity=.5] (0.1,1.15) -- (1.6,0.68) -- (2.1, -.95) -- (4, 0.5) -- (4,3) -- (0,3) -- cycle;
\end{tikzpicture} 
\hspace{4em}%%%%%%%%%%%%%%%%%%%%%%%%%%%%%%%%%%%%%%%%%%%%%%%%%%%%%%%%%%%%%%%%
\begin{tikzpicture}[thick,scale=1.2, every node/.style={scale=1.3}]
\node[anchor=east,black] at (1.2,2.65) (coord) {\large $\cC$};
\draw[->, >=stealth] (0,0) -- (4,0);
\draw[->, >=stealth] (0.5,-.5) -- (0.5,3);
%----------observations--------
    {\draw [fill] (2,2) circle [radius=0.07];
    \draw [fill] (3,2) circle [radius=0.07];
    \draw[fill] (1.8,1.2) circle [radius=0.07];
    \draw[fill] (2.7,0.8) circle [radius=0.07];
    \draw[fill] (1.2,1.6) circle [radius=0.07];
     \draw [fill, blue] (1.4,.8) circle [radius=0.07];
    } %\pause
    {\draw [fill, red] (.8,.8) circle [radius=0.07];
    \draw [fill, red] (1,2) circle [radius=0.07];
    \draw [fill, red] (3,3) circle [radius=0.07];
    \draw [fill, red] (3.5,1) circle [radius=0.07];
    \draw [fill, red] (2,-.2) circle [radius=0.07]; % we CAN have infeas obs in nominal.
    \draw [fill, red] (2.5,.2) circle [radius=0.07];
    }
%---inferred constraint--------
% \draw[red, dashed] plot [smooth] coordinates{(1.1,-.2) (1,1.5) (2,2.5) (3.2,2.1) (2.7,-.5)} ; 
%\draw[pattern=north east lines, pattern color=blue, opacity=0.5] plot [smooth] coordinates{(1.1,-.2) (1,1.5) (2,2.5) (3.2,2.1) (2.7,-.5)} ;
%----known constraint--------
\draw[gray] (0, .32) -- (4, .32); %known constraint
\draw[pattern=north east lines, pattern color=gray, opacity=0.4] plot [smooth] coordinates{(1.03,.32) (1,1.5) (2,2.5) (3.2,2.1) (2.92,.32)} ;
%--------sublevel set of f(x;c) --------
% \draw[->, blue]  (0.1,.5) coordinate -- (-.3,.15) ;
\draw[ dashed , blue] plot [smooth] coordinates{(0.1,1.15) (1.55,0.65) (2.1, -.95)} node[anchor= south east] {\footnotesize $f(\bx;\bc)$}; 
\draw[black] (0,1.5) -- (4,-0.5);
\fill[pattern=north west lines, pattern color=black,  opacity=.4] (0,1.5) -- (4,-0.5) -- (4, 0.5) -- (4,3) -- (0,3) -- cycle;
\end{tikzpicture} 
\caption{The sublevel set of $f(\bx,\bc)$ (left) and its tangent half-space (right) at $\bx^0$. }%Left: The blue shaded  area is the sublevel set of $f(\bx,\bc)$ at $\bx0$. The red dashed area is the inferred feasible set.}
    \label{fig:sublevelset1}
\end{figure}
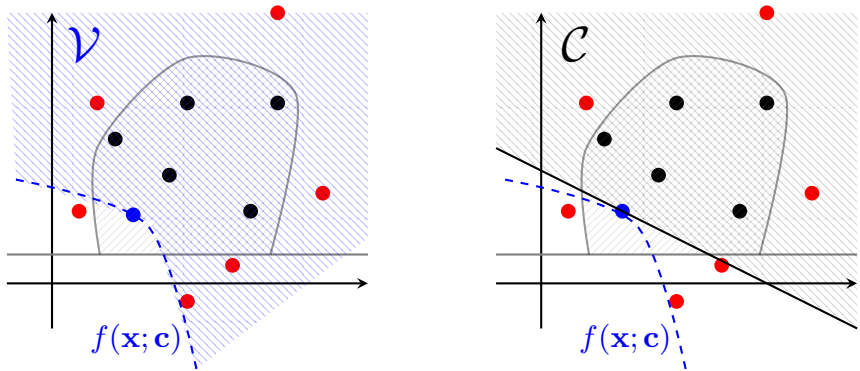

\begin{definition}\label{def:tangent}
The tangent half-space $\cC$ to the sublevel set of %\changemarker{a convex monotone} 
$f(\bx;\bc)$ at $\bx^0$ is defined as 
% \[ \cC = \{\bx \in \mathbb{R}^n  \mid (\nabla f(\bx^0;\bc))' \,\bx \geq f(\bx^0; \bc)  \}. \text{wrong old}\]
\[ \cC = \{\bx \in \mathbb{R}^n  \mid (\nabla f(\bx^0;\bc))' \,\bx \geq \nabla f(\bx^0; \bc)'\bx^0 \}.\]
\end{definition}
%
% Given that $f(\bx;\bc)$ is convex, the tangent half-space $\cC$ is equivalent to the variation inequality conditions for FO which satisfy the first-order optimality conditions of $\underset{\bx \in \cX}{\text{minimize }} f(\bx;\bc)$ 
% . These conditions are necessary and sufficient for the first-order optimality conditions of $\bx^0$ for the FO problem.  %on the optimality of $\bx^0$ 
%for the FO problem, 
Simplified schematics of a sublevel set and its tangent half-space are shown in Figure~\ref{fig:sublevelset1}. 
%
%the first-order optimality conditions of the FO problem. 
%Therefore, the tangent half-space $\cC$ can be used to 
% replace the KKT conditions since \changemarker{GIO} is convex. 
It can be seen that if $f(\bx;\bc)$ is linear, the tangent half-space is equivalent to the sublevel set of $f(\bx;\bc)$ at $\bx^0$, i.e., $\cC=\cV$. 
% Solving $\underset{\bx \in \cX}{\text{minimize }} f(\bx;\bc)$ to optimality \changemarker{for a convex mononone $f$} is equivalent to finding $\bx^0$ that satisfies these variational inequality conditions in $\cC$. 
In the rest of this paper, we use $\cC$ to denote the tangent half-space of the sublevel set of $f(\bx;\bc)$ at $\bx^0$, for brevity.  %Proposition~\ref{lem:tangentisfeasible} shows that $\cC$ always includes all accepted observations.

%Known constraint: $(\nabla f(\bx^0;\bc))' \bx - f (\bx^0; \bc) \geq 0 $. Note that $\nabla f(\bx^0; \bc)$ is the value of $\nabla f(\bx; \bc) $ at point $\bx^0$. 

% %$\bx \in \cS \in {\cX} \cap \bar{\cX}$, 
% $\bx \in{\cX} \cap \bar{\cX}$, for any $\bar{\cX}$  that has $\bx^0 \in \underset{\bx \in \bar{\cX}}{\arg\min} \{f(\bx;\bc)\}$. This means that for such $\bar{\cX}$, we always have $\bar{\cX} \subseteq \cC$. The reason is that the sub-level set at $\bx^0$ only contains points that have a worse objective value $f(\bx,\bc)$ than $\bx^0$. 
%Note that the sublevel set is not necessarily convex, but its intersection with the inferred constraints will always be convex. 
%

\subsection{Reduced Formulation}
%\hm{STOPPED HERE. Next time, read the rest from here to see if we can use the concept of VIs to better streamline the proofs. We could also add another proposition or lemma using VIs to show equivalence.}

%{In this section, we use the theoretical properties of a \changemarker{GIO} solution to enforce optimality conditions through a reduced reformulation that is less complex than the original \changemarker{GIO} model. }
One of the key complexities of the \changemarker{GIO} formulation is the inclusion of nonlinear KKT conditions for stationarity and complementary slackness, which ensure the optimality of the preferred solution. In this section, \changemarker{we leverage the concept of variational inequalities, which characterize solutions to inequalities over a feasible region, to find alternative formulation to the nonlinear KKT conditions. To this end, we use the half-space $\cC$ that is defined based on variational inequalities (see Definition~\ref{def:halfspace}). This half-space is tangent to sublevel sets of $f(\bx; \bc)$ at $\bx^0$, and hence, enforces the optimality conditions of $x^0$, as detailed in Proposition~\ref{prop:CcapS}. By replacing the KKT conditions with half-space $\cC$, we arrive at a reduced formulation with less complexity than the original \changemarker{GIO} model.}  
%we employ the tangent half-space previously introduced in Definition~\ref{def:halfspace} to enforce optimality conditions and propose reduced reformulation that is less complex than the original \changemarker{GIO} model. 
For brevity of notations, we refer to any feasible region that is inferred using the \changemarker{GIO} model as an \emph{imputed set} for the forward problem. Given the properties of \changemarker{GIO} solution outlined in Proposition~\ref{prop:DIOoptimal}, Definition~\ref{def:imputed} characterizes an imputed set for FO. 
\begin{definition}\label{def:imputed}
Any convex nominal set $\cS=\cX \cap \cXt$ such that $\bx^0 \in \underset{\bx \in\cX \cap \cXt}{\argmin}\,\{f(\bx; \bc) \}$ is an \emph{imputed set} for FO. 
\end{definition}
%
% Based on Proposition~\ref{prop:DIOoptimal} and Definition~\ref{def:imputed}, any optimal solution of \changemarker{GIO} corresponds to an imputed set. Remark~\ref{rem:DIOimputed} accurately states this claim. 

% \begin{remark}\label{rem:DIOimputed}
% An optimal solution of \changemarker{GIO} corresponds to an imputed set for FO. 
% \end{remark}

% In the rest of this paper, we refer to an FO feasible region that is based on a \changemarker{GIO} solution as an imputed set. 
We note that different \changemarker{GIO} solutions may result in the same imputed set since it is a geometric representation of the feasible region, as opposed to an algebraic one. For instance, in a \changemarker{GIO} solution, multiplying the coefficients of a linear constraint by a constant would result in a different solution, which may even be infeasible for \changemarker{GIO}, but it would represent the same feasible set for FO.
%
% Recall the definition of tangent half-space $\cC$ of the sublevel set of $f(\bx;\bc)$ at $\bx^0$ in Definition~\ref{def:halfspace}. 
Proposition~\ref{prop:imputedinC} shows that any imputed set as defined in Definition~\ref{def:imputed} is always contained within the tangent half-space $\cC$. %This property is detailed in Proposition~\ref{prop:imputedinC}.

% \begin{proposition} \label{prop:tangentisfeasible}
% \hm{repetitive} For any convex $f(\bx;\bc)$, the tangent half-space of the sublevel set of $f(\bx;\bc)$ containts all accepted observations. i.e., $\bx^k \in  \cC, \, \forall k \in \cKgood$. 
% \end{proposition}
% \begin{proof}{Proof.}
% We know that $\bx^k \in \cV, \, \forall k \in \cKgood$. Assume that $\exists \bx^{k_1} \in \cV\setminus\cC$. Then, consider the two cases: 
% {\bf Case 1:} If $f(\bx;\bc)$ is linear, then $\cC=\cV$ and $\cV\setminus\cC = \varnothing$, which is in contradition with $\bx^k \in \cV\setminus\cC$. 

% \noindent {\bf Case 2:} If $f(\bx;\bc)$ is convex nonlinear, then $\cV$ is non-convex, and hence $\cV\setminus\cC$ is non-convex. Hence, it must be that $\exists \bx^{k_2} \in \cV\setminus \cV$ such that $\bx^0 = \lambda \bx^{k_1} + (1-\lambda) \bx^{k_2}$ for $\lambda >0$. This contradicts the definition of preferred solution. 
% \end{proof}

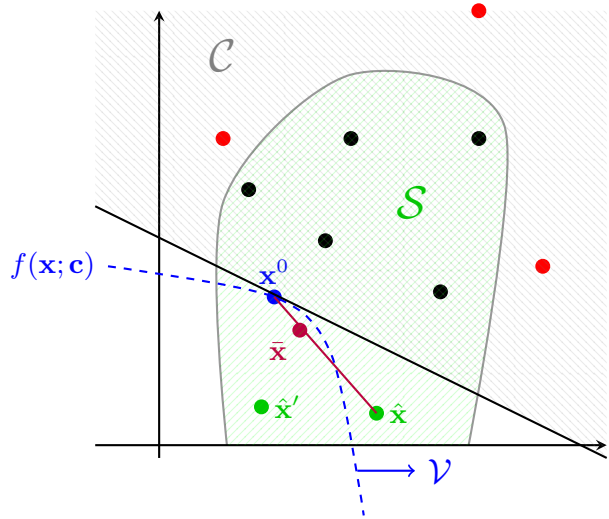
\begin{figure}[htbp]
    \centering
\begin{tikzpicture}[thick,scale=1.7, every node/.style={scale=1.2}]
% \node[anchor=east,gray] at (0.2,1.5) (coord) {\large $=$};
%[thick, scale=.8]
\node[anchor=east,gray] at (1.2,2.65) (coord) {\large $\cC$};
\node[anchor=east,myGreen] at (2.7,1.5) (coord) {\large $\cS$};
\draw[->, >=stealth] (0,-.4) -- (4,-.4);
\draw[->, >=stealth] (0.5,-.5) -- (0.5,3);
%----------observations--------
    {\draw [fill] (2,2) circle [radius=0.05];
    \draw [fill] (3,2) circle [radius=0.05];
    \draw[fill] (1.8,1.2) circle [radius=0.05];
    \draw[fill] (2.7,0.8) circle [radius=0.05];
    \draw[fill] (1.2,1.6) circle [radius=0.05];
    \draw [fill, blue] (1.4,.76) circle [radius=0.05] node[anchor= south] {\footnotesize $\bx^0$}; 
    %  \draw [fill, myGreen] (1.85,.45) circle [radius=0.05] node[anchor= north west] {\footnotesize $\hbx_2$}; 
    \draw [fill, myGreen] (2.2,-.15) circle [radius=0.05] node[anchor= west] {\footnotesize $\hbx$}; 
    \draw [fill, myGreen] (1.3,-.1) circle [radius=0.05] node[anchor= west] {\footnotesize $\hbx'$}; 
    \draw[purple] (1.4,.76) -- (2.2,-.15);
    \draw [fill, purple] (1.6,.5) circle [radius=0.05] node[anchor= north east] {\footnotesize $\bar{\bx}$};
} %\pause
    {
    \draw [fill, red] (1,2) circle [radius=0.05];
    \draw [fill, red] (3,3) circle [radius=0.05];
    \draw [fill, red] (3.5,1) circle [radius=0.05];
    % \draw [fill, red] (2,-.2) circle [radius=0.07]; % we CAN have infeas obs in nominal.
    % \draw [fill, red] (2.5,.2) circle [radius=0.07];
    }
%---inferred constraint--------
% \draw[red, dashed] plot [smooth] coordinates{(1.1,-.2) (1,1.5) (2,2.5) (3.2,2.1) (2.7,-.5)} ; 
%\draw[pattern=north east lines, pattern color=blue, opacity=0.5] plot [smooth] coordinates{(1.1,-.2) (1,1.5) (2,2.5) (3.2,2.1) (2.7,-.5)} ;
%----known constraint--------
% \draw[gray] (0, .32) -- (4, .32); %known constraint
%--------IMPUTED SET (WRONG on purpose)
\draw[pattern=north east lines, pattern color=green, opacity=0.4] plot [smooth] coordinates{(1.03,-.4) (1,1.5) (2,2.5) (3.2,2.1) (2.92,-.4)} ;
%--------sublevel set of f(x;c) --------
% \draw[->, blue]  (0.1,.5) coordinate -- (-.3,.15) ;
\draw[ dashed , blue] plot [smooth] coordinates{(0.1,1) (1.4,.76) (1.85,0.3) (2.1, -.95)};
\node[anchor= east, blue] at (0.1,1) (coord){\footnotesize $f(\bx;\bc)$}; 
% \draw[pattern=north west lines, pattern color=blue, opacity=0.4] plot [smooth] coordinates{(0.1,1.15) (1.55,0.65) (2.1, -.95) (4, 0.5) (4,4) (0,4) } ;
% \fill[pattern=north west lines, pattern color=blue,  opacity=.5] (0.1,1.15) -- (1.6,0.68) -- (2.1, -.95) -- (4, 0.5) -- (4,3) -- (0,3) -- cycle;
%--------sublevel set of f(x;c) --------
% \draw[->, blue]  (0.1,.5) coordinate -- (-.3,.15) ;
% \draw[ dashed , blue] plot [smooth] coordinates{(0.1,1.15) (1.55,0.65) (2.1, -.95)} node[anchor= south east] {\footnotesize $f(\bx;\bc)$}; 
\draw[->, blue] (2.05,-.6) -- (2.5,-.6) node [anchor = west]{$\cV$};
% \draw[->, blue] (2.05,-.6) -- (2.5,-.6) node [anchor = west]{$\cV$};
\draw[black] (0,1.47) -- (4,-0.5);
\fill[pattern=north west lines, pattern color=gray,  opacity=.4] (0,1.5) -- (4,-0.5) -- (4, 0.5) -- (4,3) -- (0,3) -- cycle;
\end{tikzpicture} 
\caption{A convex nominal set $\cS \not \subseteq \cC$, to show the intuition behind the proof of Proposition~\ref{prop:imputedinC}}%Left: The blue shaded area is the sublevel set of $f(\bx,\bc)$ at $\bx0$. The red dashed area is the inferred feasible set.}
    \label{fig:halfspaceintuition}
\end{figure}

% \kg{Do we assume all imputed sets are convex?} \hm{Yes, based on definition 5.}
\begin{proposition} \label{prop:imputedinC}
%A convex nominal set $\cS = \cX \cap \cXt$ is an imputed set for FO if and only if $\cS \subseteq \cC$. 
If a $\cS = \cX \cap \cXt$ is an imputed set for FO then $\cS \subseteq \cC$. 
\end{proposition}
Figure~\ref{fig:halfspaceintuition} shows the intuition behind Proposition~\ref{prop:imputedinC} which states any imputed set for FO must be contained within the tangent half-space $\cC$. %In Figure~\ref{fig:halfspaceintuition}, 
Consider the convex nominal set $\cS$ denoted by the dashed green area, which is not contained within $\cC$. Then $\cS$ must contain points outside of $\cC$ that are either in $\cV \setminus \cC$ or in the complement of $\cV$. If $\cS$ contains $\hbx' \not \in \cV$, then it cannot be an imputed set, by definition, since $\hbx'$ dominates $\bx^0$ and hence $\bx^0 \not \in \underset{\bx \in \cS}{\argmin}\,\{f(\bx; \bc) \}$. If $\cS$ contains a point $\hbx \in \cV\setminus\cC$, then because $f(\bx;\bc)$ is convex, there exists a convex combination of $\bx^0$ and $\hbx$ that dominates $\bx^0$ and hence, again, $\bx^0 \not \in \underset{\bx \in \cS}{\argmin}\,\{f(\bx; \bc) \}$, which contradicts the definition of an imputed set.  
Proposition~\ref{prop:imputedinC} %demonstrates an important property of imputed sets. It also 
provides a guideline for constructing an imputed set %that can be constructed 
by intersecting any convex nominal set with the tangent half-space $\cC$, as outlined in Proposition~\ref{prop:CcapS}. % states this result.  

\begin{proposition}\label{prop:CcapS}
For any convex nominal set $\cS$, the set $\cS \cap \cC$ is an imputed set. 
\end{proposition}
% \begin{proof}{\textsc{Proof of Proposition~\ref{prop:CcapS}}.}
% Both $\cS$ and $\cC$ are convex so $\cC \cap \cS$ is also convex. Because $\cS$ is a nominal set and $\cC$ includes all accepted observations, $\cC\cap\cS$ is also a convex nominal set. Finally, Since $\cS$ is a convex nominal set, it includes the convex hull of all observations, and $\bx^0 \in \cS$, and therefore, $\bx^0\in\cS\cap\cC$. Because $\cC$ is the tangent half-space to the sublevel set of $f(\bx;\bc)$ at $\bx^0$, it is also given that $\bx^0 \in \underset{\bx \in \cS\cap \cC}\argmin\{f(\bx;\bc)\}$. Therefore, $\cS \cap \cC$ meets the criteria outlined in Definition~\ref{def:imputed} and is therefore an imputed set. 
% %\hm{show by contradiction.}
% \end{proof}

Proposition~\ref{prop:CcapS} illustrates that the optimality condition on $\bx^0$ can always be guaranteed if $\cC$ is used as one of the constraints in shaping the imputed set. Because any imputed set must be a subset of the tangent half-space $\cC$, as shown in Proposition~\ref{prop:imputedinC}, the addition of $\cC$ does not exclude or cut any possible imputed sets for FO. Hence, instead of searching for a nominal set that satisfies the optimality on $\bx^0$, we can add the tangent half-space $\cC$ as one of the known constraints for FO, thereby, always guaranteeing that $\bx^0$ will be the optimal solution of FO for any inferred feasible set. This additional constraint allows us to relax the KKT conditions on the optimality of $\bx^0$ in the \changemarker{GIO} formulation and derive a reduced formulation, presented in Theorem~\ref{theorem:RDIO}. 

\begin{theorem}\label{theorem:RDIO}
If $\cC$ is appended to the known constraints of FO, solving \changemarker{GIO} is equivalent to solving the following reduced model:
% , and $\bx^0 \in \ext(\cS \cap \cC)$.
%constraints~\eqref{eq:IOstationarity}--\eqref{eq:IOcsL} are redundant in  \changemarker{GIO}. 
% If $\cC$ is appended to the known constraints of FO, then constraints~\eqref{eq:IOstationarity}--\eqref{eq:IOcsL} can be removed from \changemarker{GIO} without affecting the imputed set. 
\begin{subequations}
\begin{align} 
%hspace{-1cm}
\mathbf{\changemarker{RGIO}}: \quad \underset{\ba, b, \bq, \by}{\text{Maximize}} \quad  &   
%\underset{k \in \cKbad}{\max}(\bP ([\bA, \bb]_\text{inferred},\bx^{k})) \label{eq:eIOobj} \\ 
%\bP \left([\bA, \bb]_\text{inferred},\left\{\bx^{k},\, k\in \cK\right\}\right) \label{eq:eIOobj} \\ 
%\bP \left(\cXt,\left\{\bx^{k},\, k\in \cK\right\}\right) 
\bP \left(\bq_1, \dots, \bq_{|\cNt|}, \bA, \bb; (\bx^{1},\dots \bx^{|\cK|}) \right) \label{eq:eIOobj} \\
 \text{subject to} \quad & g_n(\bx^k; \,\bq_n) \ge \bzero, \quad \forall k\in \cKgood , n\in \cNt 
 \label{eq:eIOprimalfeasNL} \\ 
&   \ba_\ell'\, \bx^k \ge b_\ell,   \qquad \forall k\in \cKgood, \ell \in \cLt \label{eq:eIOprimalfeasL}  \\
& \ba_{\ell}\bx^k \le b_{\ell} - \epsilon + M y_{\ell k}, \quad \forall \ell \in \cL \cup \cLt, \,  k \in \cKbad \label{eq:eIOinfeasL}\\
& g_n(\bx^k; \,\bq_n) \le \bzero - \epsilon + M y_{nk}, \quad \forall n \in \cN \cup \cNt, \, k \in \cKbad \label{eq:eIOinfeasNL}\\
& \sum_{i \in \cI } y_{ik} \leq \, \mid \cI \mid - 1, \quad \forall k \in \cKbad \label{eq:eIOinfeasSum}\\
& \ba_{\ell} \in \cA_{\ell} \,, b_{\ell} \in \cB_{\ell}, \quad \forall \ell \in \cLt \qquad \label{eq:eIOnormL} \\
& \bq_n \in \cQ_n \,, \quad \forall n \in \cNt  \label{eq:eIOnormNL} \\ %||\bq_n||= 1
& y_{ik} \in \{0,1\}, \qquad \forall i \in \cI, \, k \in \cKbad. \label{eq:eIOBinaryVars}
\end{align}
\end{subequations}
\end{theorem}

We note that \changemarker{RGIO} is always feasible because it is a relaxed version of \changemarker{GIO} with fewer constraints, and we know from Proposition~\ref{prop:DIOfeas} that \changemarker{GIO} is always feasible. Based on Theorem~\ref{theorem:RDIO}, to find an imputed set, we can simply find a nominal feasible set and super-impose the known constraints including the tangent halfspace of the sublevel set of $f(\bx;\bc)$ at the preferred solution~$\bx^0$.

\begin{corollary}\label{cor:theorem}
\changemarker{RGIO} %finds constraint parameters $\ba_{\ell}, b_\ell,\, \forall \ell \in \cLt$, and $\bq_n, \, \forall n \in \cNt$ such that the inferred set $\cXt$ %\[\cX = \left\{\,\bx\in \mathbb{R}^n \, \middle\vert \, g_n(\bx; \,\bq_n) \ge \bzero, \, \forall n \in \cNt, \,\,\, \ba_{\ell}' \bx \ge b_{\ell}, \, \forall \ell \in \cLt \, \right\}\] 
infers unknown constraints $\cXt$ that
shape a nominal feasible set for FO such that $\bx^0$ is an optimal solution of 
 \begin{subequations}
\begin{align}   \label{eq:theorem}
\underset{\bx}{\text{minimize}} \quad &f(\bx; \bc)\\ 
 \text{subject to} \quad & \bx \in \cC \cap \cX \cap \cXt.
\end{align}
\end{subequations}
% \hm{Be careful with the wording of ``nominal'' and make sure we are consistent with the inclusion or exclusion of $\cap\cC$.}
% eDIO: Replace two of three KKT conditions with the halfspace constraint instead. 
\end{corollary}
% \begin{proof}{\textsc{Proof of Corollary~\ref{cor:theorem}.}}
% Due to constraints~\eqref{eq:eIOprimalfeasNL}--\eqref{eq:eIOBinaryVars}, we know that $\cXt$ is a nominal feasible set for FO. 
% %
% Note that $\bx^0 \in \cX \cap \cXt \cap \cC$ given that $\bx^0 \in \cX$ by Assumption~\ref{assumption:ObsAreFeas}, $\bx^0 \in \cXt$ because $\cXt$ is a nominal set, and $\bx^0 \in \cC$ by definition. Furthermore, $\bx^0 \in \cC = \{\bx \mid f(\bx;\bc) \geq f(\bx^0; \bc) \}$, which implies that $\bx^0 \in \underset{x\in\cC\cap\cX\cap\cXt}{\argmin} f(\bx;\bc)$ and is therefore an optimal solution of~\eqref{eq:theorem}.  
% \hfill \halmos
% \end{proof}

Corollary~\ref{theorem:RDIO} provides a method for reducing the complexity of the \changemarker{GIO} model. In Section~\ref{sec:numEx}, we provide a numerical example that illustrates how an imputed set can be constructed using Corollary~\ref{cor:theorem}. %In this example, the solution is constructed by manual inspection.  

The \changemarker{RGIO} model eliminates the need for explicitly writing the stationarity and complementary slackness constraints because the inclusion of the tangent half-space $\cC$ makes them redundant for the \changemarker{GIO} model. Because $\cC$ is a tangent half-space of the sublevel set of $f(\bx;\bc)$ at the preferred solution $\bx^0$, its inclusion ensures that the resulting inferred feasible region is an imputed set. The only constraints that are required to remain in the \changemarker{RGIO} model are those that ensure the imputed constraints form a nominal solution that includes all accepted observations and none of the rejected ones. Therefore, the size of the \changemarker{RGIO} problem is considerably lower than that of the \changemarker{GIO} problem, and it relaxes a large number of nonconvex nonlinear constraints. A comparison of the number of variables and constraints in the \changemarker{GIO} and \changemarker{RGIO} models is provided in Table~\ref{table:size}. 
\begin{table}[htbp]
\centering
    \begin{tabular}{l @{\hspace{.5cm}} l l c }
     \toprule 
    % & & & General & Example  \\ \midrule
    & Type & Model & Count \\ \midrule
    \multirow{4}{*}{\rotatebox[origin=c]{90}{Variables \quad } } & \multirow{2}{*}{Continuous }  & \changemarker{GIO} & $(m+2)|\cLt|+(\phi+1)|\cNt| + |\cN|+|\cL|$   \\ \cmidrule{3-4}
    & & \changemarker{RGIO} & $(m+1)|\cLt|+\phi|\cNt|$  \\ \cmidrule{2-4}
    & \multirow{2}{*}{Binary} & \changemarker{GIO} & $(|\cN|+|\cL|+|\cNt|+|\cLt|)\,|\cKbad|$  \\ \cmidrule{3-4}  
    & & \changemarker{RGIO} & $(|\cN|+|\cL|+|\cNt|+|\cLt|)\,|\cKbad|$ 
        %same as \changemarker{GIO} 
        \\ \bottomrule [1pt]
    \multirow{4}{*}{\rotatebox[origin=c]{90}{Constraints \quad} } & \multirow{2}{*}{Linear } & \changemarker{GIO} & $|\cLt|(|\cKgood|+1) + (|\cL|+1)|\cKbad|$ \\ \cmidrule{3-4}
    & & \changemarker{RGIO} & $|\cLt|(|\cKgood|+1) + (|\cL|+1)|\cKbad|$  
    %same as \changemarker{GIO} 
    \\ \cmidrule{2-4}
    & \multirow{2}{*}{Nonlinear} & \changemarker{GIO} & $|\cNt|(|\cKgood|+|\cKbad| + 1) + |\cN|(|\cKbad|+1) + m $  \\ \cmidrule{3-4}  
    & & \changemarker{RGIO} & $|\cNt|(|\cKgood|+|\cKbad|) + |\cN||\cKbad| $  \\ \bottomrule
        %\multirow{2}{*}{\small Constraints} &  {\small Linear}  & 10 & 10 & &  \\ \cmidrule{2-6} 
    %& {\small Nonlinear} & & & & \\ \bottomrule
    \end{tabular}
    \caption{Comparison of the number of variables and constraints in the \changemarker{GIO} and \changemarker{RGIO} models.}
    \label{table:size}
\end{table}

Table~\ref{table:size} shows that the \changemarker{RGIO} model has fewer continuous variables and nonlinear constraints compared to the \changemarker{GIO} model. The nonlinear constraints that remain in \changemarker{RGIO} primal feasibility constraints (i.e., $g_n(\bx;\bq)\geq \bzero$) which are nonlinear in $\bx$, but may be linear (or be linearized) in $\bq_n$, which is the decision variable in \changemarker{RGIO}.  
%of the same complexity of the constraints in FO. 
%Hence, the \changemarker{RGIO} retains the complexity of the FO model, notwithstanding the addition of binary variables that ensure the rejected observations are outside of the inferred feasible region.  
Particularly, if the FO is linear, then the corresponding \changemarker{RGIO} model can be linear \changemarker{integer} if a linear distance metric is used in the objective function. In what follows, we discuss one example of a distance metric that can be used in the \changemarker{RGIO} model.

\subsection{Example of Distance Metric} \label{sec:distance}
The objective function in \changemarker{GIO} and \changemarker{RGIO} maximizes a non-negative distance metric between the inferred constraints and the observations. In this section, we provide an example of such a distance metric and write the complete \changemarker{GIO} formulation based on it, which will be used in the application example presented in Section~\ref{sec:case}.

We consider an objective that aims to \changemarker{robustify against inclusion of rejected observations in the inferred feasible region by finding} constraints that are as far as possible from the rejected observations and as close as possible to the accepted ones. To formulate this objective, we use Separation Distance, defined as 
%We consider the objective function as
maximizing the maximum distance between the inferred constraints to exclude each of the rejected observations, $\bx^k, k\in \cK^-$, from the inferred feasible set. That is, among the constraints that make each rejected observation infeasible, we select the constraint furthest away from the rejected observation and push it as close to the accepted observations as possible. Hence, we define the Separation Distance  $\bP$ as \[ %\bP \left(\cXt,\left\{\bx^{k},\, k\in \cK\right\}\right) 
\bP \left(\bq_1, \dots, \bq_{|\cNt|}, \bA, \bb; (\bx^{1},\dots \bx^{|\cK|}) \right) =  \sum_{k \in \cKbad} 
\max \left\{ {\underset{\ell \in \cLt}{\max}\left\{d_{\ell k}([\ba_{\ell}, b_{\ell}],\bx^{k})\right\}, \, \max \left\{ d_{nk}(\bq_n,\bx^k)\right\}} \right\},
 \]
%\mathbf{DIO}: \quad \underset{\ba, b, \bq, \lambda, \mu,y}{\text{Maximize}} \quad  &  \bP \left(\bq_1, \dots, \bq_{|\cNt|}, \bA, \bb; (\bx^{1},\dots \bx^{|\cK|}) \right) \label{eq:IOobj} 
% \hm{ NEXT TIME: The above is correct where we calculate the distance from linear and nonlinear separately and take the maximum as the right distance to be maximized in the obj. Below, separate linear and nonlinear indices to represent what we have above. }
where $d_{ik}$ is the distance between the imputed constraint $i \in \cLt \cup \cNt$ and the rejected observation $k \in \cKbad$. This objective can be linearized using a set of additional constraints as shown in model \eqref{model:RDIO_SeparationMetric}. 
% subject to distbad_1{i in I1, k in Kbad}: slack[i,k] >= b[i] - sum {j in J} a[i,j]*xbad[k,j]; 
%subject to nonnegslack{i in I1, k in Kbad}: slack[i,k]>=0;
% subject to distbad_2{i in I1, k in Kbad}: slack[i,k]<= b[i] - sum {j in J} a[i,j]*xbad[k,j] + M*bin1[i,k];
% subject to distbad_3{i in I1, k in Kbad}: slack[i,k]<= M*(1-bin1[i,k]); 
%
\begin{subequations}
\begin{align}
%\mathbf{RDIO_{1}}: 
\quad \underset{\ba, b, \bp, y, z}{\text{Maximize}} \quad  &  \sum_{k \in \cKbad} z_k \label{eq:5a} \\ 
 \text{subject to} \quad & d_{n k} \geq -g_n(\bx^k;\bq_n), \qquad \forall n \in \cNt, k \in \cKbad \label{eq:5b} \\ 
 & d_{n k} \leq -g_n(\bx^k;\bq_n) + M y_{nk}, \qquad \forall n \in \cNt, k \in \cKbad \\
%  & d_{\ell k} \geq b_{\ell} - \sum_{m\in M} a_{\ell m} x_{km} \qquad \forall \ell \in \cLt, k \in \cKbad \\
  & d_{\ell k} \geq b_{\ell} - \ba_{\ell}' \bx^{k}, \qquad\quad  \forall \ell \in \cLt, k \in \cKbad \\
%  & d_{\ell k} \leq b_i - \sum_{m \in M} a_{\ell m} x_{km} + M r_{\ell k} \qquad \forall \ell \in \cL, k \in \cKbad \\
%  & d_{\ell k} \leq b_{\ell} - \sum_{m \in M} a_{\ell m} x_{km} + M y_{\ell k} \qquad \forall \ell \in \cLt, k \in \cKbad \\
  & d_{\ell k} \leq b_{\ell} - \ba_{\ell}' \bx^{k} + M y_{\ell k}, \qquad \forall \ell \in \cLt, k \in \cKbad \\
%  & d_{\ell k} \leq M (1-r_{\ell k})\\ %same as the one below
 & d_{ik}  \leq M(1-y_{ik}), \qquad \forall i \in \cNt \cup \cLt, k \in \cKbad \\ 
 & d_{ik} \geq \epsilon (1-y_{ik}),  \qquad \forall i \in \cNt \cup \cLt, k \in \cKbad  \label{eq:5g} \\
 & z_k \leq d_{ik} + M p_{ik}, \qquad \forall i \in \cNt \cup \cLt, k \in \cKbad \label{eq:5h}\\
 & \sum_{i \in \cI} p_{ik} \leq |\cI|-1, \qquad \forall k \in \cKbad  \\
 %& p_{ik} \leq |\cI|-1 \qquad \forall i \in \cNt \cup \cLt, k \in \cKbad\\
 & p_{ik} \geq y_{ik}, \qquad \forall i \in \cNt \cup \cLt, k \in \cKbad  \label{eq:5j}\\
&  \eqref{eq:eIOprimalfeasNL}-\eqref{eq:eIOBinaryVars}. %\\
 %& d_{i k} \geq 0, \qquad \forall i \in \cNt \cup \cLt, k \in \cKbad.  
\end{align}  \label{model:RDIO_SeparationMetric}
\end{subequations}
Constraints~\eqref{eq:5b}--\eqref{eq:5g} find the slack distance between each rejected point and each constraint. Constraints~\eqref{eq:5h}--\eqref{eq:5j} find the maximum distance of each constraint that makes each point infeasible, and the objective~\eqref{eq:5a} maximizes this maximum distance.  The original primal feasibility constraints of the \changemarker{RGIO} model are also enforced. Note that this specific distance model only requires the addition of linear constraints and binary variables, so if the FO is a linear problem, the corresponding \changemarker{RGIO} problem will be a linear integer problem. Any other metric of interest can also be used in the objective, but for simplicity, we will use this linear metric in the application example in Section~\ref{sec:case}. 

% Generalization: if other objectives are known, they can be added as a multiple objective to this model (see other paper)

% \kg{Might need more text here too before closing the methods section.}

\subsection{Illustrating Numerical Example}\label{sec:numEx}
\changemarker{In this section, we provide an illustrative two-dimensional FO problem and formulate the corresponding RGIO formulation. We formulate a nonlinear FO model and its corresponding nonlinear mixed-integer RGIO model in AMPL Version: 3.1.1 and solve them both using Gurobi 11.0.3. All results are rounded to one decimal place.}
% \vspace{0.5em}
\begin{example} \label{ex:numerical}
 %%%%%%%%%%%%%%%%% New Nonlinear example - July 2024
\changemarker{Consider the following nonlinear convex FO problem 
% \kg{we need to update the example}
%
\begin{subequations}
\begin{align}  
\underset{\bx}{\text{minimize}} \quad & 2^{x_1} +x_2  \label{eq:objNL}\\ 
 \text{subject to} \quad & \frac{(x_1-q_1)^2}{4}+\frac{(x_2-q_2)^2}{2} \leq q_3 ,  \label{eq:ex:nonlin}\\
%\text{subject to} \quad & \frac{(x_1-3)^2}{4}+\frac{(x_2-2)^2}{2} \leq 1 \\
%%  & x_1-2x_2 \leq 1 
  & a_1x_1+a_2x_2 \geq b , \label{eq:ex:lin} \\
  & x_1, x_2 \geq 0. \label{eq:ex:sign}
\end{align}
\end{subequations}
where $x_1, x_2$ are the decision variables, $a_1$, $a_2$, and $b$ are unknown parameters of the linear constraint~\eqref{eq:ex:lin}, and $q_1, q_2,$~and $q_3$ are the unknown parameters of the nonlinear convex constraint~\eqref{eq:ex:nonlin}. The non-negativity constraints~\eqref{eq:ex:sign} are the only known constraints. %Note that the objective function~\eqref{eq:objNL} is convex and monotone. 
\\
Let $\cKgood=\{1,\dots,13\}$ and $\cKbad=\{14,\dots,21\}$ indicate 21~total observations, with 13~accepted  observations 
\setcounter{MaxMatrixCols}{20}
$\begin{pmatrix} 
    \bx^1,  \dots, \bx^{13} 
\end{pmatrix}
=\begin{pmatrix}
1.5 & 2 & 2.5 &  3&  4 & 4 & 1.5 &3 &  3   & 3   & 4.8 & 2 & 5 \\
1.5 & 2 & 1   & 3 & 2  & 3 & 2.5 & 1.5 & 3.4 &  2.5 & 2.4 & 3 &2 \\
\end{pmatrix}$ 
and 7~rejected observations %as the accepted observations and
$\begin{pmatrix}
    \bx^{14}\dots \bx^{21}\\
\end{pmatrix}
=\begin{pmatrix}
4 & 4 & 1 & 3 & 1 & 1.5 &5  \\
1 & 3.5 & 3 & 4 & 1 & 3.5 &1 \\
\end{pmatrix}$. 
}
% the black and red points, respectively, as shown in Figure~\ref{fig:NLExample} and equation \eqref{??}. 
% %
% The preferred solution can be identified based on Definition~\ref{def:x0} and is shown in the figure as $\bx^0$. The sublevel set of $f(\bx;\bc)$ at $\bx^0$ is denoted as $\cV$ and can be written as $\cV= \{ \bx \in \mathbb{R}^2 \mid 2^{x_1} +x_2 \geq 4\}$. %which is the outside of the dashed circle passing through $\bx^0$. 
% Assume that a known constraint is $x_1-2x_2 \leq 1$ for which the corresponding side of the inequality is shown by $\cX$ in Figure~\ref{fig:NLExample}. }
%
 %The following example is a convex nonlinear optimization. The objective function $f(\bx;c)=2^{x_1} +x_2$ is convex and monotonically increasing with gradient of $\nabla=[2^{x_1}ln(2), 1] \geq 0$ and Hessian of $H=[2^{x_1}ln(2)^2, 0; 0, 0] \geq 0$. 
 %
\changemarker{
The preferred observation is $\bx^0 = \begin{pmatrix}
    1.5 \\ 1.5
\end{pmatrix}$ which leads to the $\bx^0$-\textit{sublevel set} of $\cV=\{(x_1,x_2) \mid  2^{x_1} +x_2 \geq 2^{1.5}+1.5\}$. Since the objective function~\eqref{eq:objNL} is monotone, convex, and differentiable in $x_1, x_2$, we can write the tangent half-space as 
\begin{align*}
&\cC = \{\bx \in \mathbb{R}^n  \mid (\nabla f(\bx^0;\bc))' \,\bx \geq \nabla f(\bx^0; \bc)'\bx^0 \}  \\
% & \cC = [2 \ln(2),1] [x_1, x_2]' \geq [2\ln(2),1] [1, 2]' \Rightarrow \\
\Rightarrow \quad & \cC = \begin{pmatrix}
    2^{1.5}\ln 2 \\ 
    1
\end{pmatrix}' 
\begin{pmatrix}
    x_1 \\ 
    x_2
\end{pmatrix}
 \geq
 \begin{pmatrix}
     2^{1.5} \ln 2 \\ 
     1
 \end{pmatrix}'
 \begin{pmatrix}
     1.5 \\
     1.5
 \end{pmatrix} \\ 
\Rightarrow \quad  & \cC = (2^{1.5}\ln 2) x_1+x_2 \geq (2^{1.5}\ln 2)(1.5) +1.5.
\end{align*}
%
%The part of IO that's missing. Assume $g_n(\bx^k; \bq_n)\geq 0$, hence, $q - \frac{(x_1-3)^2}{5}-\frac{(x_2-2)^2}{2} \geq 0$, where $q$ is the original RHS that in our example was ``1". Then we should have 
Considering this $\cC$ as a known constraint in FO, the RGIO model can be written as 
\begin{subequations}
\label{eq:nonlinearEx2}
\begin{align} 
%\mathbf{RGIO}: 
% \quad \underset{\ba, b, \bq, y}{\text{Maximize}} \quad  & \bP \left(\bq_1, \dots, \bq_{|\cNt|}, \bA, \bb; (\bx^{1},\dots \bx^{|\cK|}) \right) \label{eq:eIOobjEx} \\
\quad \underset{\ba, b, \bq, \by}{\text{Maximize}} \quad  & \bP \left(q_1, \dots, q_3, \, \ba, b; (\bx^{1},\dots \bx^{|\cK|}) \right) \label{eq:eIOobjEx} \\
% \text{subject to} \quad & (6b)--(6d) \\
 \text{subject to} \quad & q_3 - \frac{(x_1^k-q_1)^2}{4}-\frac{(x_2^k-q_2)^2}{2} \ge \bzero, \quad \forall k\in \cKgood\\ 
& q_3 - \frac{(x_1^k-q_1)^2}{4}-\frac{(x_2^k-q_2)^2}{2} \le \bzero - \epsilon + M {y}_{1k}, \quad \forall k \in \cKbad \label{eq:eIOinfeasNLEx}\\
%& \frac{(x_1^k-q_1)^2}{4}+\frac{(x_2^k-q_2)^2}{2} \le 1 - \epsilon + M y_{\ell k}, \quad \forall \ell \in \cN \cup \cNt, \,  k \in \cKbad \label{eq:eIOinfeasLEx}\\
&   \ba'\bx^k \ge b,   \qquad \forall k\in \cKgood \label{eq:eIOprimalfeasLEx}  \\
& \ba'\bx^k \le b - \epsilon + M y_{2 k}, \quad \forall  k \in \cKbad \label{eq:eIOinfeasLEx}\\
& y_{1k}+ y_{2k} \leq \, 1, \quad \forall k \in \cKbad \label{eq:eIOinfeasSumEx}\\
& y_{1k}, y_{2k} \in \{0,1\}, \qquad \forall k \in \cKbad. \label{eq:eIOBinaryVarsEx} \\
& \text{other existing constraints.}
\end{align}
\end{subequations}
% \begin{subequations}
% \label{eq:nonlinearEx2}
% \begin{align} 
% %\mathbf{RGIO}: 
% \quad \underset{\ba, b, \bq, y}{\text{Maximize}} \quad  &   
% \bP \left(\bq_1, \dots, \bq_{|\cNt|}, \bA, \bb; (\bx^{1},\dots \bx^{|\cK|}) \right) \label{eq:eIOobjEx} \\
% % \text{subject to} \quad & (6b)--(6d) \\
%  \text{subject to} \quad & q_3 - \frac{(x_1^k-q_1)^2}{4}-\frac{(x_2^k-q_2)^2}{2} \ge \bzero, \quad \forall k\in \cKgood , n\in \cNt \\ 
% & q_3 - \frac{(x_1^k-q_1)^2}{4}-\frac{(x_2^k-q_2)^2}{2} \le \bzero - \epsilon + M {y}_{nk}, \quad \forall n \in \cN \cup \cNt, \, k \in \cKbad \label{eq:eIOinfeasNLEx}\\
% %& \frac{(x_1^k-q_1)^2}{4}+\frac{(x_2^k-q_2)^2}{2} \le 1 - \epsilon + M y_{\ell k}, \quad \forall \ell \in \cN \cup \cNt, \,  k \in \cKbad \label{eq:eIOinfeasLEx}\\
% &   \ba_\ell'\, \bx^k \ge b_\ell,   \qquad \forall k\in \cKgood, \ell \in \cLt \label{eq:eIOprimalfeasLEx}  \\
% & \ba_{\ell}\bx^k \le b_{\ell} - \epsilon + M y_{\ell k}, \quad \forall \ell \in \cL \cup \cLt, \,  k \in \cKbad \label{eq:eIOinfeasLEx}\\
% & \sum_{i \in \cI } y_{ik} \leq \, \mid \cI \mid - 1, \quad \forall k \in \cKbad \label{eq:eIOinfeasSumEx}\\
% & y_{ik} \in \{0,1\}, \qquad \forall i \in \cI, \, k \in \cKbad. \label{eq:eIOBinaryVarsEx} \\
% & \text{other existing constraints}
% \end{align}
% \end{subequations}
}
% \kg{Update the rest of the example from here (August 28)\\}
\noindent\changemarker{Solving Formulation~\eqref{eq:nonlinearEx2} using the Separation Distance metric $\bP$ as defined in Section~\ref{sec:distance} yields the optimal RGIO solution $(q_1,q_2,q_3)=(3.0,2.0,1.0)$ and $(a_1,a_2)=(-1,2), b=-1$.} % to form a linear constraint. Note that in this case, 

\changemarker{Figure~\ref{fig:NLExample2} illustrates the input and output of the inverse model graphically. The feasible and infeasible observations are depicted as black and red dots, respectively. The inferred nonlinear and linear constraints are shown in the green ellipse and green linear half-space, respectively. The linear constraint makes the rejected point $(4,1)$ infeasible, and the inferred ellipse makes the rest of the rejected observations (red) infeasible. 
%It can be seen that all acceptable observations (black) are feasible with respect to the known and inferred constraints. 
The tangent half-space, shown as a solid black line marked by $\cC$, which is considered a known constraint, enables the preferred solution, $\bx^0$ to be optimal. The sublevel set of the objective function at $\bx^0$ is shown in a solid blue curve and two other sample isocost curves are illustrated in dashed blue curves.}
\begin{figure} 
    \centering
\begin{tikzpicture}[scale=1.2]
\node (img) {\includegraphics[width=0.7\linewidth]{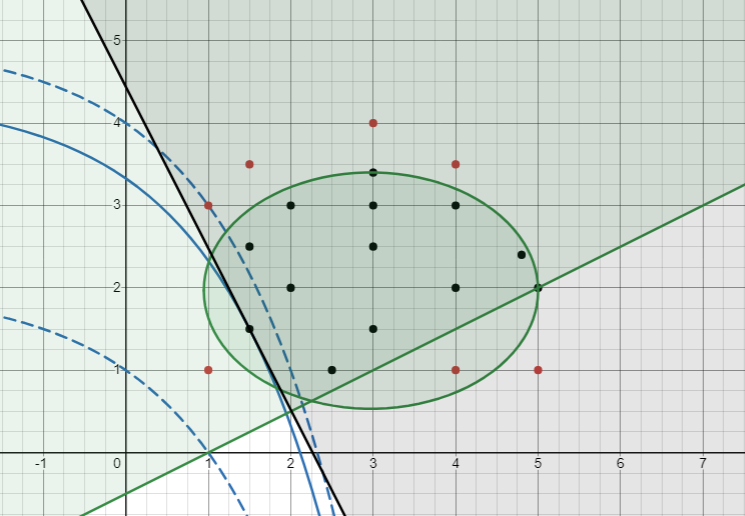}};
% \draw [->, black] (3.9,.5) -- (3.6,1);
\node [black] at (-3.2,3.2){$\cC$}; 
\node [black] at (-1.86,-0.9){$\bx^0$};
\draw [->, black] (-3.66,3.2) -- (-3,3.5);
\node [blue] at (-5.7,2.3){$2^{x_1}+x_2=5$}; %dashed
\node [blue] at (-4.8,0.9){$2^{x_1}+x_2=2^{1.5}+1.5$}; %solid
\node [blue] at (-4.8,-1.3){$2^{x_1}+x_2=2$}; %dashed
\end{tikzpicture}
    \caption{Numerical example for inferring a linear and a nonlinear constraint}
    \label{fig:NLExample2}
\end{figure}

\changemarker{It can be seen that the intersection of the inferred or known constraints, i.e., $\cX \cap \cXt$,  only includes all accepted observations and excludes all the rejected observations. However, it does not have the preferred solution $\bx^0$ on its boundary, and as a result, would not allow $\bx^0$ to be optimal for FO. On the contrary, the tangent half-space $\cC$ allows for  $\bx^0$ to be a candidate optimal solution for FO when added as a known constraint. Hence, $\cXt \cap \cX \cap \cC$ is an imputed set for FO as it satisfies the properties outlined in Definition~\ref{def:imputed}, which in this example is constructed using Corollary~\ref{cor:theorem}.
}
% and consider the convex set $\cS = \{ (x-2)^2 + 2(y-2)^2 \leq 4.84\} $. Let $\bx^0=(0.6,.8)$ and note that $\bx^0 = \underset{\bx \in \cS}{\arg\min} \{f(\bx;\bc)\}$, and $f(\bx^0;\bc) = (1)^2 + (1)^2 = 2$. %It can be seen that $\cS \subseteq \cC$ in Figure~\ref{fig:sublevelset}. 
% \begin{figure}[htbp]
%     \centering
%     \includegraphics[width=0.5\textwidth]{desmos-graph.png} %https://www.desmos.com/calculator/r2taaskltu
%     \caption{Caption}
%     \label{fig:fig:sublevelset}
% \end{figure}
\hfill $\triangle$
 \end{example} 
\section{Application: Standardizing Clinical Radiation Therapy Guidelines} 
%\section{Numerical Example} 
\label{sec:case}
%-----------------------
% \hm{Can we say something about what the optimization after IO is different from the initial IO and how the last step involves personalization because each patient has its own Dij}
In this section, we test the proposed methodology using an application example on standardizing radiation therapy treatment planning guidelines for breast cancer patients. We first introduce the problem in the context of radiation therapy. Next, we describe the data and the experimental setup. Lastly, we present and discuss the results and provide practical insights.  
%\hm{remove mentions of "case study" and "data-driven" from the paper.}

\subsection{Problem Description: Standardizing Clinical Radiation Therapy Guidelines}
\begin{figure}[htbp]
%\begin{wrapfigure}{r}{0.6\textwidth}\vspace{-1em}
    \centering
\includegraphics[width=0.6\textwidth]{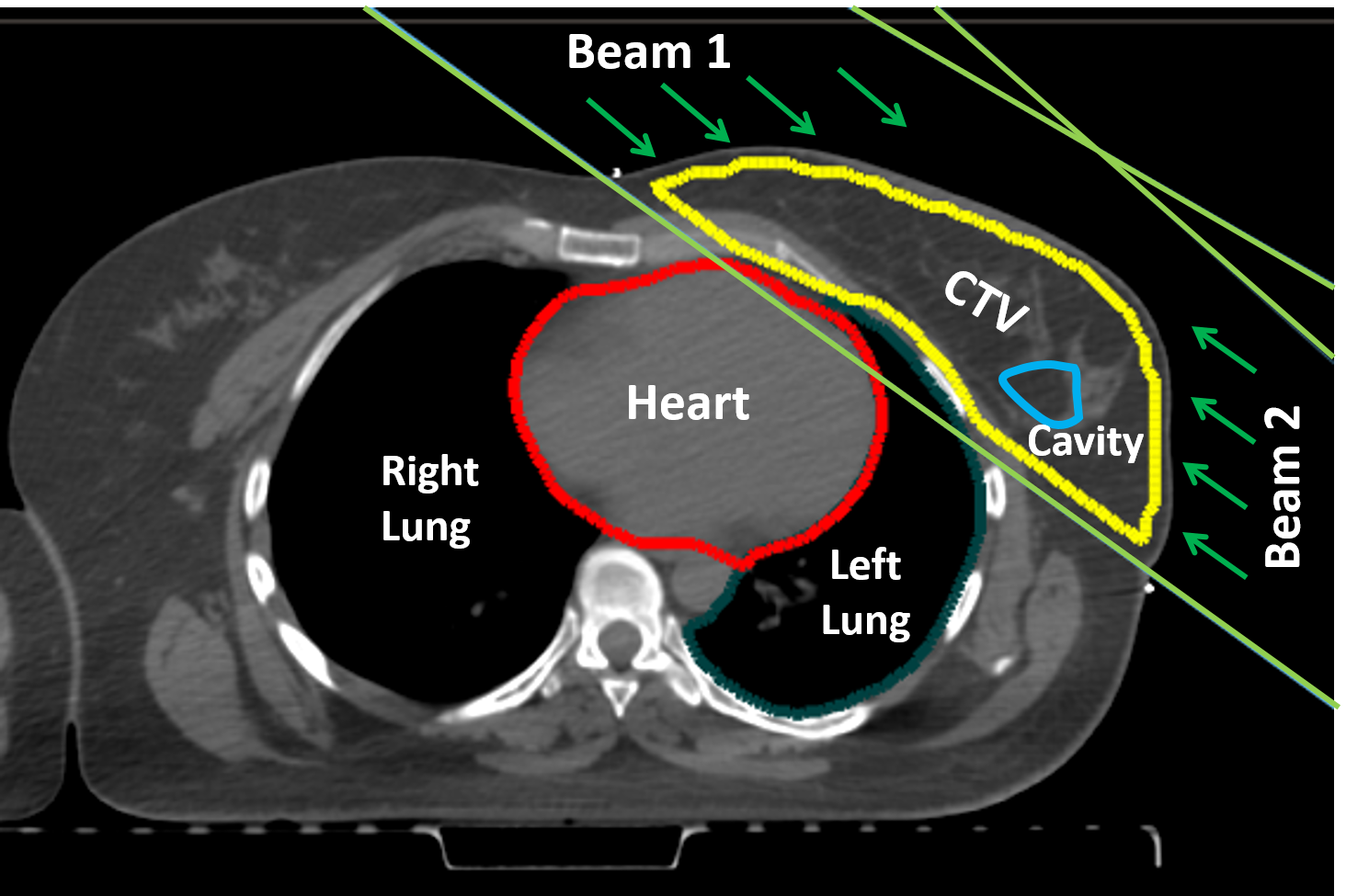}
    \caption{Computed Tomography (CT) scan of a breast cancer patient with important organs delineated. Image adapted from~\cite{mahmoudzadeh2015robust}.}
    \label{fig:CT}
    \vspace{-1em}
%\end{wrapfigure}
\end{figure}

Breast cancer is the most widely diagnosed type of cancer in women worldwide. The cancerous tumor is often removed, leaving behind a cavity, and radiation treatment is subsequently prescribed to eliminate any remaining cancer cells. Tangential intensity-modulated radiation therapy (IMRT) is a treatment modality often used as part of the treatment for most breast cancer patients. In tangential IMRT, two opposing beams that are tangent to the external body of the patient are used to deliver radiation to the breast tissue. The main organ at risk in breast cancer IMRT is the heart, particularly when the left breast is being irradiated \citep{mahmoudzadeh2015robust}. The goal is to find the radiation beam configurations from each angle such that the clinical target volume (CTV) inside the breast is fully irradiated and the heart is spared from radiation as much as possible. Figure~\ref{fig:CT} shows a computed tomography (CT) scan of a breast cancer patient along with the two tangential beams and the contours showing the organs at risk.

There is a set of acceptability guidelines in breast cancer treatment planning, which involve clinical dose-volume criteria on the CTV and the heart. A dose-volume criterion is a clinical metric that calculates the dose threshold to a certain fraction of an organ measured in units of radiation dose, Gray (Gy). For instance, the prescribed dose for the CTV is 42.4 Gy and at least 99\% of the CTV must receive lower than 95\% of the prescribed dose for a plan to be accepted. That is, if the body is discretized into three-dimensional cubes called voxels, then the 95\% quantile of the voxels must receive a dose of $42.4\times0.95=40.28$~Gy or higher. Similarly, at most 0.5\% of the CTV can receive a dose higher than 108\% of the prescribed dose, which means the dose to upper 0.5\% quantile of the CTV must be lower than $1.08 \times 42.4 = 45.792$. There are also upper bounds on the dose delivered to the heart, where the highest-dosed 10~cc and 25~cc volume of the heart must receive a dose lower than 90\% and 50\% of the prescribed dose, respectively. 

\begin{figure}[htbp]
    \centering
    \includegraphics[width=0.48\textwidth]{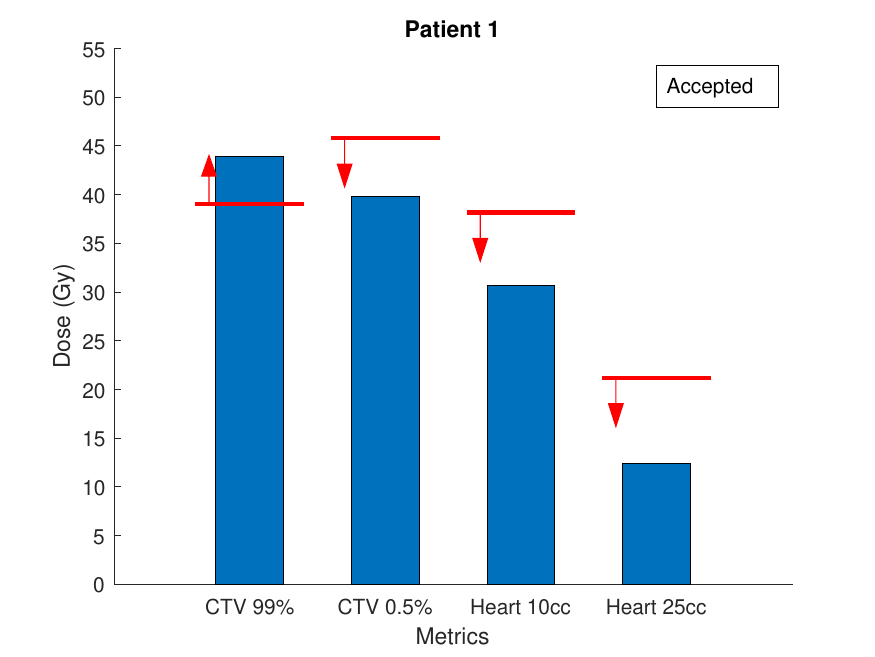}
    \includegraphics[width=0.48\textwidth]{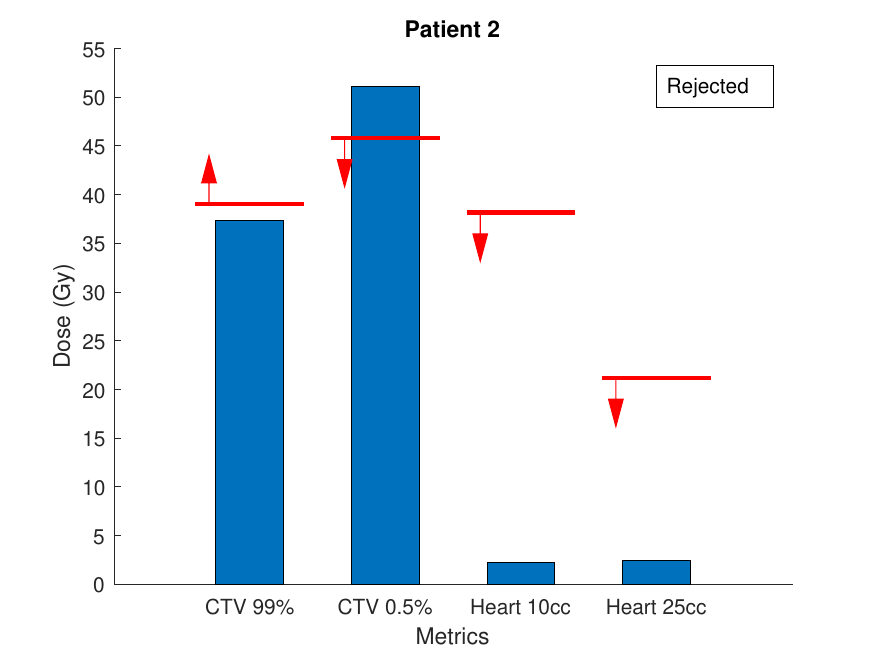}\\
    \vspace{2em}
    \includegraphics[width=0.48\textwidth]{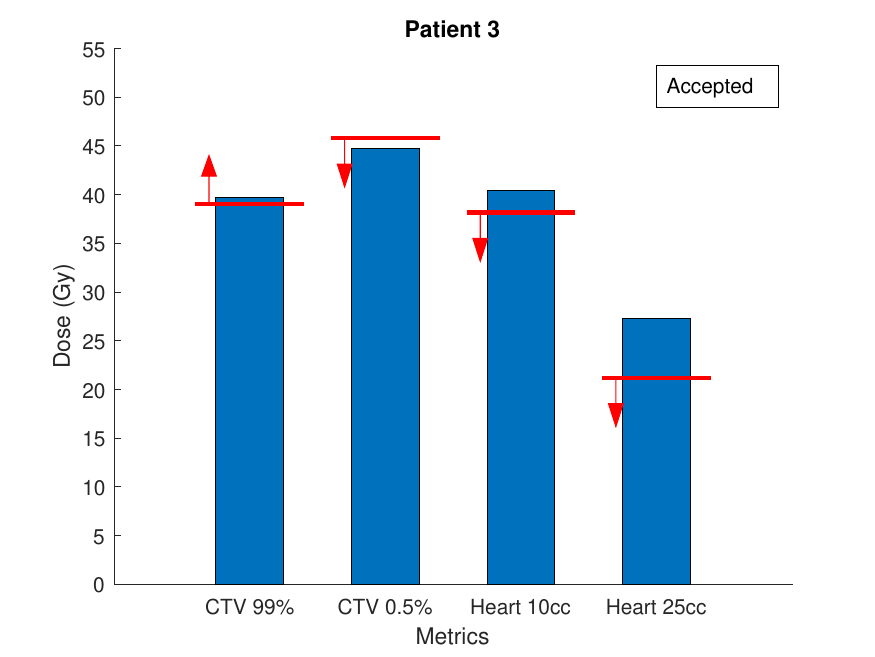}
    \includegraphics[width=0.48\textwidth]{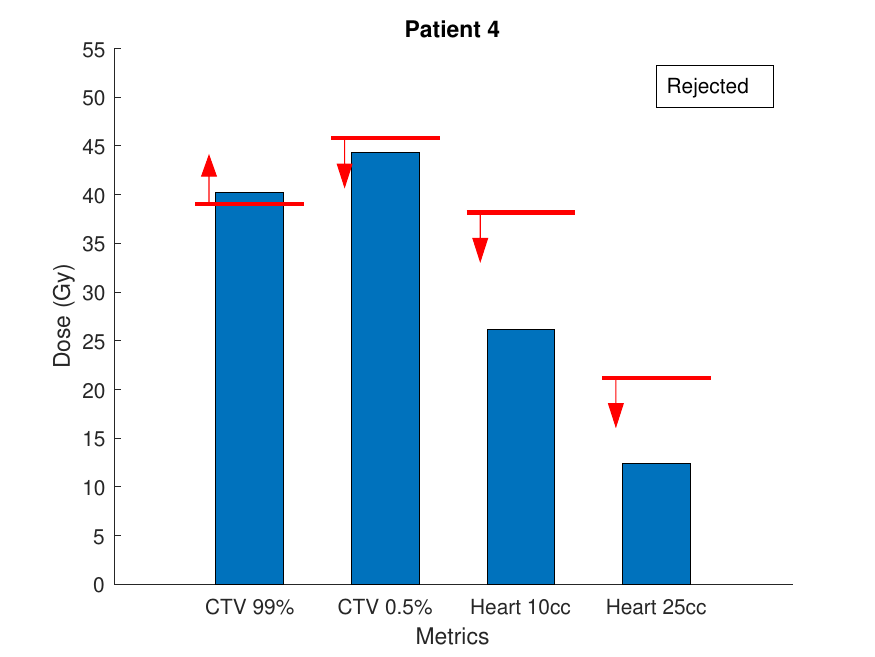}
    \caption{\changemarker{Hypothetical examples} of accepted/rejected treatment plans. The guidelines are indicated with (red) lines on the bars and arrows indicate the direction of bounds. }
    \label{fig:patientexamples}
\end{figure}

These clinical dose-volume guidelines are used as a base reference for planning, but acceptance or rejection of the plan is at the discretion of oncologists. % \changemarker{based on the practical clinical acceptance criteria for each patient's specific case}. 
%show guidelines and the first figure with 4 patients to describe the case. explain prescribed dose and DVH criteria for heart and CTV.   Add a paragraph or two to explain the current clinical guidelines. Do they explain all feasibility? Are there mislabels?
Figure~\ref{fig:patientexamples} shows an example of possible clinical guidelines and metric values for four \changemarker{hypothetical} %different 
patients that are labeled as accepted or rejected. It can be seen that the plan for Patient~1 and Patient~3 are accepted, while the plan for Patient~1 meets the guidelines but the plan for Patient~3 violates two of the dose-volume metrics on the heart.  On the contrary, the plan for Patient~2 is rejected because it does not meet the guidelines, but the plan for Patient~4 is also rejected even though it meets all the guidelines. \changemarker{These inconsistencies typically stem from discrepancies between current clinical guidelines-- which often are not data-driven, tailored to patient populations, or updated regularly -- and the practical criteria that oncologists use in their day-to-day decisions to accept or reject treatment plans. This ambiguity in true underlying criteria can lead to unnecessary back-and-forth between oncologists and planners to arrive at an acceptable plan. }
In what follows, we provide more details about the patient data and existing clinical guidelines in Section~\ref{sec:data}, and discuss the details of the corresponding FO and IO models for this problem in Section~\ref{sec:clinicalFOIO}.  % \hm{Maybe move this paragraph to beginning of Sec 5? }

% define FO and eIIO for the guidelines case (here or in the supplement). Note that, the FO here is not the inverse treatment planning optimization model that is frequently used in radiation planning. In contrast, our FO is a new optimization model to identify which treatment plans (perhaps outputs of the typical inverse planning) are accepted or not.  
% Explain FO. Assume that we have a set of features for each structure. Examples of the features can include min, max, mean dose to each structure and the structures can be left lung, CTV, and heart for a breast cancer patient. Explain model in words and write the detailed math model to the Appendix. 
%For this FO, we assume that we have a set of observations $\bx_k, k \in \cK$ which are known to be feasible for the FO problem and a set of observations $\bx^{BAD}_k, k \in \cKbad$ that are known to be infeasible for the FO. We write the following IO to infer the constraints. 
% Justify Linear because of computational complexity

\subsection{Patient Data and Clinical Guidelines}\label{sec:data}
% Explain that we have 5 patients, we are simulating 20 of each randomly, etc. One of the patients is a ``bad" case. How are we judging infeasibility of the generated (perturbed) patients. Currently this is based on $aCTV_{max} >= 45$ and $varCTV_{max} >= 44$. Which other current guidelines hold? Do we replicate any of the current guidelines with the inferred constraints or are they all new? Any insights? 
% Explain how we made the data. How did we decide which is feasible/infeasible. We first calculate a set of important metric from the available previously accepted dose for each patients. These metrics include minimum dose, maximum dose, average dose, dose-volume points (e.g. D99\% CTV dose). 

% \begin{itemize}
%             \item 100 breast cancer patient treatment plans 
%             \item Each plan labelled as either accepted or rejected
%             \item 16 dose metrics measured on each plan (min, max, dose-volume)
%             \begin{itemize}
%                 \item Targets: CTV, Cavity
%                 \item OARs: Lung, Heart
%             \end{itemize}
%\end{itemize}

% Data Generation
%\kg{2. Explain our data Generation\\}
We used retrospective clinical treatment plans for five breast cancer datasets to \changemarker{bootstrap and} simulate a population of an additional 100 plans, to a total of 105 patients. We perturbed the \changemarker{structure-based} radiation dose of each patient by 20\% in a uniformly random manner to create an additional 20 synthetic patients per original patient. The prescribed dose for all patients, synthetic or original, is 42.4~Gy. Among the five original plans, four were deemed clinically acceptable, and one was rejected (with the patient ultimately receiving an alternate mode of treatment under breath-hold \citep{wong1999use}). We assigned the same acceptable or rejected labels to the synthetic plans as the original plans that they are based on. %\footnote{The patient with the rejected plan was eventually moved to a different modality (under breath-hold) to receive their treatment.}. 
%Using this combination of acceptable and rejected patients, we perturbed the radiation dose of each patient by \%20 in a uniformly random manner to create an additional 20 synthetic patients per original patient, for a sum of 105 real and synthetic treatment plans. We assigned the same acceptable and rejected labels for the synthetic plans as the original plans that they are based on.  
%To identify acceptable and unacceptable treatments. Each of these plans was then labeled as accepted or rejected based on a \kg{set of expert-knowledge criteria} that did not match the guidelines. %These criteria were not fed into or used in any of the inverse learning steps. 

%\hm{This paragraph needs to be updated because we measure 14 metrics, but there are only 8 guidelines.} 
To inform the clinical treatment planning process, the plan for each of the patients is clinically measured through eight dose metrics, including maximum and minimum doses to different regions and dose-volume values for the CTV, the cavity, the heart, and the lung. These eight clinical guidelines and the accepted clinical limits are outlined in Table~\ref{tab:Guideline}. The last column illustrates the percentage of accepted plans that met each of the clinical criteria, %, for instance,  
which confirms that there is wide variability in how rigorously the guidelines are imposed on accepted plans. 
For instance, while the Heart max dose guidelines are met by all accepted plans, the rest of the guidelines are only met by 30-60\% of the plans. Similarly,  50\% of all accepted patients met the clinical guideline of having at least 40.28~Gy dose to one of the cavity structures, highlighting that the clinical guidelines are not strictly followed in practice. 
We note that this dataset did not include any plans that met all the criteria but were rejected based on the current guidelines, \changemarker{perhaps an indication that these particular guidelines are collectively too tight}. Similarly, there were no accepted plans that met all guideline criteria. 
In our analyses in Section \ref{sec:results}, we use the criteria listed in Table~\ref{tab:Guideline} as clinical guidelines.  

\begin{figure}[htbp]
    \centering
    \includegraphics[width=0.9\linewidth, trim={1cm 0.5cm 2cm 0.5cm},clip]{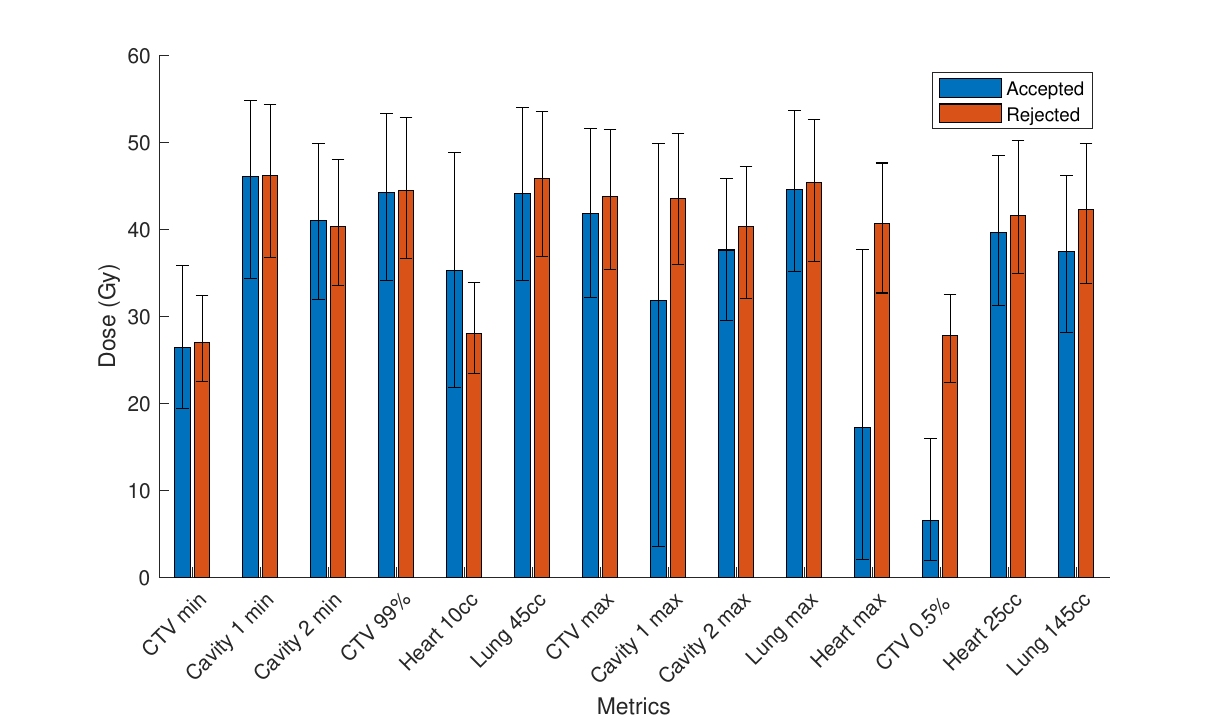}
    \caption{Visualization of a set of metrics calculated for all accepted and rejected plans. The errorbars show the range of each calculated metric across all plans in each category. %\kg{We need to update this figure (remove the two additional guidelines). Also add the max/min limit of the guidelines to the figure. }
    % \hm{update the labels to remove ``dv'' from the titles, so it matches figure 8}
    }
    \label{fig:Guideline_performance}
\end{figure}

%Here say we extracted more features on each of the past accepted and rejected plans, and looked at 14 metrics which includes those used in the clinical guidelines. 
In addition to the features used in these clinical guidelines, we consider six other features in our models that may carry implicit clinical importance in final treatment plans, bringing the total to 14 clinical features. A summary of these features and their values for accepted and rejected patients is provided in Figure~\ref{fig:Guideline_performance} where the error bars show the range of each metric for all patients in that category. 
%Note that four of the calculated metrics match the guidelines (CTV 99\%, CTV 0.5\%, Heart 10cc, and Heart 25cc). 
As the figure illustrates, there is no clear separation between accepted and rejected plans and the underlying logic behind the acceptance/rejection decision cannot be inferred by just considering these metrics. Specific patient examples from this dataset were also previously shown in Figure~\ref{fig:patientexamples}.

%\kg{(Fig 9 with the addition of the guidelines or a table with 100\% data, what percentage meets or does not meet the current guidelines for each metric.)}. 

%\hm{There are 8 clinical guidelines based on which a plan is accepted or rejected currently. These guidelines are listed in the table \ref{tab:Guideline_performance}. The prescribed dose is 42.4~Gy. The last column shows what percentage of accepted plans did not meet each of the clinical criteria. We note that this dataset did not include any plans that met all the criteria but were rejected based on the current guidelines.}
\begin{table}[]
    \centering
    \begin{tabular}{l c c c }
    \toprule
     Dose Metric  \hspace{2cm} &   Clinical Limit & \hspace{1cm} \% met in past accepted plans data \\ \midrule
      % CTV min  &  &  & \\ %%%% We ex excluded these from our guidelines on 6/28/2023 because they are too obvious/loose
       Cavity1 min     & $\geq 40.28$~Gy & 50\%\\ %modCavitiy_min
       Cavity2 min     & $\geq 40.28$~Gy & 30\%\\ %DEVcavity_min
       CTV 99\% min  & $\geq 39.01$~Gy &  39\%\\ %varCTV_min
       Heart 10cc max & $\leq 38.16$~Gy & 100\%\\%\\ %varHEART1
       Lung 45 cc max & $\leq 38.16$~Gy & 38\%\\ %varLung1
     %  CTV max  &  &  & \\
     %  Cavity 1 max  &  &  & \\
     %  Cavity 2 max  &  &  & \\
     %  Lung max  &  &  & \\
     %  Heart max  &  &  & \\
       CTV 0.5\% max   & $\leq 45.79$~Gy  & 57\%\\ %varCTV_max
       Heart 25 cc max & $\leq 21.20$~Gy  & 100\%\\ %varHEART2 
       Lung 25 cc max  & $\leq 36.04$~Gy  & 40\%\\ %varLUNG2
       \bottomrule
    \end{tabular}
    \caption{The percentage of accepted plans that meet each of the criteria in the current guidelines.% There is wide variability in how rigorously the guidelines are imposed on accepted plans. While the Heart max dose guidelines are met by all accepted plans, the rest of the guidelines are only met by 30-60\% of the plans. 
    }
    \label{tab:Guideline}
\end{table}

\subsection{Forward and Inverse Problem Description} \label{sec:clinicalFOIO}

We consider a forward optimization problem that imposes a set of clinical dosimetric criteria on the radiotherapy treatment plan for each patient. The objective of the forward problem is minimizing the maximum dose to the healthy tissue, and there is a set of known clinical criteria and unknown implicit constraints that affect the acceptability of a plan, as described in model~\eqref{model:RT_FO} in Appendix~\ref{appendix:FO}. 
%that inputs a set of treatment plans and determines which ones are feasible or infeasible based on a set of constraints that impose limitations on a given clinical metrics. The FO finds the optimal plan among the feasible set, according to a known objective function. In this application, we consider the objective of minimizing the overall dose to the entire body to reduce any unnecessary radiation exposure. %The details of the FO model are provided in Appendix~\ref{appendix:FO}. 
%
The inverse model then inputs a set of past patient treatments which are labeled as accepted or rejected and infers a set of constraints that capture the underlying implicit constraints.  %We consider that all constraints are linear, which is relevant in practical settings for radiation therapy since dose-volume criteria for breast cancer can be re-written as linear constraints~\citep{chan2014robust}. 
%We assume that there are no known constraints in FO, and for simplicity of computations, 
%that correctly classify the accepted and rejected plans as feasible and infeasible, respectively, and makes a preferred treatment plan optimal for the corresponding FO problem. 
%The inverse model then inputs a set of past patient treatments that are labeled as accepted or rejected and infers a set of constraints that captured 
%We assume that there are no known constraints in FO, and for simplicity of computations, consider that all constraints are linear, which is relevant in practical settings for radiation therapy since dose-volume criteria for breast cancer can be re-written as linear constraints~\citep{chan2014robust}. 

The inverse problem learns a set of implicit constraints to characterize the underlying feasible region of the oncologists' forward problem based on their historical accept/reject decisions. To do so, we employ the data described in Section \ref{sec:data} to infer 10 constraints using model~\eqref{model:RDIO_SeparationMetric} with the Separation Metric to maximize the distance between the inferred feasible region and the infeasible observations. 
For this application, we considered linear constraints for the FO problem, which is relevant in practical settings for radiation therapy because clinical dose-volume criteria for breast cancer can be written as linear constraints using Conditional-Value-at-risk (CVaR) metrics~\citep{chan2014robust}. The details of the \changemarker{RGIO} model and the parameters used in this application are provided in Appendix~\ref{appendix:FO}. We use these models to obtain results in Section~\ref{sec:results}.

% In what follows, we provide more details around the data and clinical guidelines in Section~\ref{sec:data}, and discuss the corresponding details of the FO and IO models and numerical results in Section~\ref{sec:results}.  %\hm{Maybe move this paragraph to beginning of Sec 5? }

%%%%%%%%%%%%%%%%%%%%%%%%%%%%%%%%%%%%%%%%%%%%%%%%%%%%%%%%%%%%%%%%%%%%%%%%
% \hm{We changed this to a section because the first figure/table in results is very important!}
\section{Results}\label{sec:results}
%%%%%%%%%%%%%%%%%%%%%%%%%%%%%%%%%%%%%%%%%%%%%%%%%%%%%%%%%%%%%%%%%%%%%%%%

%we randomly divided the data into training and testing sets. 
%We used the training set, including both accepted and rejected plans, as input to the inverse model and inferred 10 constraints for the forward problem. We then used the testing set to validate the results and to check whether the inferred constraints can accurately predict which plans in the testing set are accepted or rejected.  
%We validated 
We use a random 60\% of data detailed in Section~\ref{sec:data}, a mix of accepted and rejected plans, as input to the inverse model~\eqref{model:IO_RT} to infer the proposed feasible region. We consider the remaining 40\% of the data as future patients to test the IO methodology and validate the agreement of the inferred \changemarker{RGIO} feasible region with the historical accept/reject decision of clinicians. To benchmark our inferred feasible region against the existing guidelines, we also test these guidelines against the entire dataset of accepted and rejected plans. 
% compares to the current clinical guidelines. For these future patients, we consider those that were accepted and rejected by clinicians, and 
The goal is to test whether using the inferred \changemarker{RGIO} constraints and the clinical guidelines would render the accepted and rejected observations as feasible and infeasible for FO, respectively. Vice-versa, we also test whether the accepted and rejected plans align with the inferred constraints or the clinical guidelines. 
%they are also feasible for the feasible regions generated by the guidelines and the inferred constraints using \changemarker{RGIO}. 
Tables~\ref{tab:vennData} and~\ref{tab:vennData2} summarize these results.

\begin{table}[htbp]
    \centering
    \begin{tabular}{l c c }
    \toprule
   % \multirow{2}{*}{Future patient plans} & If meet all & If meet all \\ 
   % &  guidelines & DDIO constraints
   \textbf{Future patient plans} & \hspace{1cm} \textbf{\% Accepted} ($\mu \pm \sigma$)  & \hspace{1cm} \textbf{\% Rejected ($\mu \pm \sigma$)} \\ \midrule
   If plan is feasible for guidelines & N/A* & N/A*\\ %\midrule
   If plan is feasible for \changemarker{RGIO} & $98.4 \pm 0.03$ & $1.6 \pm 0.03$ \\ \bottomrule  %
    %& &  \\ \midrule  % **** from "CompareVenn" matlab code (line 161), "precision". **** 
    % \multirow{3}{*}{Future Plans}& \multicolumn{4}{c}
  %  {\% Accepted}  & N/A* & $(98.2 \pm 0.03) \%$   \\ \midrule    %TP
  %  {\% Rejected}  & N/A* & $(1.8 \pm 0.03)$ \%   \\  \bottomrule       %FN
\end{tabular}
    \caption{Comparison of quality of \changemarker{RGIO} vs guideline constraints when a plan meets all constraints. \newline
    % First and second columns show the percentage of accepted/rejected future patient plans when the plan meets all the guidelines and DDIO constraints, respectively.
    %Percentage of constraints met for future patients considering either the inferred constraints or the guidelines, as decided by an oncologist. 
    {\normalfont *None of the future patient plans meet all the guidelines.    } }  %precision 
    \label{tab:vennData}
\end{table}   

First, in Table~\ref{tab:vennData}, we separately consider plans that are rendered feasible/infeasible by the guidelines and then by the inferred \changemarker{RGIO} constraints. The two columns show the percentage of the feasible/infeasible plans according to each method that was historically accepted/rejected by clinicians, respectively. Upon observing the data, we noticed that none of the plans in our future patient dataset met all the guidelines, hence, we cannot calculate the required metric on an empty set. For the inferred \changemarker{RGIO} constraints, on the contrary, an average of 98\% of the plans that were feasible according to the inferred constraints were also accepted by the clinicians, with a small standard deviation of 0.03\%. 
%These future patients into those that are feasible plans according to the guidelines and/or the inferred constraints.
%\kg{3. write the rest of the text for Tables 3 \& 4. }

Next, in Table~\ref{tab:vennData2}, we perform the opposite analysis. Consider all the future patient plans that are accepted or rejected by clinicians. Table~\ref{tab:vennData2} shows the percentage of accepted plans that are also feasible with respect to the clinical guidelines and the inferred \changemarker{RGIO} constraints. In line with what we observed in Table~\ref{tab:vennData}, none of the accepted plans were deemed feasible by the guidelines, while an average of 94.2\% of the accepted plans were also deemed feasible by the inferred \changemarker{RGIO} constraints. 
\begin{table}[htbp]
    \centering
    \begin{tabular}{l c c c c}
    \toprule
    {\textbf{Future patient plans}}  & \hspace{.5cm} \textbf{\% Feasible for guidelines} & \hspace{.5cm} \textbf{\% Feasible for RDIO} \\
    %Future &  \% feasible for & \% feasible for   \\
    %patient plans & guideline & DDIO \\
    \midrule
    % \multirow{3}{*}{Future Plans}& \multicolumn{4}{c}
%    \% does not meet at least one guideline & 100\% &  100\%  \\ \midrule    
%    \% does not meet at least one DDIO constraint & $(94.2 \pm 0.14) \%$ & ($5.8 \pm 0.14)$ \%   \\  \bottomrule       
    {Accepted}  & 0 \% *& $(94.2 \pm 0.14) \%$   \\ \midrule   
    {Rejected}  & 100 \% *& ($5.8 \pm 0.14)$ \%   \\  \bottomrule       
\end{tabular}
\caption{Comparison of quality of \changemarker{RGIO} vs. guideline constraints when a plan is infeasible.   \newline
{\normalfont *All future patient plans were infeasible according to guidelines.} } %TN/{(FP+TN)}
\label{tab:vennData2}
\end{table}

The insights highlighted by Tables~\ref{tab:vennData} and~\ref{tab:vennData2} confirm that using the clinical guidelines to generate plans is not helpful in practice because guidelines are too restrictive and impossible to meet in this dataset, whereas plans that are made based on \changemarker{RGIO} inferred constraints are highly likely to be accepted by clinicians. Similarly, almost all clinically accepted plans meet all the \changemarker{RGIO} inferred constraints, which is an indication that learning the acceptability criteria through \changemarker{RGIO} can help in building plans that are acceptable to clinicians. 

\begin{figure}[htbp]
    \centering
   \includegraphics[width=.65\linewidth]{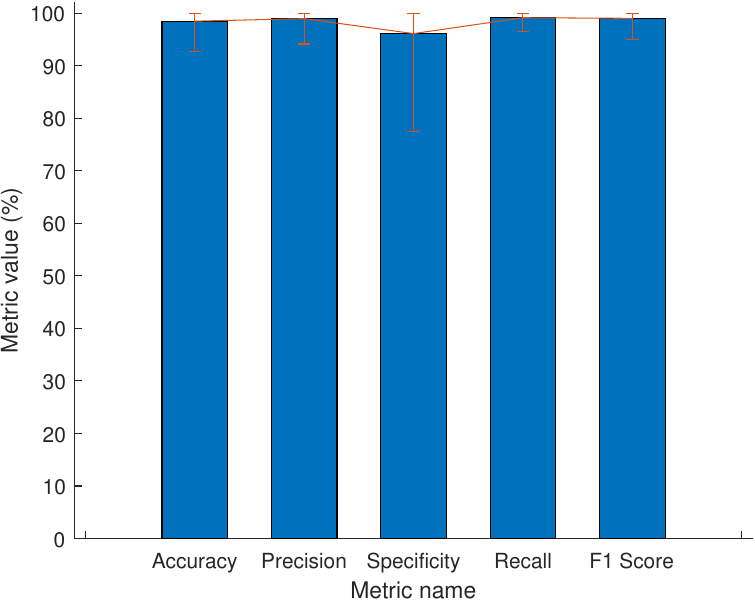}
   \caption{Out-of-sample performance of the inverse model using a 60/40 split for past/future patients. } %\hm{need to update these figures. And add guidelines for the 100\% testing set.}} %this figure is generated in "drawResults.m". Latest results from March 2023
   \label{fig:accuracy_metrics}
\end{figure}

% \hm{The rest of this needs to be updated based on what we added before here: Maybe even a new section. }
To further test the accuracy of the proposed inverse methodology, %we randomly divided the data into training and testing sets. 
%We used the training set, including both accepted and rejected plans, as input to the inverse model and inferred 10 constraints for the forward problem. We then used the testing set to validate the results and to check whether the inferred constraints can accurately predict which plans in the testing set are accepted or rejected.  
%We validated 
we calculated a set of standard metrics for out-of-sample prediction, including accuracy~(\% correct prediction), precision~(\% correct acceptable prediction), specificity~(\% correct rejection identification), recall (\% correct acceptable identification), F1 Score (harmonic average of recall and precision).
% \paragraph{Results: Out-of-sample Prediction Accuracy}
We tested the model using a similar 60/40 random split for past/future patients, repeated the split 50 times, ran the result for each split, and computed the average and the range for each metric. Figure~\ref{fig:accuracy_metrics} shows that the inverse model consistently performs well with an average performance between 95\%--100\% across all metrics. The lowest metric was specificity at 95\%, which we believe is due to a low number of rejected observations in our dataset.  

\begin{figure}[htbp]
    \centering
    \includegraphics[width=0.48\linewidth, trim={0 0 0 0}, clip=true]{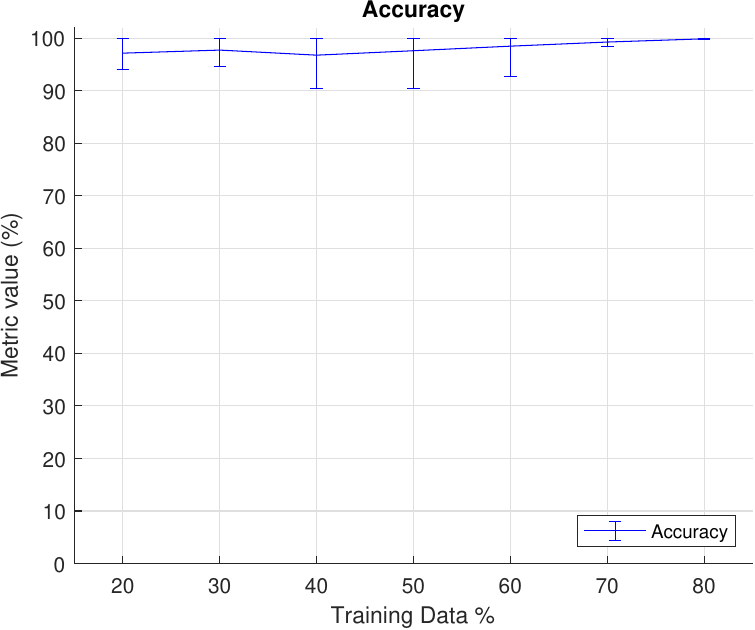} \hfill
    \includegraphics[width=0.48\linewidth, trim={0 0 0 0}, clip=true]{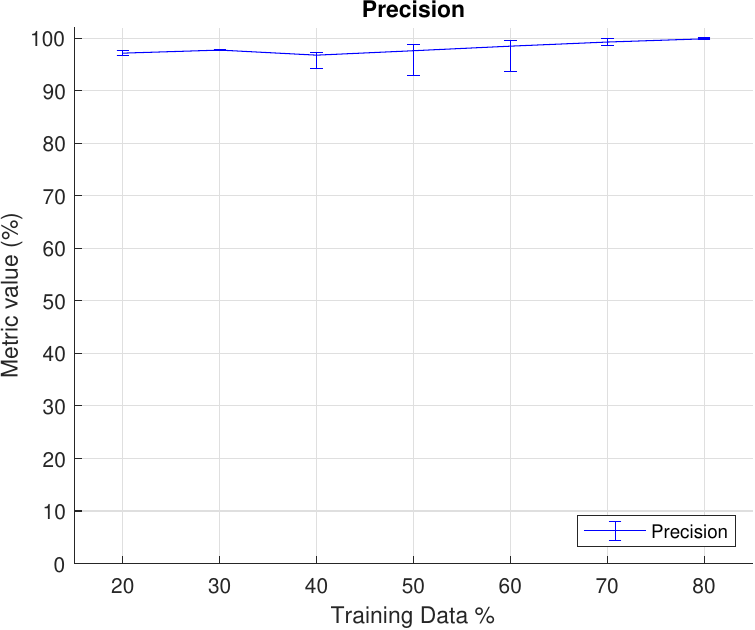}\\ [10pt]
    \includegraphics[width=0.48\linewidth, trim={0 0 0 0}, clip=true]{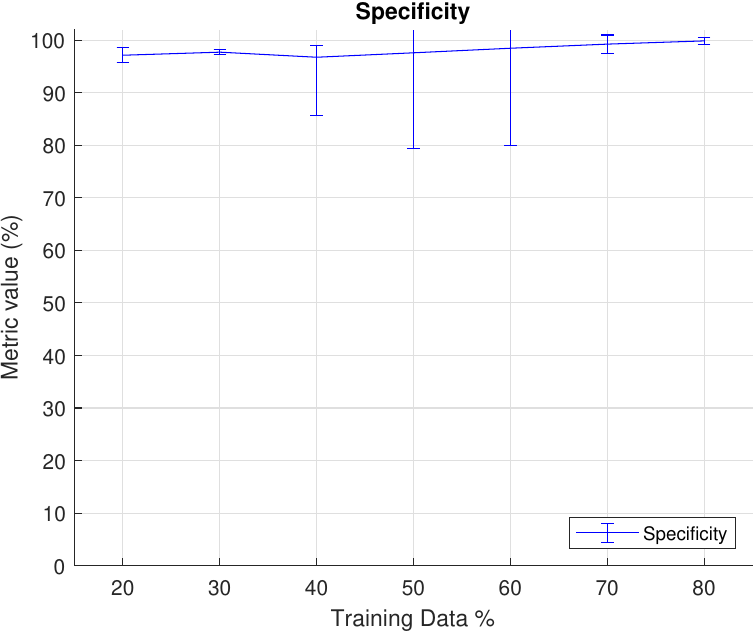} \hfill
    \includegraphics[width=0.48\linewidth, trim={0 0 0 0}, clip=true]{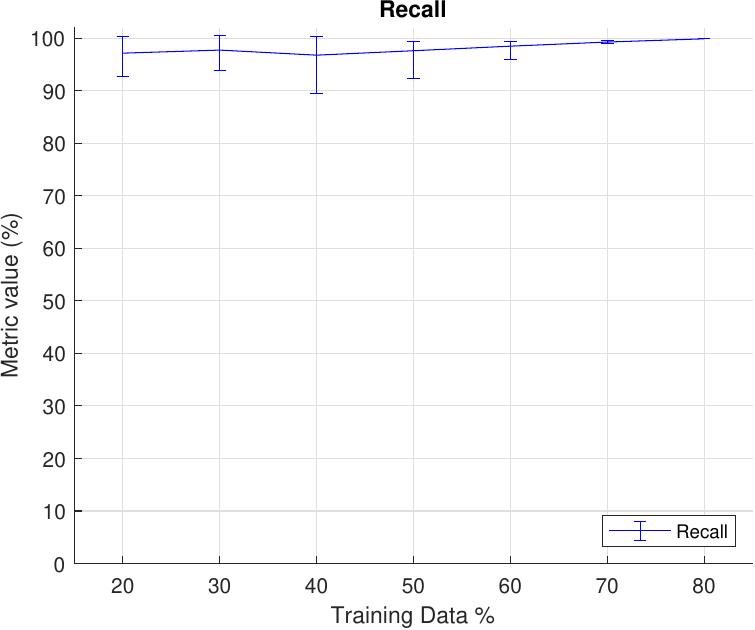}\\ [10pt]
    \includegraphics[width=0.48\linewidth, trim={0 0 0 0}, clip=true]{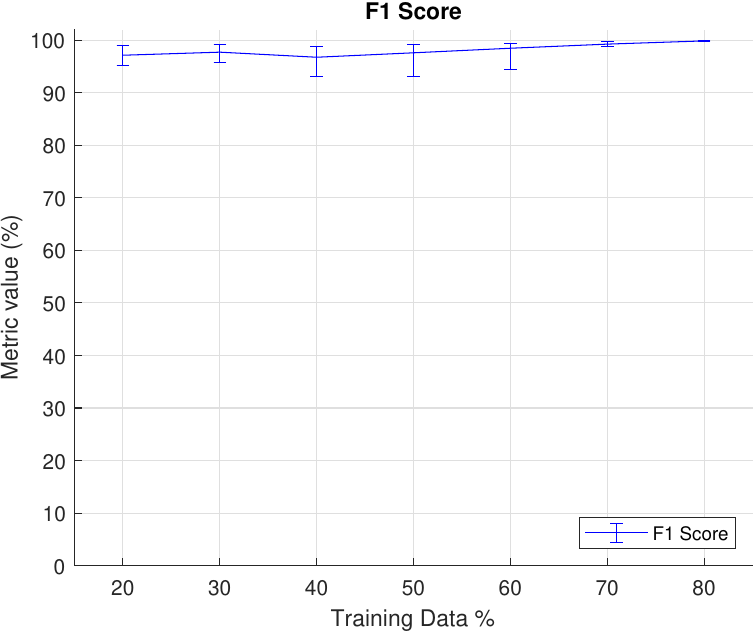}
    \caption{Sensitivity of different metrics with respect to the proportion of past/future patients} %\hm{update these figures to the newest results.}}
    \label{fig:accuracy_trends}
\end{figure}

We next tested the sensitivity of the proposed method to the required number of past patients used to train the \changemarker{RGIO} model. %The data was split randomly into training and testing in increments of 10\% starting from a 20-80 split up to an 80-20 split. Again, we ran the results 50 times and computed the average of each metric. 
%In each case, we inferred
Figure~\ref{fig:accuracy_trends} illustrates the trends of each metric when the number of past patients is varied between 20\% to 80\% of all patients in increments of 10\%, with 250 random splits for each percentage. 
The x-axis shows the \% of the data that was used as past patients (the rest used as future patients), and the y-axis shows how each metric performed (on average). %
The results show that the method works well across the board, even when using only 20\% of the available data as past patients. 
We see a slight dip right before the 50-50 split but it picks up again when we pass the ratio. While there is always randomness, the methods seem to consistently obtain good results on average (with all metrics above 95\%) at a 60-40 split. %Note that if we train too much, the risk of overfitting increases and the out-of-sample result may suffer, which for this application seems to happen at around 70-80\% and upward. 
Even in the worst case, our recall and F1 scores are above 90\%.
%   \begin{itemize}
%             \item Training set size from 20\% to 80\% (250 random splits for each)
%     \end{itemize}

% \begin{itemize}
%     \item Overall performance on all metrics $\geq 80\%$
%     \item Suggested training set size: 60\% with all metrics $\geq 95\%$
%     \item All metrics $\geq 90\%$ when training size $\geq 60\%$
% \end{itemize}

% {\red separate \ref{fig:accuracy_trends} into several figures, with trends over errorbars or boxplots}. 
           
%\subsection{Discussions and limitations}   %\hm{to edit}

% Report on how current guidelines compare with imputed guidelines. Which ones are and are not always met with our labeled constraints? Which metrics are correlated?  what are the imputed guidelines? is it explainable? can we find an example of an explainable constraint that was inferred through IO \hm{we don't have this}

%Double-check accepted/rejected vs feas/infeas.

%---------------------------
%---------------------------
%---------------------------
%---------------------------

\section{Limitations and Future Work} \label{sec:future} %future work
% ALT: Further Research Opportunities
\changemarker{In this section, we review a few limitations of the presented work in this paper and provide directions for future research to address these limitations. %In what follows, we first 
We first discuss the need for data pre-processing (Section~\ref{sec:preprocessing}), followed by methods for handling ill-posed data and allowing misclassification (Section~\ref{sec:misclassification}). We then discuss possible routes for mitigating data uncertainty (Section~\ref{sec:robust}) and overfitting (Section~\ref{sec:overfitting}), and finally highlight the key differences between the proposed approach with well-known classification approaches and motivate further benchmarking (Section~\ref{sec:benchmarking})}. 
%\kg{We should expand and add more of the "limitations" of the work and what else can be done with the data. }
% preprocessing, noise, uncertainty, benchmark against classification, solution method, extending to good/bad for objective?. 
% \hm{include headers: preprocessing, handling misclassification, noisy data, data uncertainty and overfitting, Further Benchmarking}

\subsection{Data Preprocessing} \label{sec:preprocessing}
% pre-processing can speed things up:
 \changemarker{In our framework,} %Given the data-driven nature of the framework, 
data quality and pre-processing can play an important role in the size of our models and the quality of the obtained solutions. Simple pre-processing methods can potentially reduce the size of the observations by identifying and removing those observations that are already rendered infeasible by the known constraints. A reduction in the number of these infeasible observations can largely improve the complexity and size of the inverse problem. %by reducing the required number of inferred constraints. 
Any redundant observations or outliers can be removed from the set of observations using a data-cleaning method of choice.
% Any redundant accepted or rejected observations (e.g., observation points with a distance less than a given threshold from each other) or outliers can be removed from the set of observations using the data cleaning method of choice.

Data pre-processing can also assist in ensuring acceptable data quality, removing conflicts between the data and known constraints, and increasing adherence to the properties required by the proposed models. Preprocessing methods can analyze past decisions, identify conflicting observations and criteria, and prompt the decision-maker to either remove or rectify any mismatch in the past data. 
This pre-processing results in data that is well-defined for the FO problem. \changemarker{If the data is still ill-posed, we will need to extend the model to allow inconsistencies in the data, as discussed next.}
%% We assume that the data is well-defined for the FO problem, and those data points that contradict the FO properties can be removed from the set. For instance, if a rejected observation exists within the convex hull of the accepted observations, which contradicts the convexity assumption of the FO model, it needs to be removed. 

\subsection{Ill-posed Data}\label{sec:misclassification}
%% preprocessing ensures data is well-defined
\changemarker{In Section~\ref{sec:method}, we assumed that the data is well-posed and there are no inconsistencies in the observations, particularly, $\{\not \exists k \in \cKbad | \bx^k \in \cH\}$, indicating that there are no rejected observations in the convex hull of accepted observations as captured by Assumption \ref{assumption:wellposed}. In practical settings, however, this assumption may be violated due to data inconsistency. Data pre-processing can capture and correct some of such inconsistencies, as discussed in Section~\ref{sec:preprocessing}.   %A regularization term can also be used to complement, or replace, the data pre-processing step by capturing and minimizing inconsistency error. 
%In particular, we assume that $\not \exists \cKbad \in \cK$ such that $\bx^k \in \cH$, which indicates that there are no rejected observations in the convex hull of accepted observations. 
% In practice, this can be achieved by pre-processing the data and removing any such observations that are inconsistent. 
Another approach is to modify the RGIO model to introduce a regularization term that penalizes a slack that captures the violation in the inferred constraints to the loss function of the inverse problem. This addition will change the nature of the inferred constraints from hard constraints to soft constraints that allow violations with a certain penalty. 
% \hm{More explicitly say what can be done if there is a bad obs in the middle of good ones. Like SVM, minimize misclassification error, and introduce a measure function. Talk about hard and soft margins and hinge loss function.}
In the context of our proposed RGIO model, this can be achieved via the following model {RGIO-$\xi$}. 
\begin{subequations}
\begin{align} 
%hspace{-1cm}
% \mathbf{RGIO-\xi}: 
\quad \underset{\ba, b, \bq, \by}{\text{Maximize}} \quad  &   
%\underset{k \in \cKbad}{\max}(\bP ([\bA, \bb]_\text{inferred},\bx^{k})) \label{eq:eIOobj} \\ 
%\bP \left([\bA, \bb]_\text{inferred},\left\{\bx^{k},\, k\in \cK\right\}\right) \label{eq:eIOobj} \\ 
%\bP \left(\cXt,\left\{\bx^{k},\, k\in \cK\right\}\right) 
\bP \left(\bq_1, \dots, \bq_{|\cNt|}, \bA, \bb; (\bx^{1},\dots \bx^{|\cK|}) \right) - \bP_\xi(\xi^{\cL+},\xi^{\cL-},\xi^{\cN+},\xi^{\cN-})\label{eq:eIOobjXi} \\
 \text{subject to} \quad & g_n(\bx^k; \,\bq_n) \ge -\xi_{n,k}^{N+}, \quad \forall k\in \cKgood , n\in \cNt 
 \label{eq:eIOprimalfeasNLXi} \\ 
&   \ba_\ell'\, \bx^k - b_\ell\ge -\xi_{\ell,k}^{L+},   \qquad \forall k\in \cKgood, \ell \in \cLt \label{eq:eIOprimalfeasLXi}  \\
& g_n(\bx^k; \,\bq_n) + \epsilon - M y_{nk} \le \xi_{n,k}^{N-}, \quad \forall n \in \cN \cup \cNt, \, k \in \cKbad \label{eq:eIOinfeasNLXi}\\
& \ba_{\ell}\bx^k - b_{\ell} + \epsilon - M y_{\ell k} \le \xi_{\ell,k}^{L-}, \quad \forall \ell \in \cL \cup \cLt, \,  k \in \cKbad \label{eq:eIOinfeasLXi}\\
& \xi^{\cL+},\xi^{\cL-},\xi^{\cN+},\xi^{\cN-} \ge \bzero,\\
& \eqref{eq:eIOinfeasSum}-\eqref{eq:eIOBinaryVars}. \notag
\end{align}
\end{subequations}
}

\subsection{Data Uncertainty} \label{sec:robust}
%Noise, uncertainty, overfitting:
%Similar to other inverse optimization methods, existing noise and uncertainty in the observations can impact the results. Unlike existing noise mitigation methods, noise in constraint inference  
%
% \hm{Robustness against data uncertainty} 
\changemarker{An area of future research is to consider the effect of data uncertainty on constraint inference and to explore robustness conditions for the inverse framework. In the literature on inverse for objective inference, robust optimization has been used to tackle data uncertainty where an uncertainty set is considered around each data input \citep{ghobadi2018robust}.  A similar approach can be used for constraint inference, where uncertainty sets are considered around both the accepted and rejected observations and the goal of the inverse model is to infer a feasible region that is robust against the worst-case scenario of data uncertainty. We note that these uncertainty sets may overlap, which may necessitate the use of soft margins and a notion of error as discussed in Section~\ref{sec:misclassification}. 
}

\subsection{Overfitting} \label{sec:overfitting}
\changemarker{Depending on the loss function used in our proposed inverse framework, it is possible to overfit to the data by either inferring a feasible region that is tightly wrapped around the accepted observation and disallowing other potentially acceptable solutions that are just outside this region, or on the contrary, inferring a loose feasible region that only eliminates past rejected solutions. An intuitive method to avoid overfitting is to consider a margin around observations, which would be mathematically similar to considering an uncertainty set around each observation, as discussed in Section~\ref{sec:robust}. Another possible direction is to define a measure of optimism (or pessimism) which indicates how much to expand the feasible region to allow better solution (or contract the region to deviate from the rejected data). In combination with the methods from Section~\ref{sec:misclassification}, this tunable parameter would allow the user to decide how conservative they want to be concerning overfitting to both the accepted and the rejected past data.} 

\subsection{Benchmarking}\label{sec:benchmarking}
%\hm{lead to the next section}
% \kg{try to improve this paragraph.}\\
\changemarker{We note that although this work conveys some similarities with classification methods, it is structurally different from the conventional artificial intelligence classification, e.g., support vector machines. Inferring a feasible region has connections with binary classification methods, distinguishing between the interior and exterior of the feasible region. While traditional classification methods often include assumptions on the type of separators, they rarely consider any assumptions on the recovered region or objective functions. %In this work, we focus on two categories of accepted and rejected decisions with the explicit assumption that the region containing accepted decisions is a convex set. 
In this work, we include the explicit setting where the inferred feasible region must be a convex set so that it can later be used in a forward convex optimization setting. This connection can be further explored between existing classification methods and inverse optimization, investigating distance measures that can be used and the properties of the obtained region. 
Future work can also explore the relationship between the distribution of observations and the inferred feasible region. }

% \hm{Our solution method can be similar to classification (i.e., minimizing distance to feas/infeas but is not limited to that. Future work can explore the relationship between the distribution of observations and the inferred feasible region. Beyond the scope of the current paper.}
% 

% \kg{This work of good and bad could be extended to inferring objectives and making the obj be closer to good solutions and far from bad solutions. }
% \kg{probably delete the following}

% \iffalse
% Similarly, if the preferred observation is an interior to the convex hull of the acceptable observation, then it is conceivable that either the set of observations is too limited or a measurement error exists because there exists a point in the convex hull that has a better objective value than the current preferred solution, and hence, should be optimal for the inferred FO. This new point can be identified using heuristic algorithms and added to the set of observations to mitigate the situation. 
% \fi

\section{Conclusions} \label{sec:conclusions}

In summary, this paper provides an inverse optimization framework that inputs both accepted and rejected observations and imputes a convex feasible region of a forward problem. %We consider the case where a large number of observations are available. The goal is to find a feasible set for the forward problem such that all given observations become feasible and the preferred solution(s) become optimal. 
%We demonstrate 
%{theoretical properties} of our methodology and propose a new tractable reformulation. %We also present and discuss several measure functions that can be used to derive imputed sets that have certain desired properties. 
\changemarker{The proposed inverse model is a complex nonlinear mixed integer formulation; therefore, using the properties of the constructed solutions, we propose a reduced reformulation that partially mitigates the problem complexity by using variational inequalities to render a set of nonlinear optimality conditions 
redundant.} %\kg{changed wordi ng} 
%Our numerical case studies demonstrate the differences between various modeling approaches. 

% these measure functions and serve as a basic guideline for users to choose the appropriate measure function to find imputed sets with desired properties, depending on the available data and the relevant application.
We consider the problem of {radiation therapy treatment planning} to {infer the implicit clinical guidelines} that are utilized in practice based on historical plans. Using realistic patient datasets, we randomly divide the data into past and future patients and use our inverse models to derive a feasible region for accepted plans such that rejected plans are infeasible. Our results show that on average, 98\% of the plans that are feasible for our inferred models are also historically accepted, and our obtained feasible region (i.e., the underlying clinical guidelines) performs with 95\% accuracy in predicting the accepted and rejected treatment plans for future patients, on average. 
The results also highlight an interesting property of IO models: even with a small size of past patient data, a high accuracy can be obtained in generating acceptable plans for future patients. This characteristic, perhaps one of the advantages of IO compared to conventional machine learning methods, is due to the fact that IO models have a predetermined optimization structure. 
%confirm that IO works well even with a small size of training data, which is perhaps one of advantages of IO compared to conventional machine learning methods, due to the fact that IO models have a predetermined optimization structure. 

In radiation therapy, implementing our methodology results in standardized clinical guidelines and infers the true underlying criteria based on which the accept/reject decisions are made. This information allows planners to generate higher-quality initial plans that, in turn, result in better patient care through personalized treatment plans. It also streamlines the planning process by reducing the number of times a plan is rejected and sent back for corrections. The standardization of the guidelines reduces variability and potential human errors among clinicians and health care centers, and enables personalized guidelines for patient subpopulations.  
%
%other applications
We believe the method can also be applied in other application areas where an understanding of implicit expert constraints can help streamline the decision-making processes.
%solution methods
One future direction consists of focusing on special classes of convex problems and devising more efficient solution methods for large-case data-driven applications.

\bibliographystyle{informs2014}
\bibliography{robustinverse2}

\clearpage

% \section*{Appendices}
\begin{APPENDICES} 

\section{Proofs}

% \begin{APPENDIX}{Proofs}
 % \DoubleSpacedXI
% \section{Proofs}
\begin{proof}{\textsc{Proof of Proposition~\ref{prop:DIOoptimal}.}} 
%(it's not unbounded). 
Constraints~\eqref{eq:IOprimalfeasNL}--\eqref{eq:IOprimalfeasL} ensure $\bx^k \in \cXt, \, \forall k \in \cK$,  and  %, which correspond to the feasibility requirement~\eqref{eq:nominalGood} in Definition~\ref{def:nominal}. 
constraints \eqref{eq:IOinfeasL}--\eqref{eq:IOinfeasSum} use binary variables $y_{ik}$ to ensure that at least one constraint (either inferred or known; linear or nonlinear) makes each rejected point infeasible, which, in turn, ensures $\bx^k \not \in \cX \cap \cXt$. Therefore, the conditions of Definition~\ref{def:nominal} are met and $\cXt$ is a nominal set. 
Since the set $\cXt \cap \cX$ is convex by definition, constraints~\eqref{eq:IOprimalfeasNL}--\eqref{eq:IOcsL} provide the necessary and sufficient KKT conditions to guarantee $\bx^0 \in \underset{\bx \in\cX \cap \cXt}{\argmin}\,\{f(\bx; \bc) \}$.  \hfill \halmos
%To show that $x^k, k\in \cKbad$ will be infeasible to this model, we note that based on enforced constraints it cuts out infeasible solutions too. To do so, we use auxiliary binary variables $\by$ to activate the opposite direction of each constraint and ensure that at least one of the constraints is violated among all the known and inferred linear or nonlinear constraints. 
% KKT: we need to ensure the KKT conditions for $\bx^0$ holds including stationarity conditions and complementary slackness as illustrated in constraints \eqref{eq:IOstationarity}--\eqref{eq:IOcsL} and \eqref{eq:IOstationaritySign}. We note that since $\bx^0$ is a feasible observation, the primal feasibility condition of KKT is satisfied through Req 1. 
% $g_n$ is convex based on Req 4. 
\end{proof}
\vspace{2em}

\begin{proof}{\textsc{Proof of Proposition~\ref{prop:DIOfeas}}.}
To show that the feasible region is non-empty, %it is sufficient to find one solution that satisfies all constraints. 
we construct a solution (by assigning values to all variables $\bq, \ba, b, \lambda, \mu,$ and $y$)  %and then 
that satisfies all constraints of \changemarker{GIO} and is, hence, feasible.  
Let $\bq_n = \hat{\bq}, \forall n \in \cNt $, as defined in Assumption~\ref{assumption:wellposed}. 
By Assumption~\ref{assumption:wellposed} and the Separating Hyperplane Theorem \citep{boyd2004convex}, for every rejected point $k \in \cKbad,  \exists\, \hat{\ba}_k, \hat{b}_k$ such that $\hat{\ba}_k' \bx^k < \hat{b}_k$ and $\hat{\ba}_k' \bx^p \geq \hat{b}_k, \, \forall p \in \cKgood$. Let $[\ba_{\ell}] = [\hat{\ba}_1, \dots, \hat{\ba}_{|\cKbad|}, \nabla f(\bx^0; \bc), \dots, \nabla f(\bx^0; \bc)]$ and  $\bb=[\hat{b}_1, \dots, \hat{b}_{|\cKbad|},\nabla(\bx^0;\bc)'\bx^0, \dots, \nabla(\bx^0;\bc)'\bx^0], \quad \forall \ell \in \{1,\dots, \cLt\}$. 
Let the dual variables of the nonlinear and linear constraints be $\lambda_n = 0, \forall \in \cN \cup \cNt$, and $[\mu_{\ell}] = [\bzero_{1\times |\cL \cup \cLt| -1},-1], \forall \ell \in \cL$, respectively.  
Finally, let binary variables $y_{kk} = 0, \forall i \in \cLt, k \in \cKbad$ and $y_{kk} = 1, \forall i \in \cI \setminus \cLt, k \in \cKbad$.
By substitution, it can be seen that this constructed solution satisfies all constraints \eqref{eq:IOprimalfeasL}--\eqref{eq:IOBinaryVars}. 
Note that the first $|\cKbad|$ constraints ensure that each rejected point is infeasible for at least one inferred linear constraint, and the last set of constraints (which are all identical) ensure the optimality of $\bx^0$. All constraints satisfy primal feasibility  $\forall \bx^k, \, k \in \cKgood$. 
\hfill \halmos
\end{proof}
\vspace{2em}

\begin{proof}{\textsc{Proof of Proposition~\ref{prop:imputedinC}}.}
Let $\cS = \cX \cap \cXt$ be an imputed set for FO. Assume that $\cS \not \subseteq \cC$. Then $\exists \, \hbx \in \cS \setminus \cC$. We know that $\hbx \not \in \cV$ since it contradicts the definition of the preferred solution $\bx^0$. Hence, $\hbx \in \cV\setminus\cC$. Consider the following two cases:

\noindent {\bf\it Case 1:} If $f(\bx;\bc)$ is linear, then $\cC=\cV$ and $\cV\setminus\cC = \varnothing$, which is in contradiction with $\bx^k \in \cV\setminus\cC$. 

\noindent {\bf \it  Case 2:} If $f(\bx;\bc)$ is nonlinear convex, then $\cV$ is non-convex, and hence $\cV\setminus\cC$ is non-convex. Given that $\cC$ is the tangent half-space of the sublevel set of $f(\bx;\bc)$ at $\bx^0$, then $\bx^0$ is on the boundary of $\cV$ since $f$ is monotone. Because $\cS$ is convex, for any $\hbx \in \cV \setminus \cC$, there must exist $\lambda >0$ such that $\bar{\bx} = \lambda \hbx + (1-\lambda) \bx^0$ and $\bar{\bx} \in \cS \setminus \cV$, which contradicts $\bx^0 \in \underset{\bx \in \cS}{\argmin}\,\{f(\bx; \bc) \}$. 
\hfill \halmos
\end{proof}
\vspace{2em}

\begin{proof}{\textsc{Proof of Proposition~\ref{prop:CcapS}}.}
Both $\cS$ and $\cC$ are convex so $\cC \cap \cS$ is also convex. Because $\cS$ is a nominal set and $\cC$ includes all accepted observations, $\cC\cap\cS$ is also a convex nominal set. Finally, Since $\cS$ is a convex nominal set, it includes the convex hull of all observations, and $\bx^0 \in \cS$, and therefore, $\bx^0\in\cS\cap\cC$. Because $\cC$ is the tangent half-space to the sublevel set of $f(\bx;\bc)$ at $\bx^0$, it is also given that $\bx^0 \in \underset{\bx \in \cS\cap \cC}\argmin\{f(\bx;\bc)\}$. Therefore, $\cS \cap \cC$ meets the criteria outlined in Definition~\ref{def:imputed} and is, therefore, an imputed set. 
%\hm{show by contradiction.}
\end{proof}
\vspace{2em}

\begin{proof}{\textsc{Proof of Theorem~\ref{theorem:RDIO}.}}
(\textit{i}) Any solution $\ba, b, \bq, \lambda, \mu,y$ of \changemarker{GIO} is also a solution of \changemarker{RGIO} because it has fewer constraints. (\textit{ii}) Vice-versa, for any solution $\ba, b, \bq, y$ of \changemarker{RGIO} there exists a solution of \changemarker{GIO} since $\cC$ is appended to the known constraints and the stationarity and complementary slackness conditions can be re-written as: 
\begin{align*}
    & \nabla f(\bx^0; \bc) + \sum_{n\in \cN\cup \cNt} \lambda_n \nabla g_n\, (\bx^0,\, \bq_n) + \lambda_0 \nabla f(\bx^0;\bc)  + \sum_{\ell \in \cL \cup \cLt}\mu_\ell \, \ba_{\ell} \, = \bzero \\ %\label{eq:eIOstationarity}  \\
& \lambda_0 \,\, (\nabla f(\bx^0; \bc)'\bx^0-\nabla f(\bx^0; \bc)'\bx^0) = 0  \\ 
& \lambda_n \,\, g_n(\bx^0, \bq_n) = 0 \qquad \forall n\in \cN\cup \cNt, \\%\label{eq:eIOcsNL}  \\
& \mu_\ell \,\, (b_\ell - \ba_\ell' x^0)= 0, \qquad \forall \ell\in \cL\cup \cLt \\ %\label{eq:eIOcsL}\\
& \lambda_n,\mu_\ell \leq 0, \qquad \forall n\in \cN\cup \cNt, \, \ell\in \cL\cup \cLt.%\label{eq:eIOstationaritySign}
\end{align*}
%Here,
%\[\nabla f(\bx^0; \bc) + \sum_{n\in \cN\cup \cNt} \lambda_n \nabla g_n\, (\bx^0,\, \bq_n) + \lambda_0 \nabla f(\bx^0;\bc)  + \sum_{\ell \in \cL \cup \cLt}\mu_\ell \, \ba_{\ell} \, = \bzero\] is always satisfied when 
All conditions are satisfied if we set $\lambda_0 = -1$ and $\lambda_n,\mu_\ell \leq 0, \, \forall n\in \cN\cup \cNt, \, \ell\in \cL\cup \cLt$. %This also ensures that \eqref{eq:IOcsL}--\eqref{eq:IOcsNL} are always met. 
Hence, for any solution to \changemarker{RGIO}, there exists a corresponding solution to \changemarker{GIO}. Therefore, by (\textit{i}) and (\textit{ii}), the \changemarker{GIO} and \changemarker{RGIO} are equivalent when $\cC$ is added as a known constraint.
\hfill \halmos
\end{proof}
\vspace{2em}

\begin{proof}{\textsc{Proof of Corollary~\ref{cor:theorem}.}}
Due to constraints~\eqref{eq:eIOprimalfeasNL}--\eqref{eq:eIOBinaryVars}, we know that $\cXt$ is a nominal feasible set for FO. 
Note that $\bx^0 \in \cX \cap \cXt \cap \cC$ given that $\bx^0 \in \cX$ by Assumption~\ref{assumption:wellposed}, $\bx^0 \in \cXt$ because $\cXt$ is a nominal set, and $\bx^0 \in \cC$ by definition. Furthermore, $\bx^0 \in \cC = \{\bx \mid f(\bx;\bc) \geq f(\bx^0; \bc) \}$, which implies that $\bx^0 \in \underset{x\in\cC\cap\cX\cap\cXt}{\argmin} f(\bx;\bc)$ and is therefore an optimal solution of~\eqref{eq:theorem}.  
\hfill \halmos
\end{proof}
\vspace{2em}

% \end{APPENDICES}
% \end{APPENDIX}

% \iffalse

\section{Forward and Inverse Models for the Radiation Therapy Application}\label{appendix:FO}
% \begin{APPENDIX}{Forward model for RT}\label{appendix:FO}
%For the forward problem here, 
Assume that a set of features $f \in \cF$ is given for each structure $s \in \cS$, for instance left lung, clinical target volume (CTV), and the heart for a breast cancer patient. Examples of the features in radiation therapy plans for such patients include min, max, mean dose, or dose to a certain volume of each structure. %Let's consider the set $I_1$ to be a vector of all features for all structures, therefore, $|I_1| = |\cF| \times |\cS|$, if all structures have the same number of features. 
The Forward problem we consider is a linear optimization in which the objective function is a linear measure of the features, i.e., $\bx$ which is a vector of $x_i$ where $i \in I_1$ and $I_1$ is the vector of all dosimetric features of the plans. %In this forward optimization problem, some lower and upper bound for each feature may be known as a constraint. We aim to infer a set of unknown constraints on the feature that lead to acceptable plans and are not conventionally considered as a clinical guideline. 
The ordered list of these dosimetric features is provided in Figure~\ref{fig:Guideline_performance}. Accordingly, the forward problem for the RT problem can be written as
%\begin{subequations}
\begin{align}\label{model:RT_FO}
\text{FO: Minimize} \quad & \bc' \bx \\
\text{subject to} \quad & \bG \bx \geq \bh  \qquad \text{(Known constraints)} \notag\\
& \bA \bx \geq \bb,  \qquad \text{(Unknown constraints)} \notag
\end{align}
%\end{subequations}
\noindent where, $\bG=\bI_{14\times14}$ is the identity matrix, $\bh = [10, \bzero_{1\times 13}]'$ to ensure the minimum dose to every voxel in the clinical target volume is always at least 10~Gy and the dose to all other structures are always non-negative. The objective function coefficient vector is $c=[0~0~0~0~0~0~0~0~0~1~1~0~0~0]'$ to minimize the max dose to the organs at risk, namely the heart and the left lung. In this model, we aim to recover 10 unknown underlying constraints using inverse optimization, and hence, the matrix $\bA$ and vector $\bb$ are the unknown parameters of sizes $10\times 14$ and $10 \times 1$, respectively, which will be inferred in the inverse model. 
%%%%%%% ************ Note to self %%%%%%%
%The obj vector c in the AMPL files is wrong (they have index 9 and 12 turned on) but we fix that in the Matlab files (in AnalyzeResults_HM.m and later in GeneralizedData... ) where we calculate the half-space and add it to the results. Our current results are correct
%The order of the structure is the same between matlab, AMPL, and Excel (different from the figures (Fig 9) in the paper)
%%%%%%%%%%%%%%%%%%%%%%%%%%%%%%%%%%%%%%%%%%%%%%%%%%%%%%%
Using the Separation Metric, and based on formulation~\eqref{model:RDIO_SeparationMetric}, the resulting inverse model is a mixed integer linear program and can be written as follows. 
\begin{subequations}
\begin{align}
%\mathbf{RDIO_{1}}: 
\quad \underset{\ba, b, \bp, y,z}{\text{Maximize}} \quad  &  \sum_{k \in \cKbad} z_k \label{RTeq:eIO1obj} \\ 
 \text{subject to}  \quad 
 %& d_{n k} \geq -g_n(\bx^k;\bq_n), \qquad \forall n \in \cNt, k \in \cKbad\\
 %& d_{n k} \leq -g_n(\bx^k;\bq_n) + M y_{nk}, \qquad \forall n \in \cNt, k \in \cKbad \\
%  & d_{\ell k} \geq b_{\ell} - \sum_{m\in M} a_{\ell m} x_{km} \qquad \forall \ell \in \cLt, k \in \cKbad \\
  & d_{\ell k} \geq b_{\ell} - \ba_{\ell}' \bx^{k}, \qquad\quad  \forall \ell \in \cLt, k \in \cKbad \\
%  & d_{\ell k} \leq b_i - \sum_{m \in M} a_{\ell m} x_{km} + M r_{\ell k} \qquad \forall \ell \in \cL, k \in \cKbad \\
%  & d_{\ell k} \leq b_{\ell} - \sum_{m \in M} a_{\ell m} x_{km} + M y_{\ell k} \qquad \forall \ell \in \cLt, k \in \cKbad \\
  & d_{\ell k} \leq b_{\ell} - \ba_{\ell}' \bx^{k} + M y_{\ell k}, \qquad \forall \ell \in \cLt, k \in \cKbad \\ %% this is "bin1" as in the AMPL file
%  & d_{\ell k} \leq M (1-r_{\ell k})\\ %same as the one below
 & d_{\ell k}  \leq M(1-y_{\ell k}), \qquad \forall \ell \in \cLt, k \in \cKbad \\ %% %% this is "bin1" as in the AMPL file
 & d_{\ell k} \geq \epsilon (1-y_{\ell k}),  \qquad \forall \ell \in \cLt, k \in \cKbad \\
 & z_k \leq d_{\ell k} + M p_{\ell k}, \qquad \forall \ell \in \cLt, k \in \cKbad\\
 & \sum_{\ell \in \cLt} p_{\ell k} \leq |\cLt|-1, \qquad \forall k \in \cKbad\\
 %& p_{ik} \leq |\cI|-1 \qquad \forall i \in \cNt \cup \cLt, k \in \cKbad\\
 & p_{\ell k} \geq y_{\ell k}, \qquad \forall \ell \in \cLt, k \in \cKbad\\
 & d_{\ell k} \geq 0, \qquad \forall \ell \in \cLt, k \in \cKbad  
%& g_n(\bx^k; \,\bq_n) \ge \bzero, \quad \forall k\in \cKgood , n\in \cNt 
 \label{RTeq:eIOprimalfeasNL} \\ 
&   \ba_\ell'\, \bx^k \ge b_\ell,   \qquad \forall k\in \cKgood, \ell \in \cLt \label{RTeq:eIOprimalfeasL}  \\
& \ba_{\ell}\bx^k \le b_{\ell} - \epsilon + M y_{\ell k}, \quad \forall \ell \in \cLt, \,  k \in \cKbad \label{RTeq:eIOinfeasL}\\
%& g_n(\bx^k; \,\bq_n) \le \bzero - \epsilon + M y_{nk}, \quad \forall n \in \cN \cup \cNt, \, k \in \cKbad \label{eq:eIOinfeasNL}\\
& \sum_{\ell \in \cLt } y_{\ell k} \leq \, \mid \cLt \mid - 1, \quad \forall k \in \cKbad \label{RTeq:eIOinfeasSum}\\
%& \ba_{\ell} \in \cA_{\ell} \,, b_{\ell} \in \cB_{\ell} \quad \forall \ell \in \cLt \qquad \label{RTeq:eIOnormL} \\
%& \bq_n \in \cQ_n \,, \quad \forall n \in \cNt  \label{eq:eIOnormNL} \\ %||\bq_n||= 1
& \sum_{j=1}^{m} a_{\ell j} = \alpha^+_{\ell} - \alpha^-_{\ell} \label{RTeq:eIOnorm1} \\
& \alpha^+_{\ell} + \alpha^-_{\ell} = 1  \label{RTeq:eIOnorm2} \\
& y_{\ell k},\alpha^-_{\ell},\alpha^-_{\ell} \in \{0,1\}, \qquad \forall \ell \in \cLt, \, k \in \cKbad \label{RTeq:eIOBinaryVars}  \\ 
&\ba_{\ell} \in \mathbb{R}^m, b_{\ell} \in \mathbb{R}, \qquad \forall \ell \in \cLt.  \label{RTeq:domainAandb} 
\end{align}  \label{model:IO_RT} 
\end{subequations}
In this model, to normalize the left-hand-side parameters of the inferred constraints, we use a linear proxy to the $L_1$ norm from \cite{ghobadi2021inferring} as shown in constraints~\eqref{RTeq:eIOnorm1}-\eqref{RTeq:eIOnorm2}.

\end{APPENDICES} 
\end{document}